\documentclass[graybox]{svmult}
\usepackage{mathptmx}       % selects Times Roman as basic font
\usepackage{helvet}         % selects Helvetica as sans-serif font
\usepackage{courier}        % selects Courier as typewriter font
\usepackage{type1cm}        % activate if the above 3 fonts are
\usepackage{makeidx}         % allows index generation
\usepackage{graphicx}        % standard LaTeX graphics tool
\usepackage{multicol}        % used for the two-column index
\usepackage[bottom]{footmisc}% places footnotes at page bottom
\usepackage{mathtools,amssymb,amsmath,amsfonts}
\usepackage{bm,bbm,color}
\usepackage{pifont}
\usepackage[export]{adjustbox}
\usepackage{tikz,pgfplots}
\usepgfplotslibrary{groupplots}
\pgfplotsset{compat=newest}
\usetikzlibrary{arrows,decorations,backgrounds,positioning}
\usepgfplotslibrary{colorbrewer}
\usepackage[noend]{algpseudocode}
\usepackage[ruled]{algorithm2e}
\usepackage[export]{adjustbox}
\def\Letters{A,B,C,D,E,F,G,H,I,J,K,L,M,N,O,P,Q,R,S,T,U,V,W,X,Y,Z}
\makeatletter
\@for \@l:=\Letters \do{%
  \expandafter\edef\csname\@l bb\endcsname{%
  \noexpand\ensuremath{\noexpand\mathbb{\@l}}}%
  \expandafter\edef\csname\@l bf\endcsname{{\noexpand\bf \@l}}%
  \expandafter\edef\csname\@l cal\endcsname{%
  \noexpand\ensuremath{\noexpand\mathcal{\@l}}}%
  \expandafter\edef\csname\@l eu\endcsname{%
  \noexpand\ensuremath{\noexpand\EuScript{\@l}}}%
  \expandafter\edef\csname\@l frak\endcsname{%
  \noexpand\ensuremath{\noexpand\mathfrak{\@l}}}%
  \expandafter\edef\csname\@l rm\endcsname{{\noexpand\rm \@l}}%
  \expandafter\edef\csname\@l scr\endcsname{%
  \noexpand\ensuremath{\noexpand\mathscr{\@l}}}%
}
\newcommand{\bs}[1]{{\boldsymbol#1}}
\renewcommand{\d}{\operatorname{d\!}}

\newcommand{\isdef}{\mathrel{\mathrel{\mathop:}=}}
\newcommand{\defis}{\mathrel{=\mathrel{\mathop:}}}
\newcommand{\dist}{\operatorname{dist}}

\newcommand{\spn}{{\operatorname{span}}}
\newcommand{\diam}{{\operatorname{diam}}}
\newcommand{\supp}{{\operatorname{supp}}}

\DeclareMathOperator{\kernel}{\kappa}
\definecolor{navy}{RGB}{102,153,255}
\definecolor{tuerkis}{RGB}{51,153,204}
\algrenewcommand\alglinenumber[1]{\ding{\taumexpr191 + #1}}

\definecolor{shadecolor}{rgb}{0.6, 0.6, 0.6} 
\definecolor{red}{rgb}{1,0,0.2} 
\definecolor{darkgreen}{rgb}{0, 0.6, 0}
\begin{document}
%%%%%%%%%%%%%%%%%%%%%%%%%%%%%%%%%%%%%%%%%%%%%%%%%%%%%%%%%%%%%%%%%%%%%%%%%%%%%%%%
% Title
\title*{Samplets: Wavelet concepts for scattered data}
\titlerunning{Samplets: Wavelet concepts for scattered data}
\author{Helmut Harbrecht and Michael Multerer}
\institute{Helmut Harbrecht \at
Departement Mathematik und Informatik, 
Universit\"at Basel, 
Spiegelgasse 1, 4051 Basel, Schweiz,
\email{helmut.harbrecht@unibas.ch}
\and Michael Multerer \at
Istituto Eulero,
Universit{\`a} della Svizzera italiana,
Via la Santa 1, 6962 Lugano, Svizzera,
\email{michael.multerer@usi.ch}}
\maketitle

%%%%%%%%%%%%%%%%%%%%%%%%%%%%%%%%%%%%%%%%%%%%%%%%%%%%%%%%%%%%%%%%%%%%%%%%%%%%%%%%
% abstract
\abstract{
This chapter is dedicated to recent developments in the field of wavelet 
analysis for scattered data. We introduce the concept of samplets, which are 
signed measures of wavelet type and may be defined on sets of arbitrarily
distributed data sites in possibly high dimension. By employing samplets, we 
transfer well-known concepts known from wavelet analysis, namely the fast basis 
transform, data compression, operator compression and operator arithmetics
to scattered data problems. Especially, samplet matrix compression facilitates 
the rapid solution of scattered data interpolation problems, even for kernel 
functions with nonlocal support. Finally, we demonstrate that sparsity 
constraints for scattered data approximation problems become meaningful and can 
efficiently be solved in samplet coordinates.}
%%%%%%%%%%%%%%%%%%%%%%%%%%%%%%%%%%%%%%%%%%%%%%%%%%%%%%%%%%%%%%%%%%%%%%%%%%%%%%%%
% Keywords 
% \keywords{Scattered data analysis, wavelets, reproducing kernel Hilbert space}
%%%%%%%%%%%%%%%%%%%%%%%%%%%%%%%%%%%%%%%%%%%%%%%%%%%%%%%%%%%%%%%%%%%%%%%%%%%%%%%%
% AMS subjects
% \subjclass[2000]{41A05, 41A25, 41A58, 65D05}
%%%%%%%%%%%%%%%%%%%%%%%%%%%%%%%%%%%%%%%%%%%%%%%%%%%%%%%%%%%%%%%%%%%%%%%%%%%%%%%%
%\tableofcontents

\section{Introduction}%\label{HMsec:intro}
%===============================================================================
Multiresolution methods and especially wavelet techniques have a long standing 
tradition and are a versatile tool in various fields. Applications comprise, 
among others, nonlinear approximation, image analysis, signal processing and 
machine learning, see for instance \cite{HM_Chui,HM_Dahmen,HM_Daubechies,HM_Dev98,HM_Mallat,
HM_Mallat2016} and the references therein. Starting from a signal, the pivotal idea
of wavelet techniques is the splitting of this signal into its contributions 
relative to a hierarchy of scales. Such a \emph{multiresolution analysis} starts
from an approximation on a coarse scale and successively resolves details, that 
have not been captured so far, at finer scales. Wavelet techniques naturally 
accommodate data compression and adaptivity. The great success of wavelet 
techniques has specifically been triggered by the fast wavelet transform, which
transforms a signal into its wavelet representation and back, that can be
computed with linear cost in terms of the size of the wavelet basis, see
\cite{HM_Coh03,HM_FWT}.

The original construction of wavelets is based on dilations and translations of
a given mother wavelet. This way, a nested sequence of approximation spaces is
obtained and the elements of this sequence are scaled copies of each other. 
This construction of wavelets is limited to structured data, such as uniform 
subdivisions of the real line. Adaptions to deal with bounded intervals have 
been suggested in e.g.~\cite{HM_Alp93,HM_Quak,HM_DKU}, while wavelet constructions on 
manifolds are the topic of e.g.~\cite{HM_STE,HM_HS,HM_PSS97}. An extension to surface 
triangulations, has been suggested in \cite{HM_TW03}, where (multi-)wavelets are 
constructed as linear combinations of functions at a fixed fine scale. The 
stability of the resulting \emph{Tausch-White wavelet basis} follows from its
orthonormality. Another approach to construct a multiresolution analysis on 
unstructured data, for example on graphs, are \emph{diffusion wavelets}, see 
\cite{HM_CM06}. However, there is no linear cost bound for the computation of a
diffusion wavelet basis.

This chapter is a survey on the generalization of the Tausch-White wavelet
construction towards general scattered data. This is achieve by abstracting
the construction \cite{HM_AHK14,HM_TW03} towards discrete signed measures. The
result is multiresolution analysis of discrete signed measures of localized
support that we call \emph{samplets}. Samplets are tailored to the underlying
data sites and may be computed such that the exhibit \emph{vanishing moments},
that is their associated measure integrals. Lowest order samplets, which 
resemble Haar wavelets on scattered data have been considered in the past for
data compression in \cite{HM_RE11}. It is worth mentioning that the construction
of samplets is not limited to the use of polynomial vanishing  moments Indeed,
it is easily be possible to adapt the presented concepts to primitives with
different desired properties.

We present a general construction template for samplets with an arbitrary number
of vanishing moments. This construction can always be performed with linear cost
for a balanced cluster tree for the data sites, even for non-quasi-uniform data
sites. The resulting basis is always orthonormal and hence stable. Due to the
vanishing moments, the coefficients in the representation of scattered data with
respect to samplet coordinates decay fast whenever the data values resemble a
smooth function evaluated at the data sites. This straightforwardly enables data
compression. In contrast, non-smooth regions in the data are indicated by large
samplet coefficients. This, in turn, enables feature detection and extraction. 

Similar to wavelet matrix compression \cite{HM_BCR,HM_DHS,HM_DPS,HM_SCHN,HM_PS}, samplets are
applicable for the compression of kernel matrices. Such matrices arise in the
scattered data approximation context, compare \cite{HM_Fasshauer2007,HM_HSS08,
HM_Schaback2006,HM_Wendland2004,HM_Rasmussen2006,HM_Williams1998}. Kernel matrices are
typically densely populated since the underlying reproducing kernels are 
nonlocal. Nonetheless, the kernels are typically \emph{asymptotically smooth},
meaning that they behave like smooth functions apart from the diagonal. Hence, 
the discretization of such kernel matrices by samplets with vanishing moments
results in quasi-sparse kernel matrices. They can be compressed such that only
a sparse matrix remains. The resulting compression pattern has been derived in
in \cite{HM_HM2}. Furthermore, as shown in \cite{HM_HMSS}, a respective matrix algebra
can be defined.

Samplets provide a meaningful interpretation of \emph{sparsity constraints} for
scattered data, since the representation of data itself becomes sparse. Such 
constraints have a wide applicability in machine learning, statistics, as well
as in signal processing. Examples for the latter are deblurring, feature 
selection and compressive sensing, see \cite{HM_Candes,HM_BasisPursuit,HM_Donoho,
HM_learning,HM_Tao}. In practice, sparsity constraints are imposed by adding an 
$\ell^1$-regularization term to the objective function. Dealing with this 
regularization in an efficient way is especially mandatory for \emph{basis
pursuit}, that is, for decomposing given data into an optimal superposition of
dictionary elements, where optimal means having the smallest $\ell^1$-norm of 
coefficients among all such decompositions, see, for example, \cite{HM_CDS98,HM_MZ,
HM_tropp}. We demonstrate that the basis pursuit problem can efficiently be solved 
within a the samplet basis and reconstruct scattered data using a dictionary of
multiple kernels.

The rest of this chapter is organized as follows. In 
Section~\ref{HM_section:Samplets}, multiresolution analyses for scattered data 
and the concept of samplets are introduced. The change of basis by means of the
fast samplet transform and resulting applications are the topic of 
Section~\ref{HMsec:Applications}. Section~\ref{HMsec:kernelCompression} deals
with scattered data approximation in reproducing kernel Hilbert spaces. Here, we
especially introduce the samplet compression of kernel matrices and the 
efficient treatment of sparsity constraints by the semi-smooth Newton method.

Throughout this chapter, in order to avoid the repeated use of generic but 
unspecified constants, by \(C\lesssim D\) we indicate that $C$ can be bounded by
a multiple of $D$, independently of parameters which $C$ and $D$ may depend on.
Moreover, \(C\gtrsim D\) is defined as \(D\lesssim C\) and \(C\sim D\) as 
\(C\lesssim D\) and \(D\lesssim C\).

%%%%%%%%%%%%%%%%%%%%%%%%%%%%%%%%%%%%%%%%%%%%%%%%%%%%%%%%%%%%%%%%%%%%%%%%%%%%%%%%
\section{Multiresolution analysis of scattered data}
\label{HM_section:Samplets}
%%%%%%%%%%%%%%%%%%%%%%%%%%%%%%%%%%%%%%%%%%%%%%%%%%%%%%%%%%%%%%%%%%%%%%%%%%%%%%%%
\subsection{Samplet bases}
Let \(X\isdef\{{\bs x}_1,\ldots,{\bs x}_N\}\subset\Omega\) denote a set of 
\emph{data sites} within some bounded or unbounded region 
\(\Omega\subset\Rbb^d\). Associated to these data sites, we introduce the 
Dirac-$\delta$-distributions
\begin{equation}\label{HM_eg:diracs} 
\delta_{{\bs x}_1},\ldots,\delta_{{\bs x}_N}\in [C(\Omega)]',
\end{equation}
where we endow \(C(\Omega)\) with the \(\sup\)-norm. The 
Dirac-$\delta$-distributions satisfy
\[
(f,\delta_{{\bs x}_i})_\Omega\isdef
\delta_{{\bs x}_i}(f)
=f({\bs x}_i)\ \text{for all}\ f\in C(\Omega)
\]
and serve as \emph{information functionals} of some, maybe unknown, function
\( f\in C(\Omega)\). In this setting, the \emph{data values} 
\(f_i\isdef(f,\delta_{{\bs x}_i})_\Omega\), \(i=1,\ldots,N\), amount to the 
available information of the function \(f\), compare \cite{HM_MR77}. 

By the Riesz representation theorem, the dual space \([C(\Omega)]'\) is 
isometrically isomorphic to the space of regular, finitely additive signed 
measures of bounded variation on $\Omega$, equipped with the total variation
norm. We call this space by \(\operatorname{rba}(\Omega)\), as usual. For each 
\(\ell\in[C(\Omega)]'\), there exists a signed measure 
\(\mu_\ell\in\operatorname{rba}(\Omega)\) such that
\[
\ell_{\mu}(f) = \int_{\Omega} f({\bs x}) \d\mu({\bs x})\ \text{for all}\ 
f\in C(\Omega),
\]
and vice versa, see \cite{HM_DS58}. In view of this identification, we 
interchangeably consider the Dirac $\delta$-distributions from 
\eqref{HM_eg:diracs} as linear functionals or as measures.

The Dirac-$\delta$-distributions span the linear subspace
\[
\Xcal\isdef\spn\{\delta_{{\bs x}_1},
\ldots,\delta_{{\bs x}_N}\}\subset [C(\Omega)]'
\]
containing all discrete and finite signed measures supported at the data sites
in \(X\). Endowed with the bilinear form 
\(\langle\cdot,\cdot\rangle_\Xcal\colon\Xcal\times\Xcal\to\Rbb\)
with
\(\langle\delta_{{\bs x}_i},\delta_{{\bs x}_j}\rangle_\Xcal=m_{i,j}\),
where \({\bs M}=[m_{i,j}]_{i,j=1}^N\in\Rbb^{N\times N}\) is a symmetric and 
positive definite matrix, the space \(\Xcal\) becomes a Hilbert space. For 
arbitrary \(u',v'\in\Xcal\), the inner product is given by
\[
\langle u',v'\rangle_\Xcal=\sum_{i,j=1}^N m_{i,j}u_iv_j,\ \text{where}\
u'=\sum_{i=1}^Nu_i\delta_{{\bs x}_i},\ v'=\sum_{i=1}^Nv_i\delta_{{\bs x}_i}.
\]
The canonical choice for \({\bs M}\) is the identity matrix, which renders the
space \(\Xcal\) is isometrically isomorphic to \(\Rbb^N\) with  the Euclidean
inner product. 

Given a \emph{multiresolution analysis}
\[
\Xcal_0\subset\Xcal_1\subset\cdots\subset\Xcal_J=\Xcal,
\]
we keep track of the increment of information between two consecutive levels $j$
and $j+1$. Since there holds $\Xcal_{j}\subset \Xcal_{j+1}$, we may
(orthogonally) decompose 
\[
\Xcal_{j+1} = \Xcal_j\overset{\perp}{\oplus}\Scal_j
\]
for a certain \emph{detail space} $\Scal_j$. In analogy to wavelet nomenclature,
we call the elements of a basis of \(\Xcal_0\) \emph{scaling distributions} and
the elements of a basis of one of the spaces $\Scal_j$ \emph{samplets}. This
name is motivated by the idea that the basis distributions in $\Scal_j$ are 
supported at a small subsample or samplet of the data sites in \(X\). The
collection of the bases of $\Scal_j$ for \(j=0,\ldots,J-1\)\/ together with a 
basis of \(\Xcal_0\) is called a \emph{samplet basis} for \(\Xcal\). 

To enable data compression, we assume that the samplets have
\emph{vanishing moments} of order \(q+1\), viz.,
\begin{equation}\label{HM_eg:vanishingMoments}
%===============================================================================
(p,\sigma_{j,k})_\Omega
 = 0\ \text{for all}\ p\in\Pcal_q,
\end{equation}
where \(\Pcal_q\isdef\operatorname{span}
\{{\bs x}^{\bs\alpha}:\|\bs\alpha\|_1\leq q\}\)
denotes the space of all polynomials with total degree at most \(q\). In view of
larger spatial dimension \(d\), we observe that the dimension of \(\Pcal_q\) is
\begin{equation}\label{HM_eg:mq}
%============================================
m_q\isdef\dim(\Pcal_q)=
\binom{q+d}{d} =
\frac{(d+q)\cdots(d+1)}{q!}=\Ocal(d^q).
\end{equation}
Hence, the number of vanishing moment conditions only grows polynomially with
respect to \(d\). 

Following the discrete signed measure interpretation of samplets, equation
\eqref{HM_eg:vanishingMoments} implies that a samplet amounts to a quadrature 
formula, which annihilates polynomials up to order \(q+1\). 

%===============================================================================
\subsection{Cluster tree}
%===============================================================================
To construct an orthonormal samplet basis with vanishing moments, see
\eqref{HM_eg:vanishingMoments}, we adapt the wavelet construction from 
\cite{HM_TW03} to the scattered data framework. The first step is to construct a 
hierarchy of subspaces of signed measures. To this end, we perform a 
hierarchical clustering of the set \(X\) with respect to the  Euclidean 
distance. This approach amounts to a clustering of the 
Dirac-$\delta$-distributions spanning \(\Xcal\) by means of the Wasserstein 
distance
\(
  W_p(\delta_{{\bs x}_i},\delta_{{\bs x}_j})=\|{\bs x}_i-{\bs x}_j\|_2
\)
for any \(p\geq 1\).

\begin{definition}\label{def:cluster-tree}
%===============================================================================
Let $\mathcal{T}=(V,E)$ be a tree with vertices $V$ and edges $E$. 
For \(\tau\in V\), we denote the set of \emph{children} of \(\tau\) by 
\(\operatorname{child}(\tau)\). The set of \emph{leaves} is defined as 
$\mathcal{L}(\mathcal{T})\isdef\{\tau\in V:\operatorname{child}(\tau)
=\emptyset\}$. 

The tree $\mathcal{T}$ is a \emph{cluster tree} for the set 
$X$, iff $X$ is the root of $\mathcal{T}$ and all
$\tau\in V\setminus\mathcal{L}(\mathcal{T})$ are disjoint unions of their 
children, viz., \(\tau=\cup_{\tau'\in\operatorname{child}(\tau)}\tau'\).
The \emph{level} \(\operatorname{level}(\tau)\) of $\tau\in\mathcal{T}$
is its distance from the root. Finally, the \emph{depth} \(J\) of \(\Tcal\) is 
the maximum level of all clusters. 
\end{definition}

\begin{remark}
There exist different possibilities to construct a cluster tree for the set 
\(X\). We focus on binary trees and remark that other options, such as 
\(2^d\)-trees, may be considered with the obvious modifications.
\end{remark}

Definition~\ref{def:cluster-tree} provides a hierarchical structure for the set
\(X\). Even so, it does not provide guarantees for the cardinalities of the 
clusters. To account for this issue, we introduce the concept of a 
\emph{balanced binary tree}.

\begin{definition}
Let $\Tcal$ be a cluster tree for $X$ with depth $J$. We call $\Tcal$ a 
\emph{balanced binary tree}, iff all clusters $\tau$ satisfy the following two
conditions:
\begin{enumerate}
\item[(\emph{i})]
The cluster $\tau$ has exactly two children
if $\operatorname{level}(\tau) < J$. It has no children 
if $\operatorname{level}(\tau) = J$.
\item[(\emph{ii})]
There holds $|\tau|\sim 2^{J-\operatorname{level}(\tau)}$, where $|\tau|$ 
denotes the number of points contained in \({\tau}\).
\end{enumerate}
\end{definition}

A balanced binary tree can be constructed by \emph{cardinality balanced 
clustering}. This means that the root cluster \(X\) is split into two child 
clusters of identical (or similar) cardinality. This process is repeated 
recursively for the resulting child clusters until their cardinality falls 
below a given threshold. For the subdivision, the cluster's bounding box is 
split along its longest edge such that the resulting two boxes both contain an
equal number of points. Thus, as the cluster cardinality halves with each level,
we obtain $\mathcal{O}(\log N)$ levels in total. The cost per level is 
$\mathcal{O}(N)$ if the subdivision is performed with respect to the median,
compare \cite{HM_BFPR+73}. The overall cost for constructing the cluster
tree is therefore $\mathcal{O}(N \log N)$. 

We remark that a balanced tree is only required to guarantee the cost bounds for
the presented algorithms. The error and compression estimates presented later on
are robust in the sense that they are formulated directly in terms of the actual
cluster sizes rather than the introduced cluster level.

%===============================================================================
\subsection{Construction of the samplets}
%===============================================================================
Given a cluster tree, we start by introducing a \emph{two-scale transform} 
between basis elements on a cluster $\tau$ of level $j$. To this end, we 
represent scaling distributions $\bs{\Phi}_{j}^{\tau}=\{ \varphi_{j,k}^{\tau}\}$
and samplets ${\bs\Sigma}_{j}^{\tau}=\{\sigma_{j,k}^{\tau}\}$ recursively as 
linear combinations of the scaling distributions ${\bs \Phi}_{j+1}^{\tau}$ of 
$\tau$'s child clusters. This amounts to the \emph{refinement relations}
\[
\varphi_{j,k}^{\tau}=\sum_{\ell=1}^{n_{j+1}^\tau}q_{j,\Phi,\ell,k}^{\tau}
\varphi_{j+1,\ell}^{\tau}
\ \text{and}\ 
\sigma_{j,k}^{\tau}=\sum_{\ell=1}^{n_{j+1}^\tau}q_{j,\Sigma,\ell,k}^{\tau}
\varphi_{j+1,\ell}^{\tau}\ \text{with}\ n_{j+1}^\tau\isdef|
{\bs \Phi}_{j+1}^{\tau}|,
\]
which may be written in matrix notation as
\begin{equation}\label{HM_eg:refinementRelation}
%===============================================================================
   [ {\bs \Phi}_{j}^{\tau}, {\bs \Sigma}_{j}^{\tau} ] 
 = {\bs \Phi}_{j+1}^{\tau}
 {\bs Q}_j^{\tau}=
 {\bs \Phi}_{j+1}^{\tau}
 \big[ {\bs Q}_{j,\Phi}^{\tau},{\bs Q}_{j,\Sigma}^{\tau}\big].
\end{equation}

To obtain vanishing moments and orthonormality, the transform 
\({\bs Q}_{j}^{\tau}\) is derived from an orthogonal decomposition of the 
\emph{moment matrix} ${\bs M}_{j+1}^{\tau}\in\Rbb^{m_q\times n_{j+1}^\tau}$,
given by
\[
  {\bs M}_{j+1}^{\tau}\isdef
  \begin{bmatrix}({\bs x}^{\bs 0},\varphi_{j+1,1})_\Omega&\cdots&
  ({\bs x}^{\bs 0},\varphi_{j+1,n_{j+1}^\tau})_\Omega\\
  \vdots & & \vdots\\
  ({\bs x}^{\bs\alpha},\varphi_{j+1,1})_\Omega&\cdots&
  ({\bs x}^{\bs\alpha},\varphi_{j+1,n_{j+1}^\tau})_\Omega
  \end{bmatrix}=
  [({\bs x}^{\bs\alpha},{\bs \Phi}_{j+1}^{\tau})_\Omega]_{|\bs\alpha|\le q}.
\]
Herein, $m_q$ is given as in \eqref{HM_eg:mq} and denotes the dimension of the 
space \(\Pcal_q(\Omega)\) of total degree polynomials.

In the original construction by Tausch and White, the matrix 
\({\bs Q}_{j}^{\tau}\) is obtained from the singular value decomposition of 
\({\bs M}_{j+1}^{\tau}\). For the construction of samplets, we follow the idea
from \cite{HM_AHK14} and rather employ the QR decomposition, which results in 
samplets with an increasing number of vanishing moments. We compute
\begin{equation}\label{HM_eg:QR} 
%===============================================================================
  ({\bs M}_{j+1}^{\tau})^\intercal  = {\bs Q}_j^\tau{\bs R}
  \defis\big[{\bs Q}_{j,\Phi}^{\tau} ,
  {\bs Q}_{j,\Sigma}^{\tau}\big]{\bs R}
 \end{equation}
The moment matrix for the cluster's scaling distributions and samplets is now 
given by 
\begin{equation*}%\label{HM_eg:vanishingMomentsQR}
%===============================================================================
  \begin{aligned}
  \big[{\bs M}_{j,\Phi}^{\tau}, {\bs M}_{j,\Sigma}^{\tau}\big]
  &= \left[({\bs x}^{\bs\alpha},[{\bs \Phi}_{j}^{\tau},
  	{\bs \Sigma}_{j}^{\tau}])_\Omega\right]_{|\bs\alpha|\le q}\\
  &= \left[({\bs x}^{\bs\alpha},{\bs \Phi}_{j+1}^{\tau}[{\bs Q}_{j,\Phi}^{\tau}
  , {\bs Q}_{j,\Sigma}^{\tau}])_\Omega
  	\right]_{|\bs\alpha|\le q}\\
  &= {\bs M}_{j+1}^{\tau} [{\bs Q}_{j,\Phi}^{\tau},{\bs Q}_{j,\Sigma}^{\tau}]\\
  &= {\bs R}^\intercal.
  \end{aligned}
\end{equation*}
Since ${\bs R}^\intercal$ is a lower triangular matrix, the first $k-1$ entries
in its $k$-th column are zero. This amounts to $k-1$ vanishing moments for the
$k$-th distribution generated by the orthogonal transform
${\bs Q}_{j}^{\tau}=[{\bs Q}_{j,\Phi}^{\tau} , {\bs Q}_{j,\Sigma}^{\tau} ]$. 
Defining the first $m_{q}$ distributions as scaling distributions and the 
remaining ones as samplets, we obtain samplets with vanishing moments at least
of order $q+1$. 

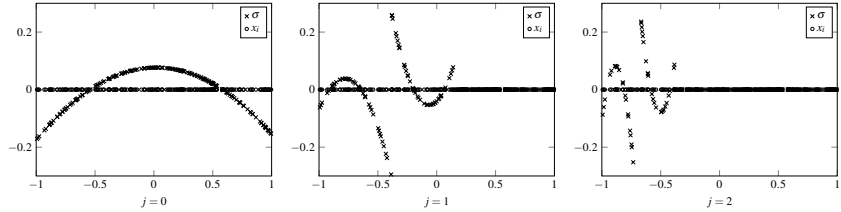
\begin{figure}[htb]
\begin{center}
\scalebox{0.52}{
\begin{tikzpicture}
\begin{axis}[width=7.6cm, height=6cm, xmin = -1, xmax=1,
 ymin=-0.3, ymax=0.3, ylabel={}, xlabel ={$j=0$}, ytick={-0.2,0,0.2}]
\addplot[color=black, only marks, mark=x,line width=0.6pt,mark size=2pt]
  table[x index={0},y index = {1}]{%
./Images/samplets.txt};
\addlegendentry{$\sigma$};
\addplot[color=black, mark=o, line width=0.4pt,mark size=1.2pt, only marks]
  table[x index={0},y index = {1}]{%
./Images/dpts.txt};
\addlegendentry{$x_i$};
\end{axis}
\end{tikzpicture}}
\scalebox{0.52}{
\begin{tikzpicture}
\begin{axis}[width=7.6cm, height=6cm, xmin = -1, xmax=1,
 ymin=-0.3, ymax=0.3, ylabel={}, xlabel ={$j=1$},ytick={-0.2,0,0.2}]
\addplot[color=black, only marks, mark=x,line width=0.6pt,mark size=2pt]
  table[x index={0},y index = {2}]{%
./Images/samplets.txt};
\addlegendentry{$\sigma$};
\addplot[color=black, mark=o, line width=0.4pt,mark size=1.2pt, only marks]
  table[x index={0},y index = {1}]{%
./Images/dpts.txt};
\addlegendentry{$x_i$};
\end{axis}
\end{tikzpicture}}
\scalebox{0.52}{
\begin{tikzpicture}
\begin{axis}[width=7.6cm, height=6cm, xmin = -1, xmax=1,
 ymin=-0.3, ymax=0.3, ylabel={}, xlabel ={$j=2$},ytick={-0.2,0,0.2}]
\addplot[color=black, only marks, mark=x,line width=0.6pt,mark size=2pt]
  table[x index={0},y index = {3}]{%
./Images/samplets.txt};
\addlegendentry{$\sigma$};
\addplot[color=black, mark=o, line width=0.4pt,mark size=1.2pt, only marks]
  table[x index={0},y index = {1}]{%
./Images/dpts.txt};
\addlegendentry{$x_i$};
\end{axis}
\end{tikzpicture}}
\caption{\label{HM_fig:samplets0}Scaling distribution (left), samplet on level 
\(j=1\) (middle) and samplet on level \(j=2\) (right) for \(N=200\) uniformly
distributed data sites and \(q+1=3\).}
\end{center}
\end{figure}

The orthogonality of basis elements of two different clusters results from the
non-overlapping supports, while orthogonality within a given branch is a
consequence of the orthogonality of the transforms. A visualization of a scaling
distribution and samplets on levels \(j=1\) and \(j=2\) with \(q+1=3\) vanishing
moments on \(N=200\) uniformly distributed data sites is shown in
Figure~\ref{HM_fig:samplets0}.

If we choose a minimum leaf size \(|\tau|\geq m_{\hat{q}}\geq 2m_q\) for the
cluster tree for some polynomial degree \(\hat{q}> q\), we can even construct 
samplets whose number of vanishing moments increase successively from order
\(q+1\) up to order \(\hat{q}+1\) without additional cost. This is advantageous
since more vanishing moments typically improve the a-posteriori compression
ratios of a given signal. Moreover, the preceding algebraic construction
\eqref{HM_eg:QR} of vanishing moments can easily be adapted to other primitives,
for example, anisotropic total degree polynomial spaces. This particularly makes
sense in the higher dimensional setting, when dimension weights are available.

\begin{remark}
%\label{remark:introCQ}
Each cluster has at most a constant number of scaling distributions and 
samplets. For leaf clusters, this number is bounded by the leaf size. For
non-leaf clusters, it is bounded by the number of scaling distributions from its
child clusters. As there are two child clusters with a maximum of $m_q$ scaling
distributions each, we obtain the bound $2 m_q$ for non-leaf clusters. If the
cluster tree is balanced, resulting in a leaf size of \(\Ocal(1)\), the above
construction of the samplet basis has hence cost \(\Ocal(N)\).
\end{remark}

For leaf clusters, we define the scaling distributions by the Dirac measures
supported at the points \({\bs x}_i\), viz.,
${\bs \Phi}^{\tau}\isdef\{ \delta_{{\bs x}_i} : {\bs x}_i\in\tau,
\tau\in\Lcal(\Tcal) \}$.
The scaling distributions of all clusters on a specific level $j$ then generate
the spaces
\begin{equation}\label{HM_eg:Vspaces}
%===============================================================================
	\Xcal_{j}\isdef \spn\big\{ \varphi_{j,k}^{\tau} : 
 k\in I_j^{\Phi,\tau},\ \operatorname{level}(\tau)=j \big\}.
\end{equation}
In contrast, the samplets span the detail spaces
\begin{equation}\label{HM_eg:Wspaces}
%===============================================================================
\Scal_{j}\isdef
	\spn\big\{ \sigma_{j,k}^{\tau} : k\in I_j^{\Sigma,\tau},\ 
	\operatorname{level}(\tau)=j\big\}.
\end{equation}
Combining the scaling distributions of the root cluster with all clusters' 
samplets gives rise to the samplet basis
\begin{equation}\label{HM_eg:Wbasis}
%===============================================================================
  {\bs \Sigma}\isdef{\bs \Phi}^{X} 
  	\cup \bigcup_{\tau \in\Tcal} {\bs \Sigma}^{\tau}.
\end{equation}

%===============================================================================
\subsection{Properties of the samplets}
%===============================================================================
The properties of the samplet basis constructed in the previous paragraph are 
summarized in the next theorem, which can be inferred by adapting the 
corresponding results from \cite{HM_HKS05,HM_TW03}.

\begin{theorem}%\label{theo:waveletProperties}
%===============================================================================
The spaces $\Xcal_{j}$ defined in equation \eqref{HM_eg:Vspaces} form a 
multiresolution analysis
\[
  \Xcal_0\subset\Xcal_1\subset\cdots\subset\Xcal_J =\Xcal,
\]
where the respective complement spaces $\Scal_{j}$ from \eqref{HM_eg:Wspaces} 
satisfy 
\[
\Xcal_{j+1}=\Xcal_j\overset{\perp}{\oplus}\Scal_{j}
\ \text{for all}\ j=0,1,\ldots,J-1.
\]
The associated samplet basis ${\bs \Sigma}$ defined in \eqref{HM_eg:Wbasis} is an
orthonormal basis in $\Xcal$. In particular, there holds:
\begin{enumerate}
\item[(\emph{i})] The number of all samplets on level $j$ behaves like $2^j$.
\item[(\emph{ii})] The samplets have $q+1$ vanishing moments.
\item[(\emph{iii})] 
Each samplet is supported in a specific cluster $\tau$. 
\end{enumerate}
\end{theorem}

The key for data compression and feature detection is the following estimate
which shows that the samplet coefficients decay with respect to the samplet's
support size provided that the data result from the evaluation of a smooth 
function on the samplet's support, compare \cite{HM_HM2}. 

\begin{lemma}\label{lemma:decay}
%===============================================================================
Given a samplet \(\sigma_{j,k}\) with support \(\tau\) and let 
\(f\in  C(\Omega)\) with $f\in C^{q+1}(O)$ for any open set \(O\supset\tau\). 
Then, there holds
\begin{equation*}%\label{HM_eg:decay}
%===============================================================================
 |(f,\sigma_{j,k})_\Omega|\le\sqrt{|\tau|} \bigg(\frac{d}{2}\bigg)^{q+1}
  	\frac{\diam(\tau)^{q+1}}{(q+1)!}\|f\|_{C^{q+1}(O)}.
\end{equation*}
\end{lemma}

Hence, in case of smooth data, the samplet coefficients are small and can be set
to zero without compromising the accuracy. Vice versa, a large samplet 
coefficient indicates that the data are singular in the region of the samplet's
support. If the set \(X\) is \emph{quasi-uniform} in the sense that the
\emph{separation radius}
\(
q_X\isdef\min_{i\neq j}\|{\bs x}_i-{\bs x}_j\|_2\)
and the
\emph{fill distance}
\(
h_{X,\Omega}\isdef\sup_{{\bs x}\in\Omega}\min_{{\bs x}_i\in X}
\|{\bs x}-{\bs x}_i\|_2
\)
satisfy \(q_X\sim h_{X,\Omega}\), then there holds
$\diam(\tau)\sim 2^{-\operatorname{level}(\tau)/d}$. 
In this case, Lemma~\ref{lemma:decay} even guarantees the exponential decay of
the samplet coefficients with the level provided that the underlying signal is
smooth.

%===============================================================================
\section{Scattered data analysis}\label{HMsec:Applications}
%===============================================================================
\subsection{Fast samplet transform}\label{HMsec:FST}
%===============================================================================
To transform between the samplet basis and the basis of Dirac measures, we 
introduce the \emph{fast samplet transform} and its inverse. In accordance with
\cite{HM_TW03}, this transform and its inverse can be performed in linear cost. 
This result is well-known in case of wavelets and was crucial for their rapid 
development.

The mapping 
\(f\mapsto\big[(f,\delta_{{\bs x}_1})_{\Omega},\ldots,
(f,\delta_{{\bs x}_N})_{\Omega}\big]^\intercal\) defines an operator 
\(I\colon C(\Omega)\to\Rbb^N\), called \emph{information operator} in the 
context of optimal recovery, see \cite{HM_MR77}. By employing samplets, this
information is represented in a multiresolution fashion
\(f\mapsto\big[(f,\sigma_{j,k})_{\Omega}\big]^\intercal_{j,k}\in\Rbb^N\).
More precisely, letting 
\begin{equation}\label{HM_eg:InfOp}
%=========================================
{\bs f}^\Delta\isdef If\quad\text{and}\quad 
{\bs f}^\Sigma\isdef[(f,\sigma_{j,k})_{\Omega}\big]^\intercal,
\end{equation}
the fast samplet transform computes
\[
{\bs f}^{\Sigma}={\bs T}{\bs f}^{\Delta},
\]
where \({\bs T}\in\Rbb^{N\times N}\) is the orthogonal matrix containing the
expansion coefficients of the samplet basis.

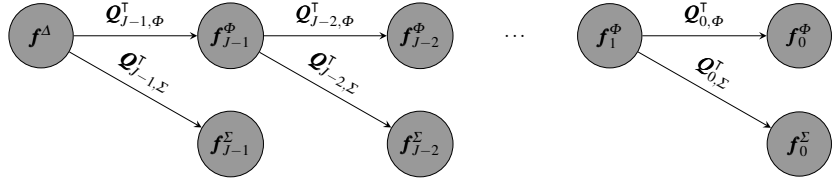
\begin{figure}[htb]
\begin{center}
\scalebox{0.9}{
\begin{tikzpicture}[x=0.35cm,y=0.4cm]
\tikzstyle{every node}=[circle,draw=black,fill=shadecolor,
minimum size=0.95cm]%
\tikzstyle{ptr}=[draw=none,fill=none,above]%
\node at (0,5) (1) {\textcolor{black}{\small${\bs f}^{\Delta}$}};
\node at (8,5) (2) {\textcolor{black}{\small${\bs f}_{J-1}^{\Phi}$}};
\node at (8,1) (3) {\textcolor{black}{\small${\bs f}_{J-1}^{\Sigma}$}};
\node at (16,5) (4) {\textcolor{black}{\small${\bs f}_{J-2}^{\Phi}$}};
\node at (16,1) (5) {\textcolor{black}{\small${\bs f}_{J-2}^{\Sigma}$}};
\node at (24,5) (6) {\textcolor{black}{\small${\bs f}_{1}^{\Phi}$}};
\node at (32,5) (7) {\textcolor{black}{\small${\bs f}_{0}^{\Phi}$}};
\node at (32,1) (8) {\textcolor{black}{\small${\bs f}_{0}^{\Sigma}$}};
\tikzstyle{forward}=[draw,-stealth]%
\tikzstyle{every node}=[style=ptr]
\draw
(1) edge[forward] node[above,sloped]{\textcolor{black}{\small${\bs Q}_{J-1,\Phi
}^\intercal$}} (2)
(1) edge[forward] node[above,sloped]{\textcolor{black}{\small${\bs Q}_{J-1,
\Sigma}^\intercal$}}%
 (3)
(2) edge[forward] node[above,sloped]{\textcolor{black}{\small${\bs Q}_{J-2,\Phi
}^\intercal$}} (4)
(2) edge[forward] node[above,sloped]{\textcolor{black}{\small${\bs Q}_{J-2,
\Sigma}^\intercal$}}%
 (5)
(6) edge[forward] node[above,sloped]{\textcolor{black}{\small${\bs Q}_{0,\Phi
}^\intercal$}} (7)
(6) edge[forward] node[above,sloped]{\textcolor{black}{\small${\bs Q}_{0,
\Sigma}^\intercal$}}%
  (8);
\tikzstyle{every node}=[style=ptr]%
\tikzstyle{ptr}=[draw=none,fill=none]%
\node at (20,5) (16) {$\hdots$};
\end{tikzpicture}}
\caption{Fishbone scheme of the fast samplet transform.}
\label{HM_fig:haar}
\end{center}
\end{figure}

The actual implementation of the fast samplet transform is recursive and follows
the fishbone scheme of the fast wavelet transform, see Figure~\ref{HM_fig:haar} for
a respective illustration of this scheme. The underlying recursion is based on
the refinement relation \eqref{HM_eg:refinementRelation}, which translates to
\begin{equation}\label{HM_eg:refinementRelationInnerProducts}
%===============================================================================
(f, [ {\bs \Phi}_{j}^{\tau}, {\bs \Sigma}_{j}^{\tau} ])_\Omega
=(f, {\bs \Phi}_{j+1}^{\tau} [{\bs Q}_{j,\Phi}^{\tau} ,{\bs Q}_{j,\Sigma}^{\tau}
])_\Omega
=(f, {\bs \Phi}_{j+1}^{\tau})_\Omega [{\bs Q}_{j,\Phi}^{\tau} ,
{\bs Q}_{j,\Sigma}^{\tau} ].
\end{equation}
On the finest level, the entries of the vector $(f,{\bs \Phi}_{J}^{\tau}
)_\Omega$ are exactly those of ${\bs f}^{\Delta}$. Recursive application of
\eqref{HM_eg:refinementRelationInnerProducts} therefore yields all the 
coefficients $(f, {\bs \Sigma}_{j}^{\nu})_\Omega$, including 
$(f, {\bs \Phi}_{0}^{X})_\Omega$, required for the representation of 
${\bs f}^\Delta$ in the samplet basis. 

The inverse transform is obtained in complete analogy by reversing the steps of
the fast samplet transform: For each cluster, we compute
\[
(f, {\bs \Phi}_{j+1}^{\tau})_\Omega
= (f, [ {\bs \Phi}_{j}^{\tau}, {\bs \Sigma}_{j}^{\tau} ] 
)_\Omega[{\bs Q}_{j,\Phi}^{\tau} ,{\bs Q}_{j,\Sigma}^{\tau} ]^\intercal
\] 
to either obtain the coefficients of the child clusters' scaling distributions 
or, for leaf clusters, the coefficients ${\bs f}^{\Delta}$.

%===============================================================================
\subsection{Data compression}
%===============================================================================
In view of Lemma~\ref{lemma:decay}, samplets enable data compression by means of
so-called \emph{hard-thresholding}. Given the coefficients \({\bs f}^\Sigma\) of
a transformed signal, we define \(\operatorname{HT}_{w}({\bs f}^\Sigma)\) as the
operator which sets all entries \(f^{\Sigma}_i\) with \(|f^{\Sigma}_i|<w\) to
zero. To showcase this application of samplets, we consider the monthly ERA5
temperature data set, which provides the temperature two meters above the
earth's surface. ERA5 is a reanalysis by the European Center for Medium-range 
Weather Forecasts (ECMWF) of global climate and weather for the past eight 
decades. Data is available from 1940 onwards, see \cite{HM_ERA5}. We use the 
monthly data from 2022, which comprises, after removing \texttt{NaN} values,
1\,038\,240 data points per month. This results in 12\,458\,880 points in total.
The original data format is World Geodesic System 1984 (WGS 84), which we 
transform by using the Robinson projection, see \cite{HM_Rob74}, to reduce the 
distortion and the projection angles. The bounding box for the data is 
\([-1.70\cdot 10^7, 1.70\cdot 10^7]\times [-8.63\cdot 10^6, 8.63\cdot 10^6]
\times[0,11]\). In particular, the data points are not located on a regular 
spatial grid and cluster towards the North and South pole. The fill distance
with respect to the convex hull of the data points is \(h_{X,\Omega}= 23\,700\) 
for each time slice, while the separation radius is \(q_X=7\,260\) for each time
slice.

\begin{figure}[htb]
\begin{center}
\begin{tikzpicture}
\draw(0,0)node{\includegraphics[scale=0.04,clip,trim=500 350 500 350]{
  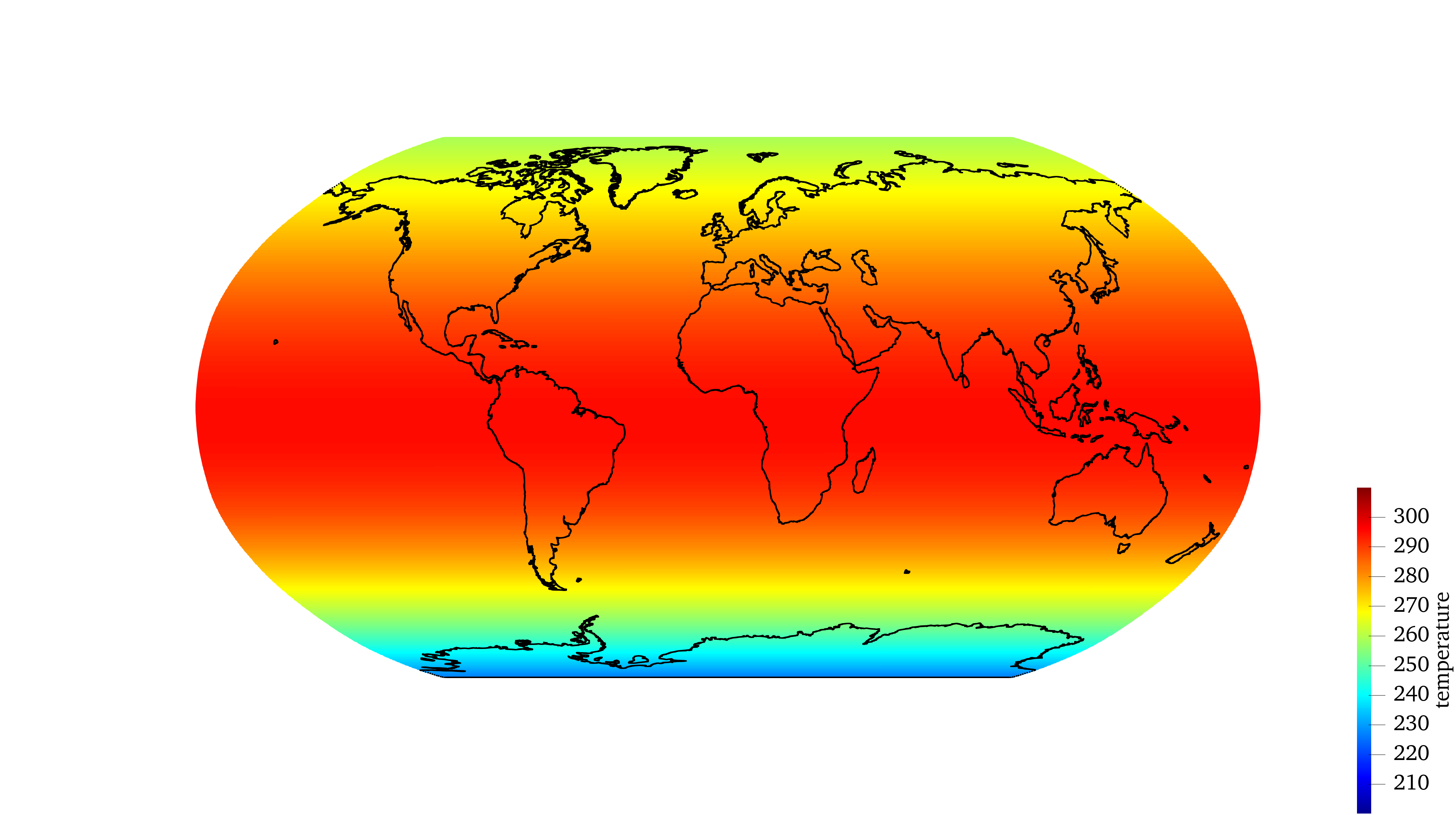}};
\draw(0,-1.25)node{\tiny Space saving: \(99.9996\%\)};
\draw(5,0)node{\includegraphics[scale=0.04,clip,trim=500 350 500 350]{
  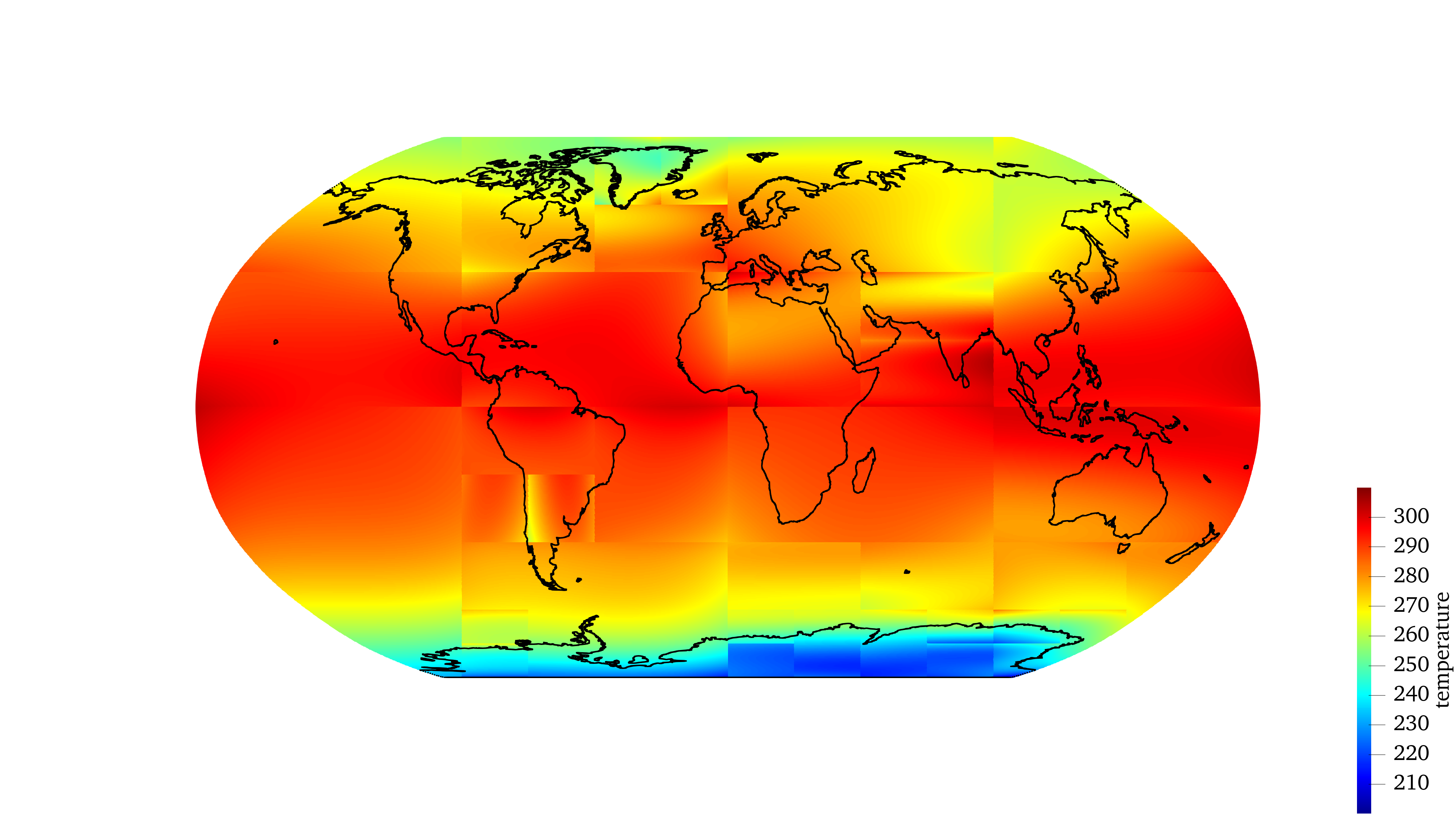}};
\draw(5,-1.25)node{\tiny Space saving: \(99.9916\%\)};
\draw(0,-3)node{\includegraphics[scale=0.04,clip,trim=500 350 500 350]{
  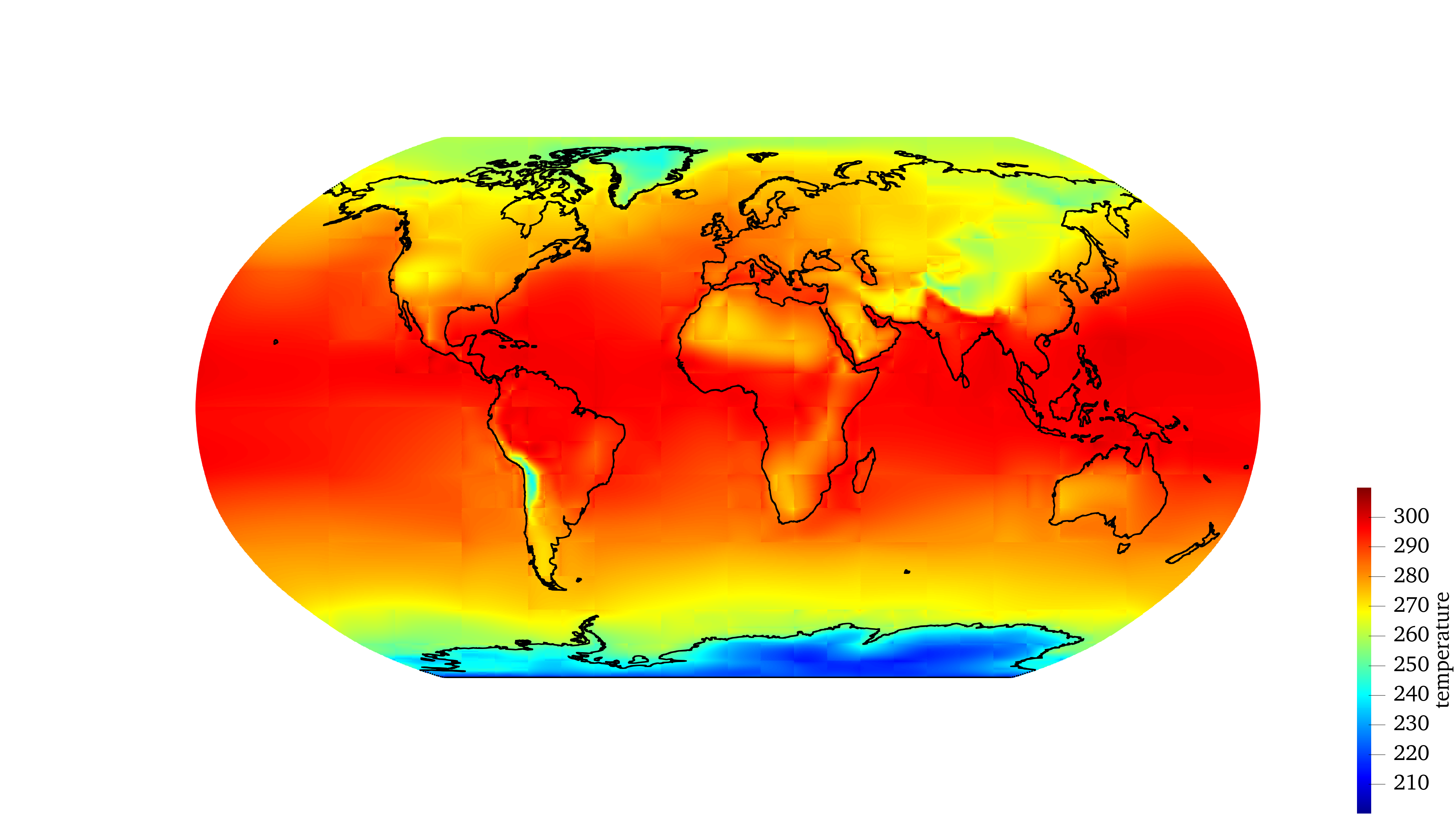}};
\draw(0,-4.25)node{\tiny Space saving: \(99.8600\%\)};
\draw(5,-3)node{\includegraphics[scale=0.04,clip,trim=500 350 500 350]{
  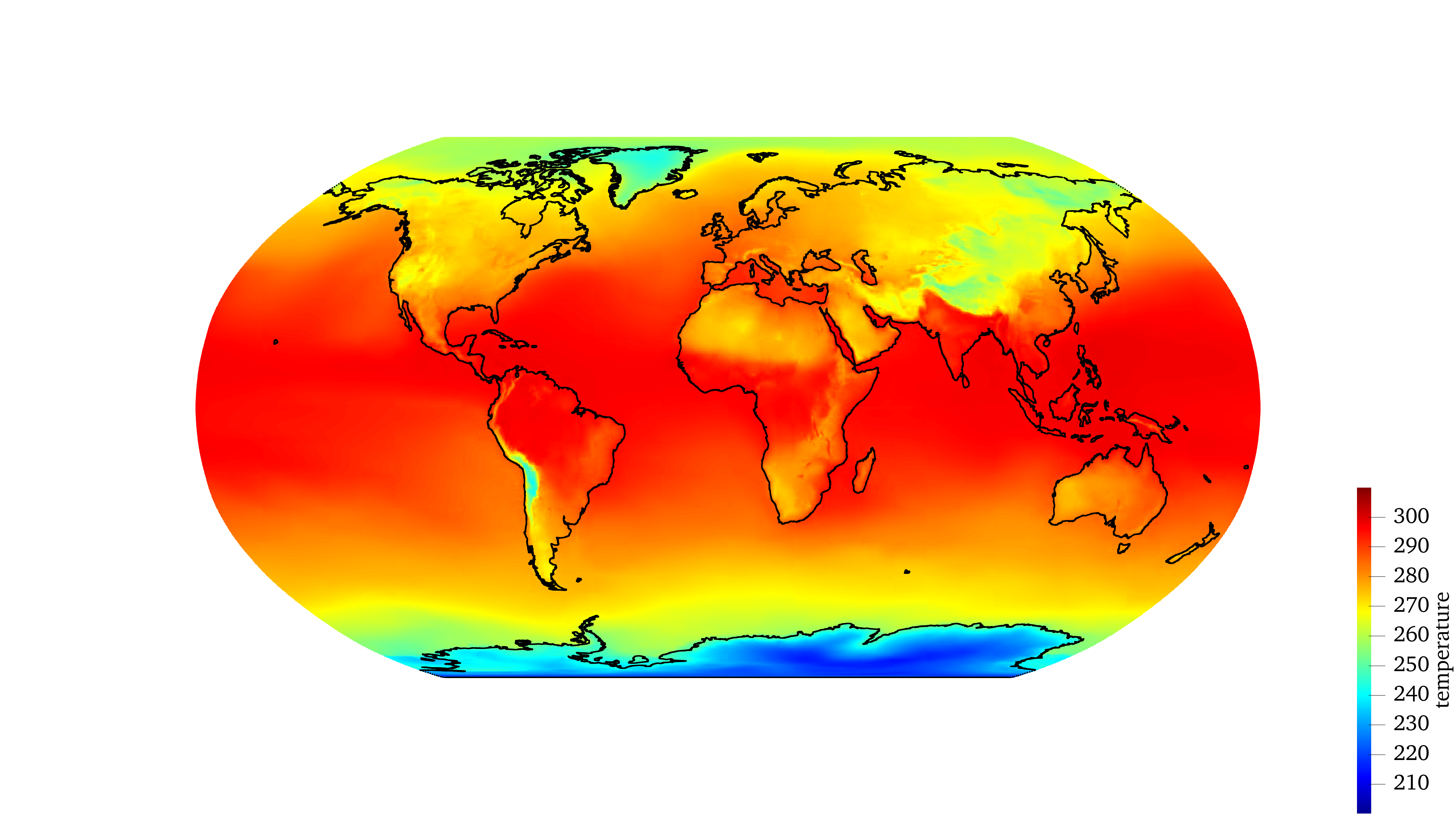}};
\draw(5,-4.25)node{\tiny Space saving: \(98.489\%\)};
\draw(8,-1.5)node{\includegraphics[scale=0.12,clip,trim=3620 0 0 1250]{
  ./Images/HT_rec_1e-2.png}};
\end{tikzpicture}
\caption{\label{HM_fig:temp_HT_rec}Reconstruction of the temperature based on the
  hard-thresholded coefficient vector for the relative thresholds \(10^{-k}\),
  \(k=2,3,4,5\) (from top left to bottom right). Notice that in case of $k=2$
  only the coarse level samplets remain after thresholding.}
\end{center}
\end{figure}

Figure~\ref{HM_fig:temp_HT_rec} shows the reconstruction \({\bs T}^\intercal
\big(\operatorname{HT}_{w}({\bs f}^\Sigma)\big)\) of the hard-thresholded
temperature for October 2022 with thresholds \(w=10^{-k}\|{\bs f}\|_{2}\) for
\(k=2,3,4,5\) and samplets with \(q+1=4\) vanishing moments. The threshold
values result in \(4\), \(87\), \(1\,454\), \(15\,688\) non-zero coefficients,
which amount to space savings of \(99.9996\%\), \(99.992\%\), \(99.86\%\), and
\(98.49\%\), respectively. The associated relative errors in the Euclidean norm
are \(2.4\%\), \(1.2\%\), \(0.39\%\), and \(0.11\%\).

\begin{figure}[htb]
\begin{center}
\begin{tikzpicture}
\draw(0,0)node{\includegraphics[scale=0.04,clip,trim=500 350 500 350]{
  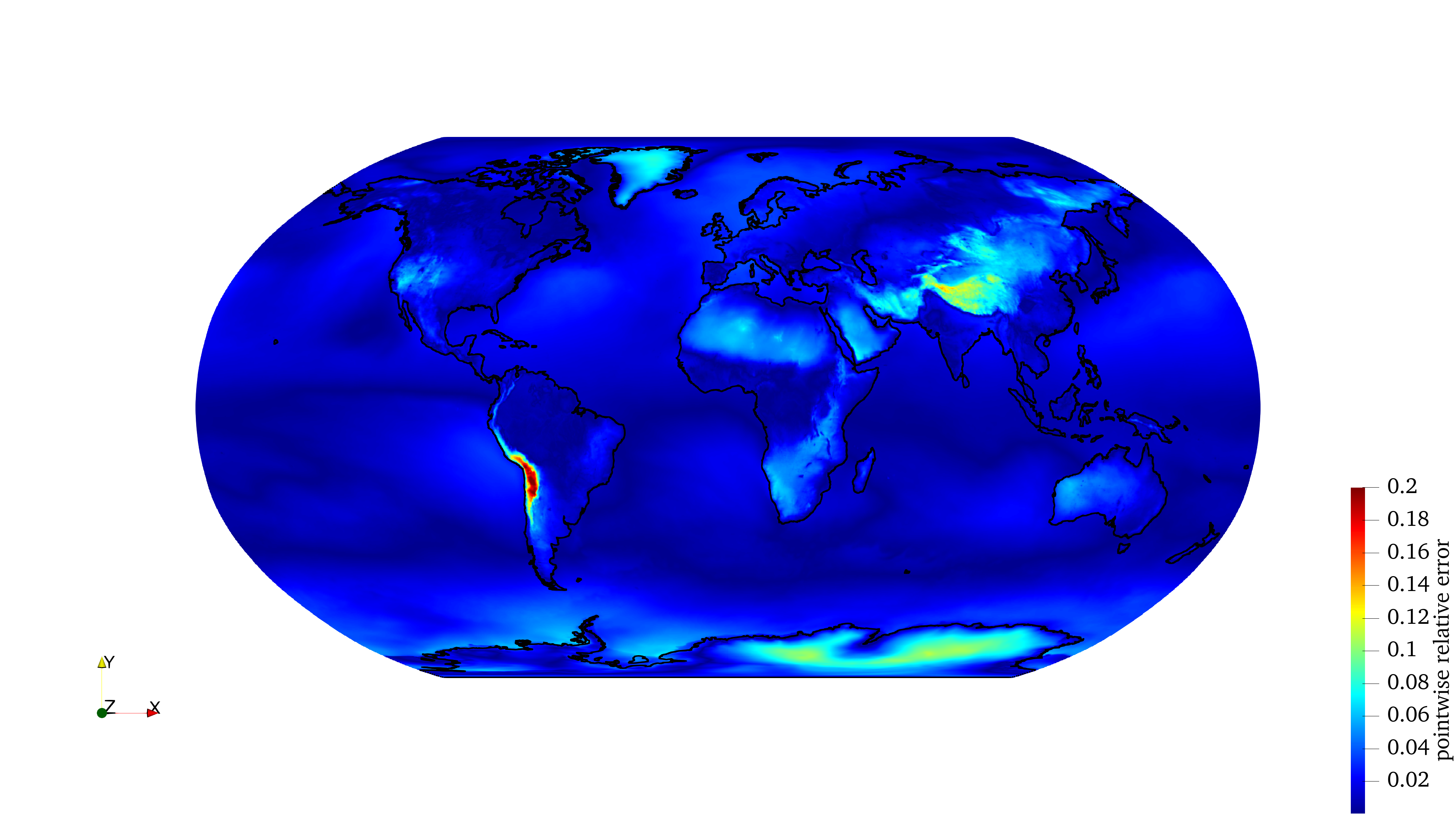}};
\draw(0,-1.25)node{\tiny Space saving: \(99.9996\%\)};
\draw(5,0)node{\includegraphics[scale=0.04,clip,trim=500 350 500 350]{
  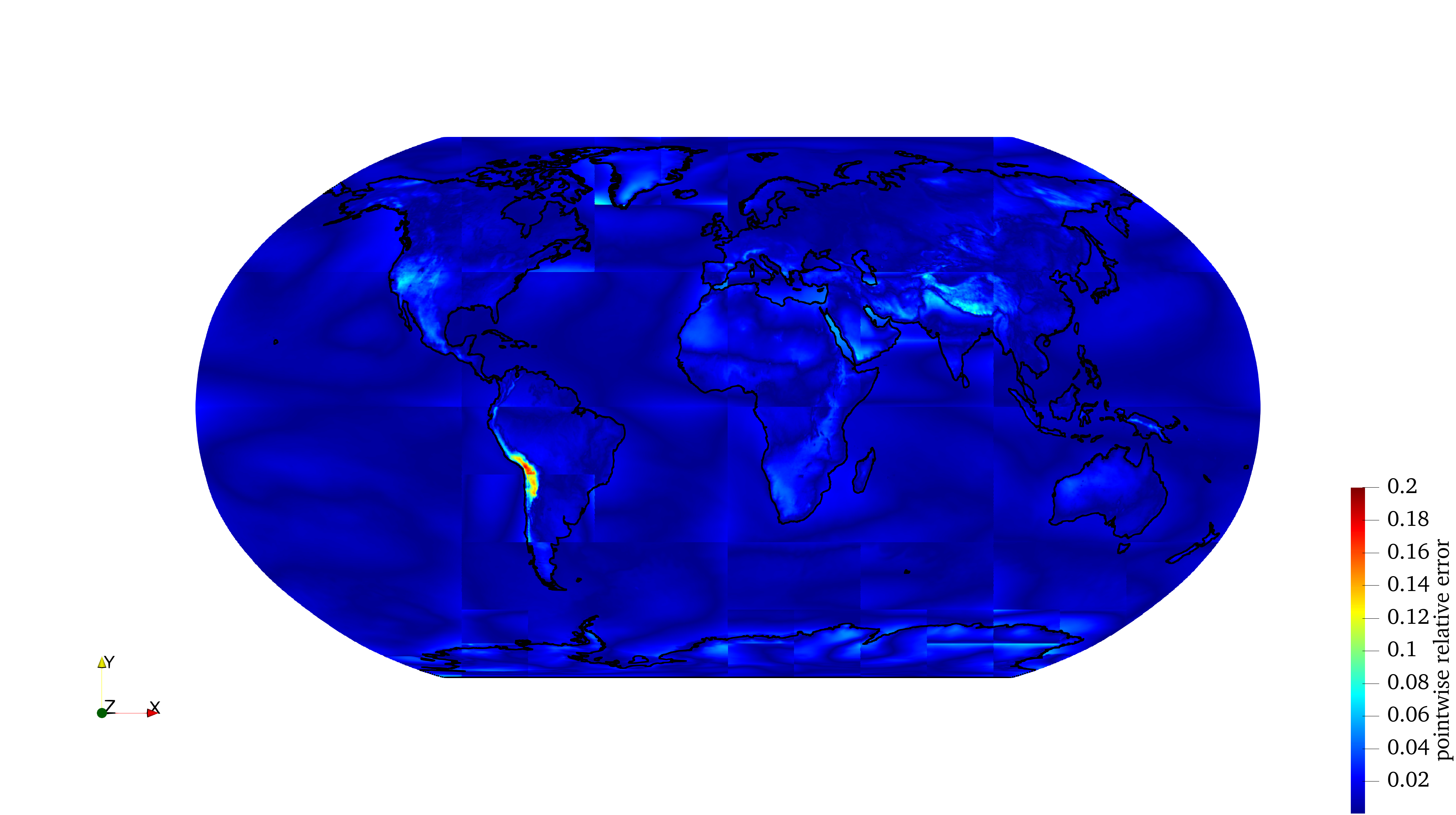}};
\draw(5,-1.25)node{\tiny Space saving: \(99.9916\%\)};
\draw(0,-3)node{\includegraphics[scale=0.04,clip,trim=500 350 500 350]{
  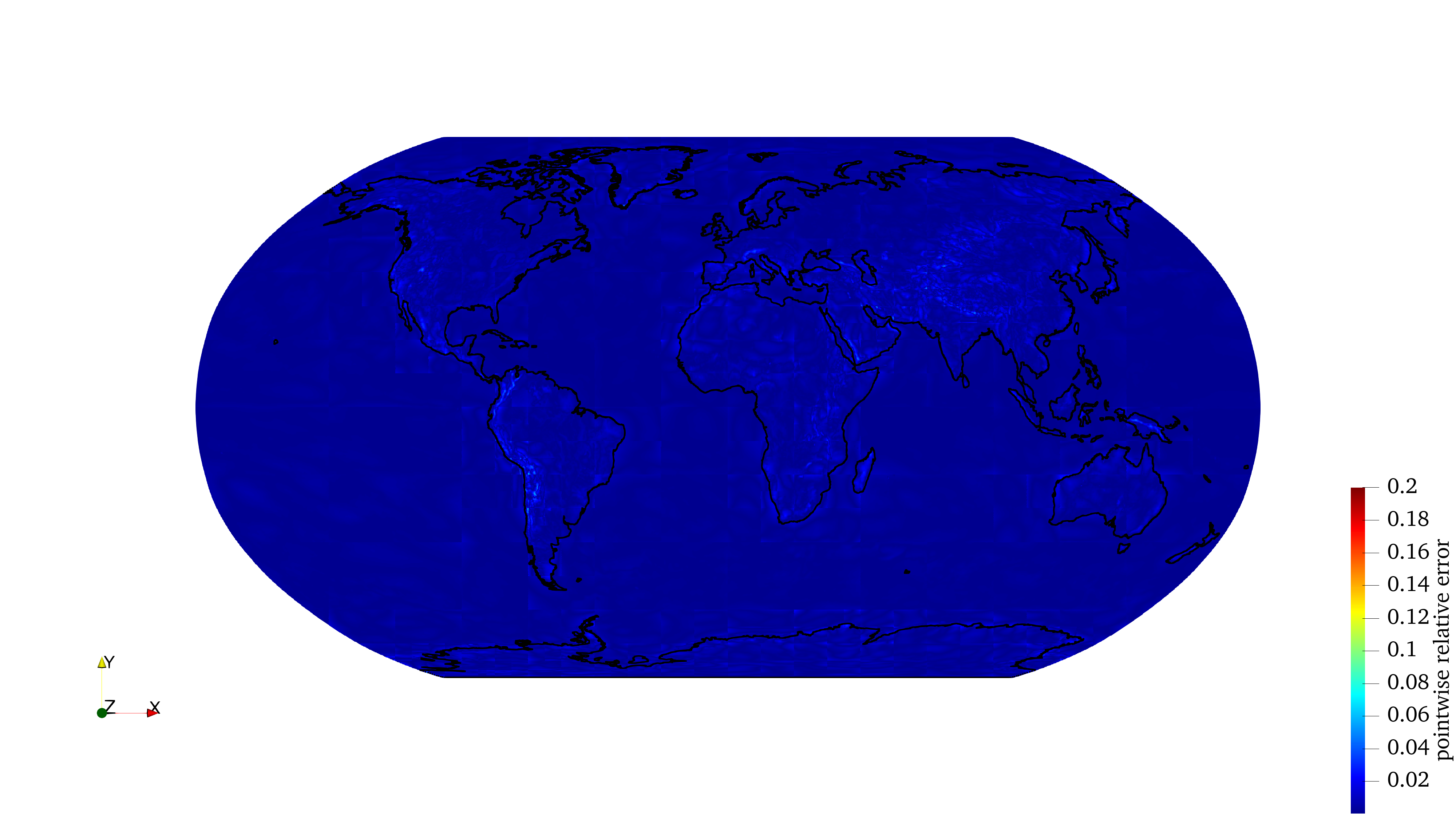}};
\draw(0,-4.25)node{\tiny Space saving: \(99.8600\%\)};
\draw(5,-3)node{\includegraphics[scale=0.04,clip,trim=500 350 500 350]{
  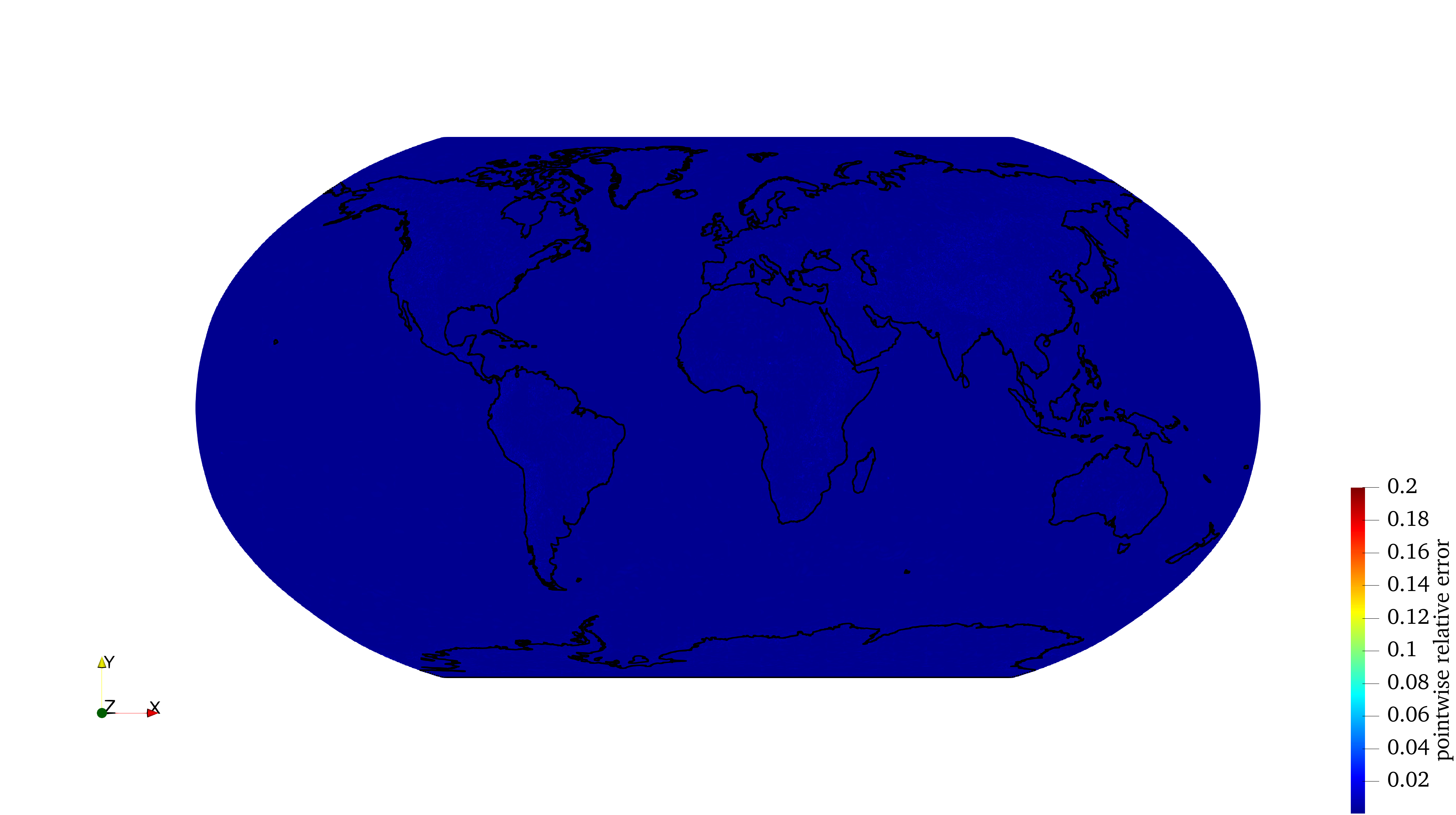}};
\draw(5,-4.25)node{\tiny Space saving: \(98.489\%\)};
\draw(8,-1.5)node{\includegraphics[scale=0.12,clip,trim=3620 0 0 1250]{
  ./Images/HT_err_1e-2_lins.png}};
\end{tikzpicture}
\caption{\label{HM_fig:temp_HT_err}Reconstruction error of the 
temperature with respect to the hard-thresholded coefficient 
vector for the relative thresholds \(10^{-k}\), \(k=2,3,4,5\)
(from top left to bottom right).}
\end{center}
\end{figure}

The associated pointwise relative errors of the reconstruction are visualized
in Figure~\ref{HM_fig:temp_HT_err}. For the largest threshold, the maximum of the
pointwise relative errors is \(21\%\), while it is \(2.2\%\) for the smallest
threshold.

%===============================================================================
\subsection{Adaptive subsampling}\label{HMsubsec:AdaptiveSubs}
%===============================================================================
Consider data sites \(X=\{{\bs x}_1,\ldots,{\bs x}_N\}\subset\Omega\) 
and corresponding data values \({\bs f}^\Delta=If\). The cluster tree for \(X\)\
is denoted by \(\Tcal\) as before. We propose an adaptive subsampling strategy 
based on the idea of the \emph{adaptive tree approximation} from \cite{HM_BD97}. 
Then, given a relevant subtree, we pursue an \emph{entropy based approach} for 
sampling uniformly from the leaves of that subtree. 

We start by generating a suitable subtree of \(\Tcal\) from the transformed data
values ${\bs f}^\Sigma = {\bs T}{\bs f}^\Delta$. To this end, we apply 
tree-coarsening as developed in \cite{HM_BD97} and directly adapt it to our setting.
The \emph{energy} contained in a cluster \(\tau\in\Tcal\) is defined as the sum
of energies of its children plus the squared Euclidean norm of the samplet 
coefficients belonging to \(\tau\), viz.,
\begin{equation}\label{HM_eg:energy}
%============================================
e(\tau)\isdef\sum_{\tau'\in\operatorname{child}(\tau)}e(\tau')
+\sum_{\sigma\in{\bs\Sigma}^{\tau}}(f,\sigma)_\Omega^2.
\end{equation}
Herein, we made the convention that \({\bs\Sigma}^X\) also contains the scaling
distributions the coarsest level. The quantity \(e(\tau)\) is the contribution
of the subtree with root \(\tau\) to the squared norm 
\(\big\|{\bs f}^\Sigma\big\|_2^2\). In particular, we have 
\(e(X)=\big\|{\bs f}^\Sigma\big\|_2^2\).

Based on the energies \eqref{HM_eg:energy}, we next define
\[
\tilde{e}(\tau')\isdef q(\tau)\isdef
\frac{\sum_{\mu\in\operatorname{child}(\tau)}e(\mu)}
{e(\tau)+\tilde{e}(\tau)}\tilde{e}(\tau)
\ \text{for all}\ \tau'\in\operatorname{child}(\tau), 
\]
where we set \(\tilde{e}(X)\isdef e(X)\) for the root of the cluster tree. Given
this modified energy, we perform the thresholding version of the second 
algorithm from \cite{HM_BD97} with threshold \(w = \varepsilon^2
\|{\bs f}^\Sigma\|_2^2\). This results in a subtree\(\Tcal_w\) that approximates
\({\bs f}^\Sigma\) up to a relative error of \(\varepsilon\) in the Euclidean
norm. Since the algorithm always selects either none or all children of a given
cluster, \(\Tcal_w\) is a cluster tree and its leaves \(\Lcal(\Tcal_w)\) form a
partition of \(X\).

\begin{figure}[htb]
\begin{center}
\begin{tikzpicture}
\draw(0,0)node{\includegraphics[scale=0.05,clip,trim=500 350 500 350]{
  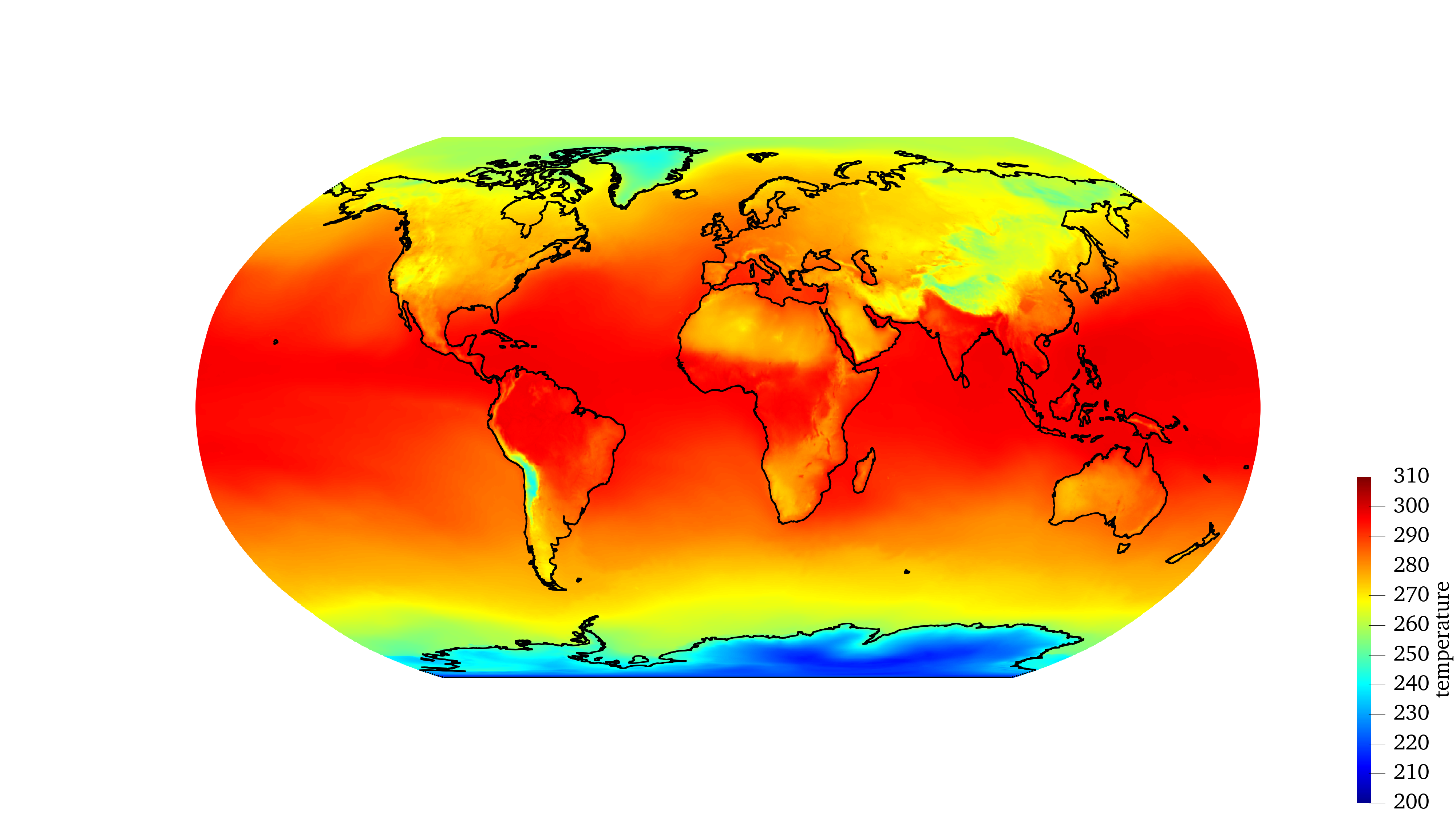}};
\draw(6,0)node{\includegraphics[scale=0.05,clip,trim=500 350 500 350]{
  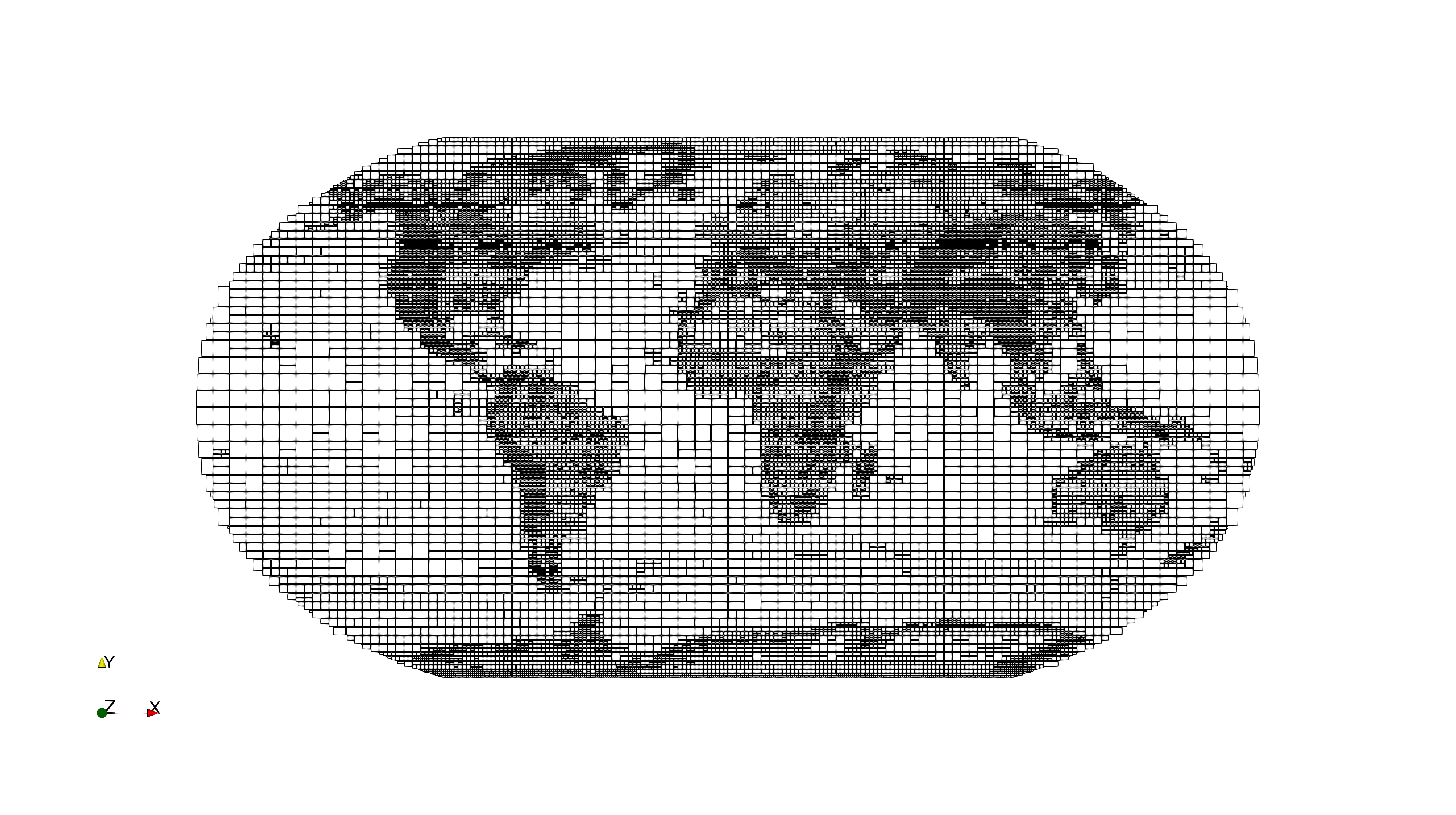}};
\end{tikzpicture}
\caption{\label{HM_fig:AdaptiveTree} 
Average temperature in October 2022 (left) and leafs of the
corresponding adaptive tree with relative threshold 
\(\varepsilon=10^{-4}\) (right).}
\end{center}
\end{figure}

Given the partition generated by \(\Lcal(\Tcal_w)\), we intend to generate a
subsample \((\bs x_{i_j},f_{i_j})\), \(j=1,\ldots,n\leq N\), adapted to the
features of the data. To this end, we pursue an \emph{entropy based approach}:
Let \(\mathbb{P}\) be a probability distribution on a discrete finite set \(Z\).
Then, the \emph{information content} of a sample \(z\in Z\) is defined as 
$I_{\mathbb{P}}(z)\isdef-\log\mathbb{P}(z)$. The \emph{entropy} $H_{\mathbb{P}}$
is the expected information content
\begin{equation}\label{HM_eg:entropy}
%==================================================
H_{\mathbb{P}}\isdef\Ebb_{\mathbb{P}}[I_{\mathbb{P}}]
    = -\sum_{z\in Z}\mathbb{P}(z)\log\mathbb{P}(z).
\end{equation}

It is well-known that the entropy \eqref{HM_eg:entropy} is maximized if each
element of \(Z\) is chosen with the same probability
\[
\mathbb{P}(z)=\frac{1}{|Z|},
\]
see the seminal work \cite{HM_Sha48} for all the details. Therefore, to generate
the adaptive subsample, we shall maximize the entropy for randomly selecting
points from leafs by choosing a point ${\bs x}\in X$ from a certain leaf 
\(\tau\in\Lcal(\Tcal_w)\) with probability
\begin{equation}\label{HM_eg:goal}
%===============================================================================
\mathbb{P}({\bs x}\in\tau)=\frac{1}{|\Lcal(\Tcal_w)|}.
\end{equation}

We like to emphasize that \eqref{HM_eg:goal} can also be used to adaptively
generate new samples ${\bs x}\in\Omega$ instead of considering only available
points ${\bs x}\in X$.

\begin{figure}[htb]
\begin{center}
\begin{tikzpicture}
\draw(6,0)node{\includegraphics[scale=0.05,clip,trim=500 350 500 350]{
  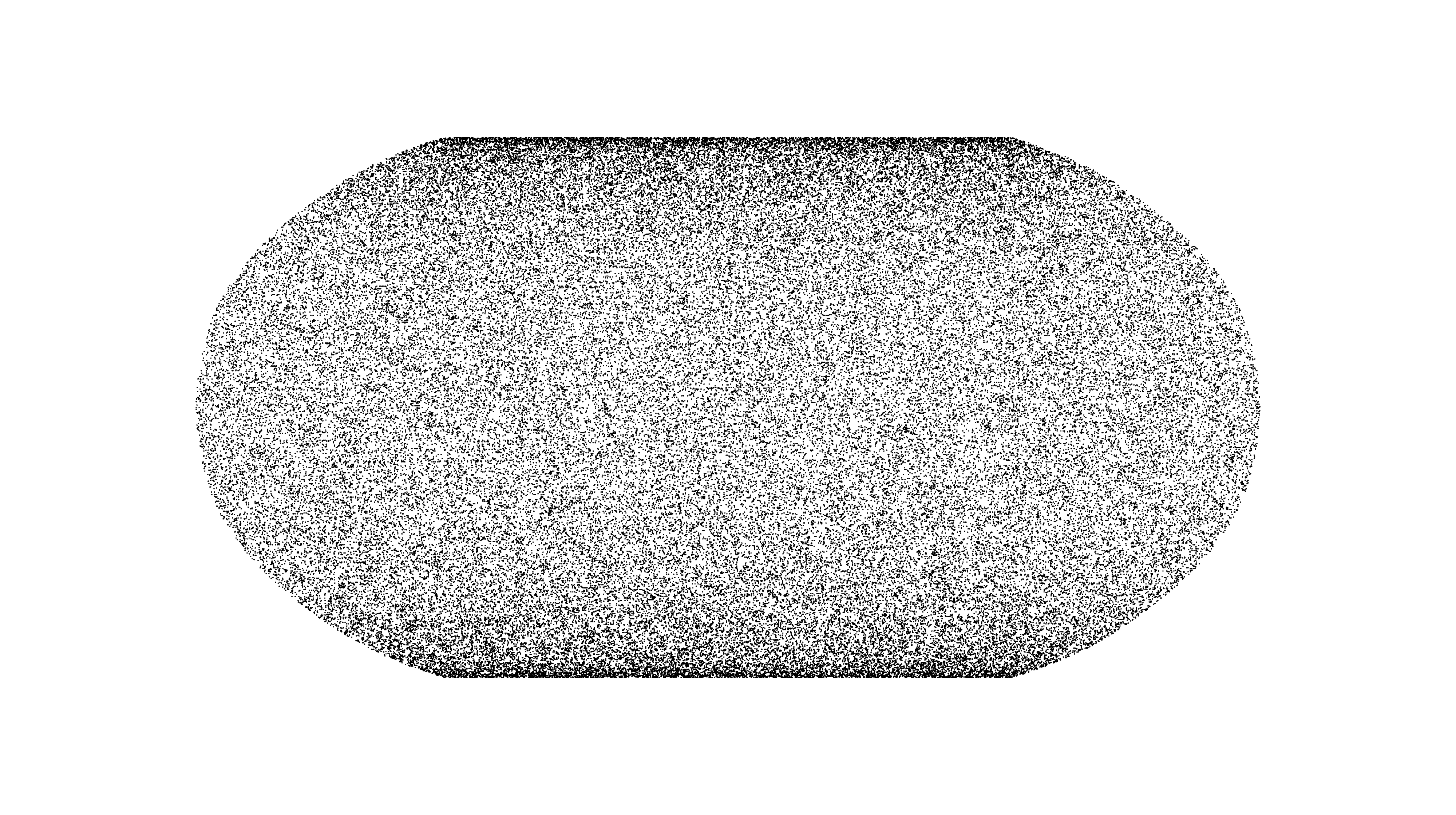}};
\draw(0,0)node{\includegraphics[scale=0.05,clip,trim=500 350 500 350]{
  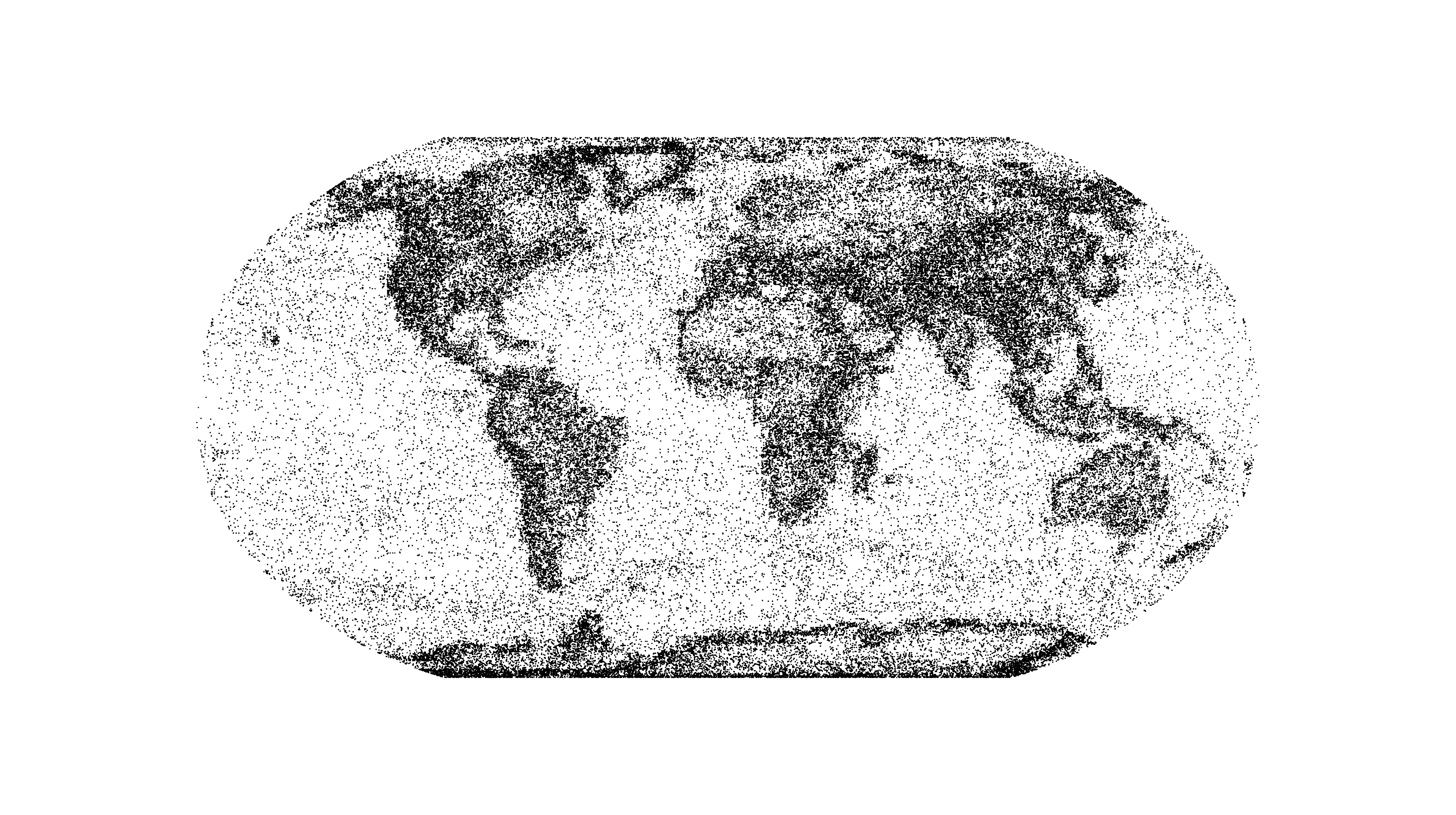}};
\end{tikzpicture}
\caption{\label{HM_fig:AdaptiveSample} 
Adaptive random subsample (left) and uniform random subsample (right).}
\end{center}
\end{figure}

To demonstrate the adaptive subsampling approach, we consider again the ERA5
temperature data set. Tree coarsening with relative threshold 
\(\varepsilon=10^{-4}\) for samplets with \(q+1=4\) vanishing moments yields
the tree found in the right plot of Figure~\ref{HM_fig:AdaptiveTree}. The
associated data are visualized in the accompanying left plot. The result of the
subsampling using 100\,000 points for the month October yields the point
distribution seen in the left plot of Figure \ref{HM_fig:AdaptiveSample}. In 
contrast, using the same amount of uniformly distributed samples yields the 
distribution in the right plot of Figure~\ref{HM_fig:AdaptiveSample}.

%===============================================================================
\section{Scattered data approximation}\label{HMsec:kernelCompression}
%===============================================================================
\subsection{Reproducing kernel Hilbert spaces}
%===============================================================================
Let \(\Hcal\subset C(\Omega)\) be a Hilbert space with inner product 
\(\langle\cdot,\cdot\rangle_{\Hcal}\). Then, there holds
\(\delta_{\bs x}\in\Hcal'\) for every \({\bs x}\in\Omega\).
Consequently, by the Riesz representation theorem, there exists 
\((R\delta_{\bs x})\in\Hcal\) such that 
\begin{equation}\label{HM_eg:RepProperty}
%====================================================
\langle (R\delta_{\bs x}), f\rangle_{\Hcal}=f({\bs x})
\ \text{for all}\ f\in\Hcal.
\end{equation}
The function $\kernel({\bs x},{\bs y})\isdef(R\delta_{\bs x})({\bs y})$ is the 
\emph{reproducing kernel} for \(\Hcal\), which renders it a \emph{reproducing 
kernel Hilbert space}, see \cite{HM_Aronszajn50,HM_Wendland2004} for instance. 
Especially, the reproducing kernel is symmetric and positive definite in the 
following sense.

\begin{definition}%\label{def:poskernel}
%=======================================================================
A symmetric kernel
$\kernel\colon\Omega\times\Omega\rightarrow\Rbb$ is 
\textit{positive definite} on $\Omega\subset\Rbb^d$,
iff the \emph{kernel matrix} 
\({\bs K}\isdef[\kernel({\bs x}_i,{\bs x}_j)]_{i,j=1}^N\)
is symmetric and positive semi-definite for every choice
of mutually distinct points
${\bs x}_1, \ldots,{\bs x}_N\in\Omega$
and any $N\in\mathbb{N}$. The kernel \(\kernel\) is
\emph{strictly positive definite}, iff \({\bs K}\)
is even positive definite.
\end{definition}

Associated to \(X=\{{\bm x}_1,\ldots,{\bm x}_N\}\subset\Omega\), we introduce
the subspace of kernel translates
\begin{equation}\label{HM_eg:HX}
%============================================
\Hcal_X\isdef
\spn\{\phi_1,\ldots,\phi_N\}
\subset\Hcal,
\end{equation}
where
\[
\phi_i\isdef\kernel({\bm x}_i,\cdot)\text{ for }
i=1,\ldots,N.
\]
The subspace \(\Hcal_X\) is isometrically isomorphic
to \(\Xcal=\spn\{\delta_{{\bm x}_1},\ldots,\delta_{{\bm x}_N}\}
\subset\Hcal'\) by the Riesz isometry. We identify
\[
u'=\sum_{i=1}^Nu_i\delta_{{\bm x}_i}\in\Xcal
\]
with
\[
Ru=\sum_{i=1}^Nu_i\phi_i\in\Hcal_X.
\]

To approximate a general element \(h\in\Hcal\) by an element in \(\Hcal_X\), we
compute its orthogonal projection. The latter is obtained by the solution of the
Galerkin formulation
\begin{equation}\label{HM_eg:Galerkin}
%===============================================================================
\langle s_h,v\rangle_\Hcal
=\langle h,v\rangle_\Hcal\ \text{for all}\ v\in\Hcal_X.
\end{equation}
Making the ansatz
\begin{equation}\label{HM_eg:RKHSinterp}
%===============================================
s_h=\sum_{i=1}^N\alpha_i\phi_i
\end{equation}
and choosing the basis of kernel translates as test functions in 
\eqref{HM_eg:Galerkin}, we arrive at the equivalent interpolation problem
$s_h({\bm x}_i) = h({\bm x}_i)$ for all ${\bm x}_i\in X$ due to the reproducing
property \eqref{HM_eg:RepProperty}. The expansion coefficients 
\({\bm\alpha}=[\alpha_i]_{i=1}^N\) from \eqref{HM_eg:RKHSinterp} can be retrieved
by solving the linear system of equations
\begin{equation}\label{HM_eg:LSE}
%===============================================================================
{\bm K}{\bm\alpha}={\bm h},\ \text{where}\ 
{\bm h}\isdef Ih=[(h,\delta_{{\bs x}_i})_\Omega]_{i=1}^N,
\end{equation}
Here, \({\bm K}=[\kernel({\bm x}_i,{\bm x}_j)]_{i,j=1}^N\) is the kernel matrix
and $I$ is the information operator from \eqref{HM_eg:InfOp}. It is known from
both, Galerkin theory and optimal recovery, that \(s_h\) is a minimum norm 
solution in \(\Hcal\), see \cite{HM_Bra13} and \cite{HM_Mic84}, respectively.

Depending on the choice of the kernel function, the linear system 
\eqref{HM_eg:LSE} of equations can be ill conditioned and a suitable 
regularization is required to obtain a solution. Traditionally, Tikhinov 
regularization is used and the system
\begin{equation}\label{HM_eg:LSEreg}
%===============================================================================
({\bm K}+\mu{\bs I}){\bm\alpha}={\bm h}
\end{equation}
is solved for an appropriate regularization parameter $\mu>0$. In Section 
\ref{HMsec:sparsity}, we address also \(\ell^1\)-regularization with respect to
the samplet basis, which is known to impose sparsity to the solution.

%===============================================================================
\subsection{Samplets in reproducing kernel Hilbert spaces}
%===============================================================================
Employing the Riesz isometry, we can embed a given samplet basis 
\({\bs\Sigma}=\{\sigma_{j,k}\}_{j,k}\) into a reproducing kernel Hilbert space. 
This idea follows the vast literature on the embedding of empirical 
distributions into reproducing kernel Hilbert spaces, see \cite{HM_MFSS17} for 
example. 

Consider a samplet 
\[
\sigma_{j,k}=\sum_{i=1}^N \omega_{j,k,i}\delta_{{\bs x}_{i}}\subset\Hcal'
\]
where \(\omega_{j,k,i}\), \(i=1,\ldots,N\), are the expansion coefficients of 
\(\sigma_{j,k}\) with respect to the Dirac-$\delta$-distributions in \(\Xcal\).
The samplet can be identified with the function
\[
\psi_{j,k}\isdef
\sum_{i=1}^N\omega_{j,k,i}\phi_i\in\mathcal{H}
\]
by means of the Riesz isometry. The vanishing moment property 
\eqref{HM_eg:vanishingMoments} translates to
\[
\langle\psi_{j,k},h\rangle_\Hcal=0
\]
for any \(h\in\Hcal\) which satisfies \(h|_{O}\in\Pcal_q\) for any open and 
convex set \(O\) with \(\supp(\sigma_{j,k})\subset O\cap\Omega\). Herein, we
define the support of a samplet in the context of the support of measures
according to
\[
\supp(\sigma_{j,k})\isdef\{{\bs x}_i\in X : \omega_{j,k,i}\neq 0\}.
\]

The functions \(\{\psi_{j,k}\}_{j,k}\) span the subspace \(\Hcal_X\), see 
\eqref{HM_eg:HX}. In particular, defining 
\[
\Wcal_j\isdef\operatorname{span}\{\psi_{j,k}\}_k,
\]
we obtain the primal multiresolution analysis
\[
\Hcal_X=\bigoplus_{j}\Wcal_j.
\]
Using the embedded samplets \(\psi_{j,k}\) as ansatz- and test functions in the
Galerkin formulation \eqref{HM_eg:Galerkin} yields the linear system of equations
\begin{equation}\label{HM_eg:Tlinsys}
%===============================================================================
{\bs K}^\Sigma{\bs\beta} = {\bs h}^\Sigma,
\end{equation}
where
\begin{equation}\label{HM_eg:samplet-matrix}
%===============================================================================
\begin{aligned}
{\bs K}^\Sigma &=
\big[(\kernel,\sigma_{j,k}\otimes\sigma_{j',k'})_{\Omega\times\Omega}
\big]_{j,j',k,k'}\\
&= [\langle\psi_{j,k},\psi_{j',k'}\rangle_\Hcal]_{j,j',k,k'}\\
&= {\bs T}{\bs K}{\bs T}^\intercal
\end{aligned}
\end{equation}
and
\[
{\bs h}^\Sigma 
= \big[(\sigma_{j,k},h)_{\Omega}]_{j,k} 
= [\langle\psi_{j,k},h\rangle_\Hcal]_{j,k} 
= {\bs T}{\bs h}.
\]
The solution ${\bs\beta}$ of the linear system \eqref{HM_eg:Tlinsys} of equations
is equivalent to the one of \eqref{HM_eg:LSE} by the transform
\[
  {\bs\beta}={\bs T}{\bs\alpha} = {\bs T}{\bs K}^{-1}{\bs h}.
\]

Noticing that the Gramian satisfies
\[
[\langle\kernel({\bs x}_i,\cdot),\kernel({\bs x}_j,\cdot)
\rangle_\Hcal]_{i,j=1}^N
=[\kernel({\bs x}_i,{\bs x}_j)]_{i,j=1}^N={\bs K},
\]
we find
\[
\langle\psi_{j,k},\psi_{j',k'}\rangle_\Hcal
={\bs\omega}_{j,k}^\intercal
{\bs K}{\bs\omega}_{j',k'},
\]
where we set \({\bs\omega}_{j,k}\isdef[\omega_{j,k,i}]_{i=1}^N\)
and \({\bs\omega}_{j',k'}\isdef[\omega_{j',k',i}]_{i=1}^N\).
Hence, the \emph{dual basis} is given by
\[
\widetilde{\psi}_{j,k}=\sum_{i=1}^N\widetilde{s}_{j,k,i}\phi_i,
\ \text{where}\  
\widetilde{\bs\omega}_{j,k}\isdef{\bs K}^{-1}{\bs\omega}_{j,k}.
\]
Defining the spaces
\[
\widetilde{\Wcal}_j\isdef\operatorname{span}\{\widetilde{\psi}_{j,k}\}_k
\]
yields the dual multiresolution analysis
\[
\Hcal_X=\bigoplus_{j}\widetilde{\Wcal}_j,
\]
where $\Wcal_j\perp\widetilde{\Wcal}_{j'}$ for $j\neq j'$
since $\big\langle{\psi}_{j,k},\widetilde{\psi}_{j',k'}
\big\rangle_\Hcal =\delta_{j,j'}\delta_{k,k'}$. 

With respect to the dual basis, the interpolant 
\eqref{HM_eg:RKHSinterp} can be written as 
\begin{equation}\label{HM_eg:dualBasisExpansion}
%===============================================================================
s_h=\sum_{j,k}\beta_{j,k}\psi_{j,k}=\sum_{j,k}
\big\langle\widetilde{\psi}_{j,k},h\big\rangle_\Hcal\psi_{j,k}=
\sum_{j,k}
\big\langle{\psi}_{j,k},h\big\rangle_\Hcal\widetilde{\psi}_{j,k}.
\end{equation}
In view of the multiresolution representation of the information operator in 
\(Ih=[(h,\sigma_{j,k})_\Omega]_{j,k}=[\langle h,\psi_{j,k}\rangle_\Hcal]_{j,k}\),
compare \eqref{HM_eg:InfOp}, the expansion with respect to the dual basis in
\eqref{HM_eg:dualBasisExpansion} amounts to the multiresolution variant of the 
optimal recovery algorithm from \cite{HM_MR77}.

\begin{figure}[htb]
\begin{center}
\scalebox{0.52}{
\begin{tikzpicture}
\begin{axis}[width=7.8cm, height=6cm, xmin = -1, xmax=1,
 ymin=-0.5, ymax=1.5, ylabel={}, xlabel ={$j=0$}, ytick={-0.5,0,0.5}]
\addplot[color=black, line width=1pt]table[x index={0},y index = {1}]{%
./Images/esamplets.txt};
\addlegendentry{$\psi$};
\addplot[color=gray, densely dotted, line width=1pt]
  table[x index={0},y index = {1}]{%
./Images/etsamplets.txt};
\addlegendentry{$\widetilde{\psi}$};
\addplot[color=black, mark=o, line width=0.4pt,mark size=1.2pt, only marks]
  table[x index={0},y index = {1}]{%
./Images/dpts.txt};
\addlegendentry{$x_i$};
\end{axis}
\end{tikzpicture}}
\scalebox{0.52}{
\begin{tikzpicture}
\begin{axis}[width=7.8cm, height=6cm, xmin = -1, xmax=1,
 ymin=-0.25, ymax=0.25, ylabel={}, xlabel ={$j=1$},ytick={-0.2,0,0.2}]
\addplot[color=black, line width=1pt]table[x index={0},y index = {2}]{%
./Images/esamplets.txt};
\addlegendentry{$\psi$};
\addplot[color=gray, densely dotted, line width=1pt]
  table[x index={0},y index = {2}]{%
./Images/etsamplets.txt};
\addlegendentry{$\widetilde{\psi}$};
\addplot[color=black, mark=o, line width=0.4pt,mark size=1.2pt, only marks]
  table[x index={0},y index = {1}]{%
./Images/dpts.txt};
\addlegendentry{$x_i$};
\end{axis}
\end{tikzpicture}}
\scalebox{0.52}{
\begin{tikzpicture}
\begin{axis}[width=7.8cm, height=6cm, xmin = -1, xmax=1,
 ymin=-0.25, ymax=0.25, ylabel={}, xlabel ={$j=2$},ytick={-0.2,0,0.2}]
\addplot[color=black, line width=1pt]table[x index={0},y index = {3}]{%
./Images/esamplets.txt};
\addlegendentry{$\psi$};
\addplot[color=gray, densely dotted, line width=1pt]
  table[x index={0},y index = {3}]{%
./Images/etsamplets.txt};
\addlegendentry{$\widetilde{\psi}$};
\addplot[color=black, mark=o, line width=0.4pt,mark size=1.2pt, only marks]
  table[x index={0},y index = {1}]{%
./Images/dpts.txt};
\addlegendentry{$x_i$};
\end{axis}
\end{tikzpicture}}
\caption{\label{HM_fig:sampletViz}\(H^1(\Rbb)\)-embedded primal and 
dual scaling distribution (left), samplet on level \(j=1\) (middle)
and samplet on level \(j=2\) (right) for \(N=200\) uniformly 
distributed data sites and \(q+1=3\).}
\end{center}
\end{figure}
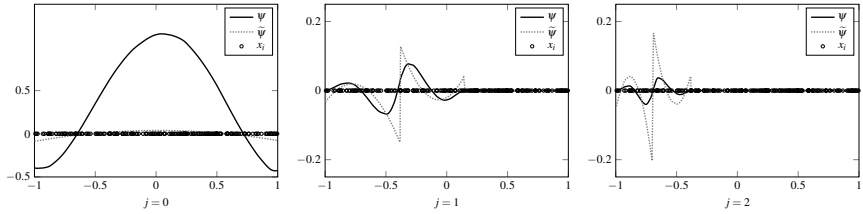

In order to give a visual idea of such \emph{embedded samplets} and the 
respective dual basis, we consider the Sobolev space \(H^1(\Rbb)\), equipped 
with the usual norm 
\[
\|u\|_{H^1(\Rbb)}^2 
= \|u\|_{L^2(\Rbb)}^2 + \|u'\|_{L^2(\Rbb)}^2.
\] 
Its reproducing kernel is given by \(\kernel(s,t) = \frac{1}{2}\exp(|s-t|)\),
see \cite{HM_BA04} for instance. Figure~\ref{HM_fig:sampletViz} shows an embedded
scaling distribution (left plot) and two embedded samplets (middle and right
plots) with \(q+1=3\) vanishing moments, constructed for \(N=200\) uniformly
distributed data sites.

%===============================================================================
\subsection{Asymptotically smooth kernels}
%===============================================================================
It is well-known that nonlocal operators of Calder\'on-Zygmund type can be 
compressed in wavelet coordinates, see \cite{HM_BCR,HM_DHS,HM_SCHN,HM_PS,HM_PSS97} for example.
Analogous results are available for samplets in the context of kernel matrices
of positive definite kernels.

The essential ingredient for the samplet compression of kernel matrices is the
\emph{asymptotical smoothness} property of the kernel \(\kernel\). This means 
that there exists a $\rho>0$ such that for all $({\bs x},
{\bs y})\in(\Omega\times\Omega)\setminus\Delta$ 
\begin{equation}\label{HM_eg:kernel_estimate}
%===============================================================================
  \bigg|\frac{\partial^{|\bs\alpha|+|\bs\beta|}}
  	{\partial{\bs x}^{\bs\alpha}
  	\partial{\bs y}^{\bs\beta}} \kernel({\bs x},{\bs y})\bigg|
  		\lesssim \frac{(|\bs\alpha|+|\bs\beta|)!}
		{\rho^{|\bs\alpha|+|\bs\beta|}
		\|{\bs x}-{\bs y}\|_2^{|\bs\alpha|+|\bs\beta|}}
\end{equation}
uniformly in $\bs\alpha,\bs\beta\in\mathbb{N}^d$ apart 
from the diagonal $\Delta\isdef \{({\bs x},{\bs y})\in\Omega
\times\Omega:{\bs x} = {\bs y}\}$.

A particular class of positive definite kernels which are asymptotically smooth
are the \emph{Mat\'ern kernels}
\(\kernel({\bs x},{\bs y})\isdef k_\nu(\|{\bs x}-{\bs y}\|_2)\).
They are known to be the reproducing kernels of the Sobolev spaces, 
see~\cite{HM_Wendland2004} for example, and are given by
\begin{equation}\label{HM_eg:MaternKernels}
%===============================================================================
k_\nu(r)\isdef\frac{2^{1-\nu}}{\Gamma(\nu)}\bigg(\frac {\sqrt{2\nu}r}
{\ell}\bigg)^\nu
K_\nu\bigg(\frac {\sqrt{2\nu}r}{\ell}\bigg),
\end{equation}
where $r=\|{\bs x}-{\bs y}\|_2\geq 0$ is the Euclidean distance and $\ell>0$ is
the lengthscale parameter. Herein, $K_{\nu}$ is the modified Bessel function of
the second kind of order $\nu$ and $\Gamma$ is the gamma function. The parameter
$\nu$ controls the smoothness of the kernel function, see, for example,
\cite{HM_Rasmussen2006}. In particular, we have the exponential kernel 
\begin{equation}\label{HM_eg:1/2}
%===============================================================================
k_{1/2}(r)=\exp\bigg(-\frac{r}{\ell}\bigg)
\end{equation}
and the Gaussian kernel
\[
k_{\infty}(r)=\exp\bigg(-\frac{r^2}{2\ell^2}\bigg).
\]

For asymptotically smooth kernels, compare \eqref{HM_eg:kernel_estimate}, we
obtain the following result, which is the basis for the matrix compression
introduced thereafter. The proof of this results is obtained in the traditional
way by applying Taylor's expansion of the kernel function under consideration,
using \eqref{HM_eg:kernel_estimate} to bound the remainder terms, compare 
\cite{HM_HM2}.

\begin{lemma}\label{lem:kernel_decay}
%===============================================================================
Consider two samplets $\sigma_{j,k}$ and $\sigma_{j',k'}$ which exhibit $q+1$ 
vanishing moments and let the associated clusters \(\tau\) and \(\tau'\) be such
that $\dist(\tau,\tau') > 0$. Then, for kernels satisfying 
\eqref{HM_eg:kernel_estimate}, there holds that
\begin{equation}\label{HM_eg:kernel_decay}
%===============================================================================
  (\kernel,\sigma_{j,k}\otimes\sigma_{j',k'})_{\Omega\times\Omega}\lesssim
  	\sqrt{|\tau||\tau'|}\frac{\diam(\tau)^{q+1}\diam(\tau')^{q+1}}
		{(\rho\dist(\tau,\tau')/d)^{2(q+1)}}.
\end{equation}
\end{lemma}

%===============================================================================
\subsection{Matrix compression}
\label{HMsec:matrix_compression}
%===============================================================================
The decay estimate \eqref{HM_eg:kernel_decay} immediately yields a compression
strategy for kernel matrices in samplet representation, the proof of which is 
found in \cite{HM_HM2}. This compression differs from the wavelet matrix 
compression introduced in \cite{HM_DHS}, since we do not exploit the decay of the 
samplet coefficients with respect to the level in case of smooth data. This 
enables us to also cover scattered data sets with arbitrarily distributed
points. As a result, we use the same accuracy on all levels, which is similar to
the setting in \cite{HM_BCR}, but implies that the number $q+1$ of vanishing 
moments needs to be increased to arrive at a higher accuracy of the matrix 
compression. Fixing the accuracy seems however not to be an issue in view of the
regularization of the kernel matrices in scattered data approximation.

\begin{theorem}\label{thm:compression}
Set all coefficients of the kernel matrix ${\bs K}^\Sigma$ from 
\eqref{HM_eg:samplet-matrix} to zero which satisfy the admissibility condition
\begin{equation}\label{HM_eg:cutoff}
%===============================================================================
   \dist(\tau,\tau')\ge\eta\max\{\diam(\tau),\diam(\tau')\},\quad\eta>0,
\end{equation}
where \(\tau\) is the cluster supporting \(\sigma_{j,k}\) and \(\tau'\) is the 
cluster supporting \(\sigma_{j',k'}\), respectively. Then, there holds
\[
  \big\|{\bs K}^\Sigma-{\bs K}^\Sigma_\varepsilon\big\|_F
	\lesssim m_q (\rho\eta/d)^{-2(q+1)} N \log(N),
\]
where \(m_q\) is given by \eqref{HM_eg:mq}.
\end{theorem}

\begin{remark}
In case of quasi-uniform points ${\bs x}_i\in X$, we have 
$\big\|{\bs K}^\Sigma\big\|_F\sim N$. Moreover, the term $\log(N)$ can be 
removed by a refined analysis in this case and we arrive at
\[
  \frac{\big\|{\bs K}^\Sigma-{\bs K}^\Sigma_\varepsilon\big\|_F}
  {\big\|{\bs K}^\Sigma\big\|_F} \lesssim  m_q (\rho\eta/d)^{-2(q+1)}
\]
while the compressed matrix has $\mathcal{O}(m_q^2N\log N)$ nonzero
coefficients.
\end{remark}

\begin{figure}[htb]
\begin{center}
\begin{tikzpicture}
\draw(0,0)node{\includegraphics[scale=0.26,clip,trim= 0 0 0 12.25,frame]{
  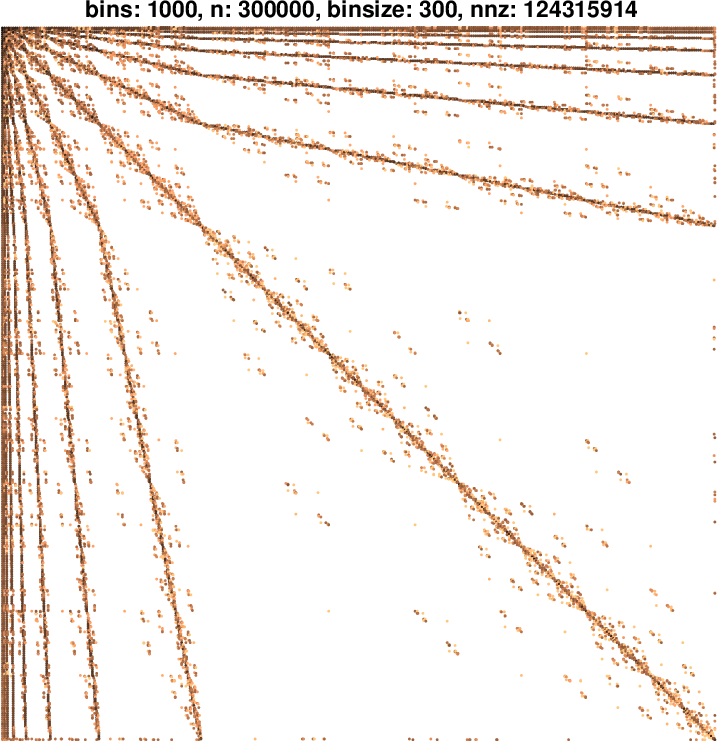}};
\draw(4,0)node{\includegraphics[scale=0.26,clip,trim= 0 0 0 12.25,frame]{
  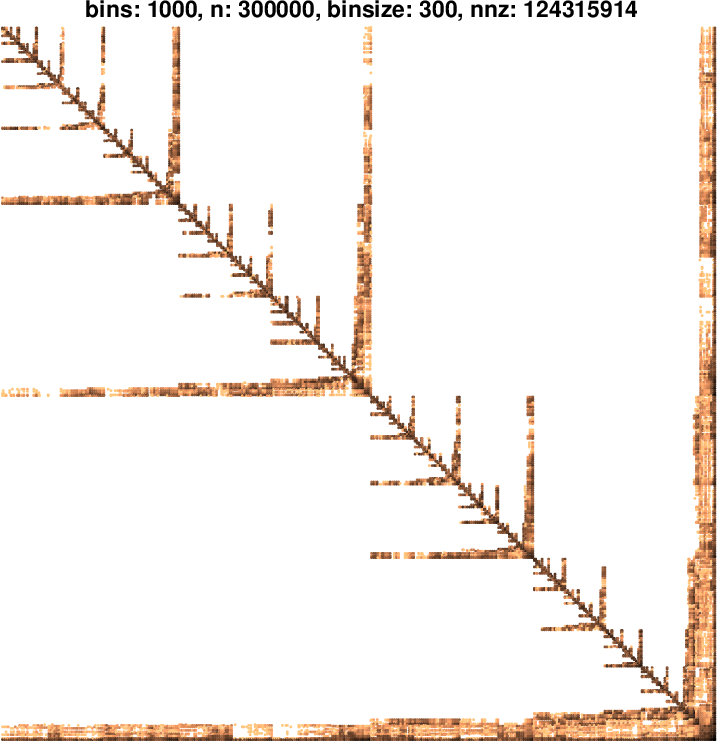}};
\draw(8,0)node{\includegraphics[scale=0.26,clip,trim= 0 0 0 12.25,frame]{
  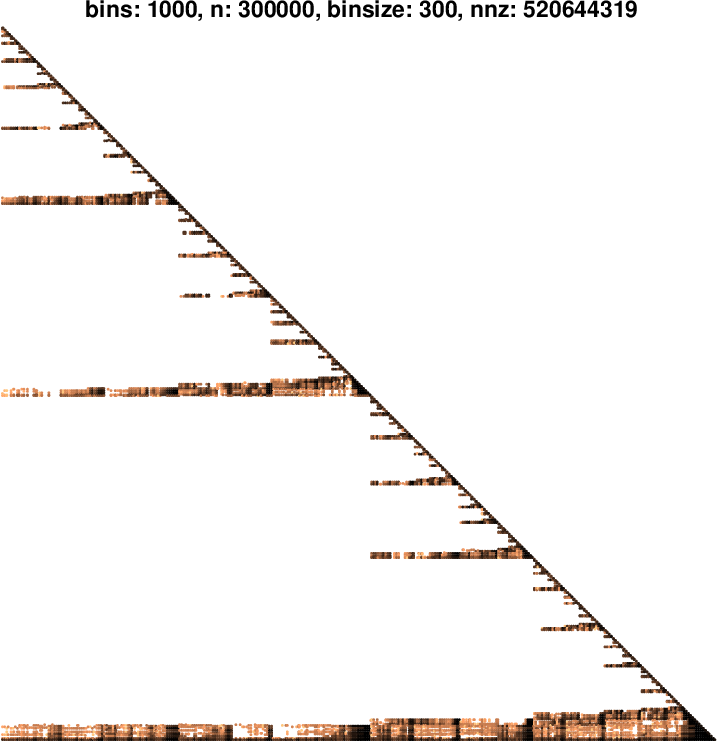}};
\end{tikzpicture}
\caption{\label{HM_fig:KernelArithmetics}Sparsity patterns of the
samplet compressed exponential kernel on the unit square (left), 
the nested dissection reordering (middle), and the Cholesky factor 
(right). Each dot represents a matrix block of size \(300\times 300\). 
The entries per block are color coded, where lighter blocks have less entries.
}
\end{center}
\end{figure}

An illustration of the matrix pattern in case of the exponential kernel 
$k_{1/2}$ from \eqref{HM_eg:1/2} with lengthscale parameter \(\ell=0.1\) on the
unit square using quasi-uniform samplets with \(q+1=4\) vanishing moments is
found in the left plot of Figure \ref{HM_fig:KernelArithmetics}. Especially, it
has been observed in \cite{HM_HM1} that the compressed matrix can be reordered by
means of nested dissection \cite{HM_Geo73} such that a very sparse Cholesky
decomposition is obtained. This is illustrated in the middle and 
right plots of Figure~\ref{HM_fig:KernelArithmetics}.

%===============================================================================
\subsection{Matrix assembly}
\label{HMsec:matrix_assembly}
%===============================================================================
Using the compression rule \eqref{HM_eg:cutoff}, we can now determine for a given
pair of clusters whether the corresponding entries need to be calculated. As
there are $\mathcal{O}(N)$ clusters, naively checking the cut-off criterion for
all pairs would still take $\mathcal{O}(N^{2})$ operations. Hence, we require
smarter means to determine the non-negligible cluster pairs. For this purpose,
we first state the transitivity of the admissibility condition to child 
clusters, compare \cite{HM_DHS} for a proof.

\begin{lemma}
Let $\tau$ and $\tau'$ be clusters satisfying the admissibility condition
\eqref{HM_eg:cutoff}. Then, for the child clusters $\tau_{\mathrm{child}}$ 
of $\tau$ and $\tau_{\mathrm{child}}'$ of $\tau'$, we have 
\begin{align*}
\dist(\tau,\tau_{\mathrm{child}}')
&\ge\eta\max\{\diam(\tau),\diam(\tau_{\mathrm{child}}')\},\\
\dist(\tau_{\mathrm{child}},\tau')
&\ge\eta\max\{\diam(\tau_{\mathrm{child}}),\diam(\tau')\},\\
\dist(\tau_{\mathrm{child}},\tau_{\mathrm{child}}')
&\ge\eta\max\{\diam(\tau_{\mathrm{child}}),\diam(\tau_{\mathrm{child}}')\}.
\end{align*}
\end{lemma}

The lemma tells us that we may omit cluster pairs whose parent clusters already
satisfy the admissibility condition. This is essential for the efficient
assembly of the compressed kernel matrix by means of \(\Hcal^2\)-matrix
techniques, see \cite{HM_HB02,HM_Gie01}. This idea has already been proposed earlier
in \cite{HM_AHK14,HM_HKS05,HM_Kae07} in case of Tausch-White wavelets.

\(\Hcal^2\)-matrices approximate the kernel interaction for sufficiently distant
clusters \(\tau\) and \(\tau'\) in the sense of the admissibility condition
\eqref{HM_eg:cutoff} by means of a polynomial interpolant, see \cite{HM_Boe10}. More
precisely, given a suitable set of interpolation points 
\(\{{\bs\xi}_t^\tau\}_t\) for each cluster \(\tau\) with associated Lagrange
polynomials \(\{\mathcal{L}_{t}^{\tau}({\bs x})\}_t\), we can approximate 
an admissible matrix block by interpolation:
\[
[(\kernel,\delta_{\bs x}\otimes
\delta_{\bs y})_{\Omega\times\Omega}]_{{\bs x}\in\tau,{\bs y}\in\tau'}
\approx\sum_{s,t} \kernel({\bs\xi}_{s}^{\tau}, {\bs\xi}_{t}^{\tau'}) 
		[(\mathcal{L}_{s}^{\tau},\delta_{\bs x})_\Omega]_{{\bs x}\in\tau}
		[(\mathcal{L}_{t}^{\tau'},\delta_{\bs y})_\Omega]_{{\bs y}\in\tau'}.
\]
However, different from the \(\Hcal^2\)-matrix setting, we shall consider this
expansion also when the clusters \(\tau\) and \(\tau'\) are located on different
levels of the cluster tree. By running recursively in a clever way through the 
samplet matrix, we arrive at an algorithm scheme that computes the compressed
kernel matrix in loglinear overall cost. We skip further details here and refer
the reader to \cite{HM_HM2} instead.

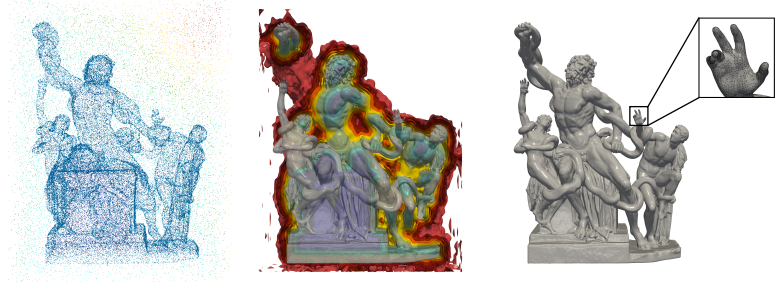
\begin{figure}[htb]
\begin{center}
\begin{tikzpicture}
\draw(-3.2,0)node{\includegraphics[scale=0.042,clip,trim= 1000 100 1000 100]{
  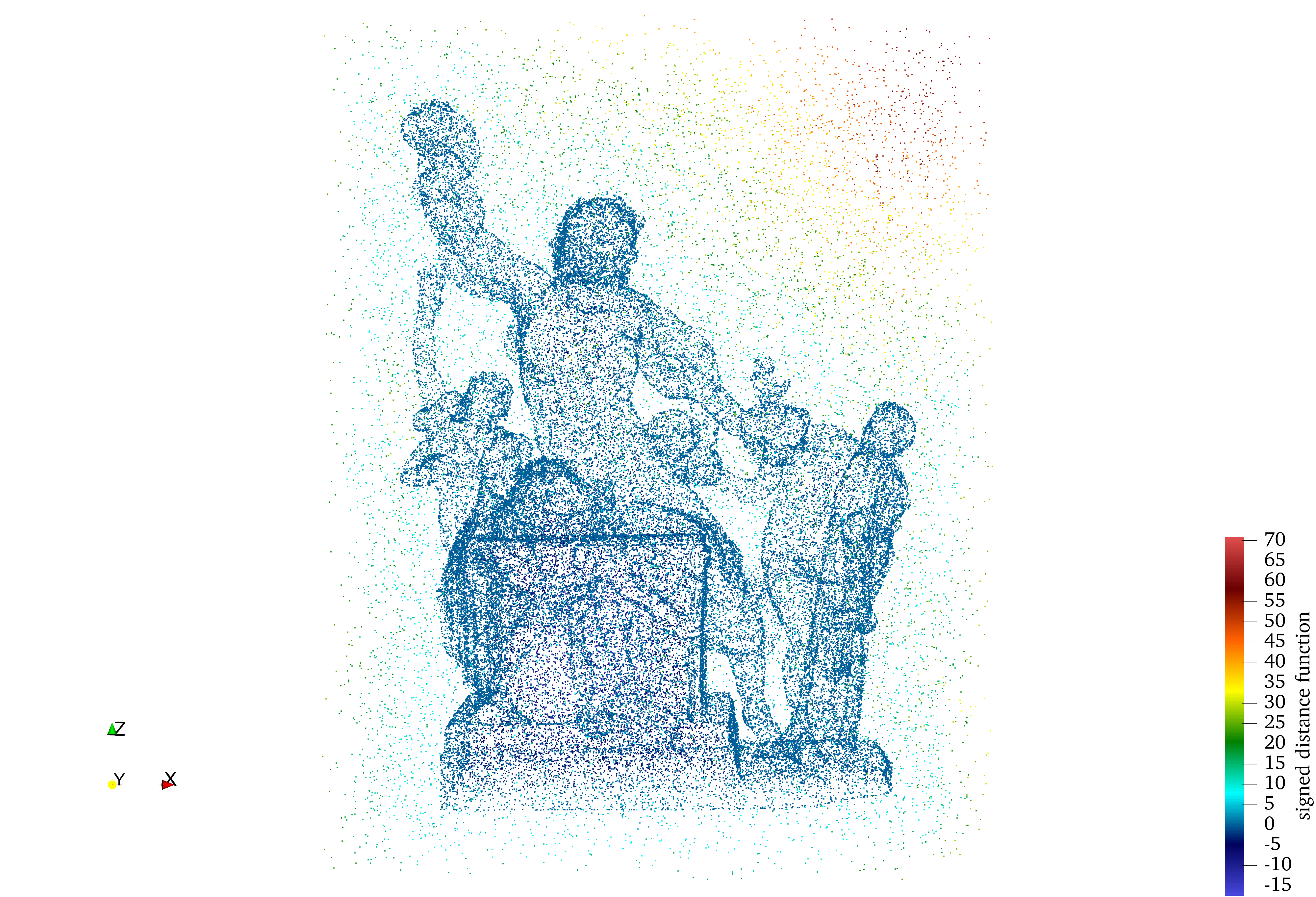}};
\draw(0,0)node{\includegraphics[scale=0.042,clip,trim= 1000 100 1000 100]{
  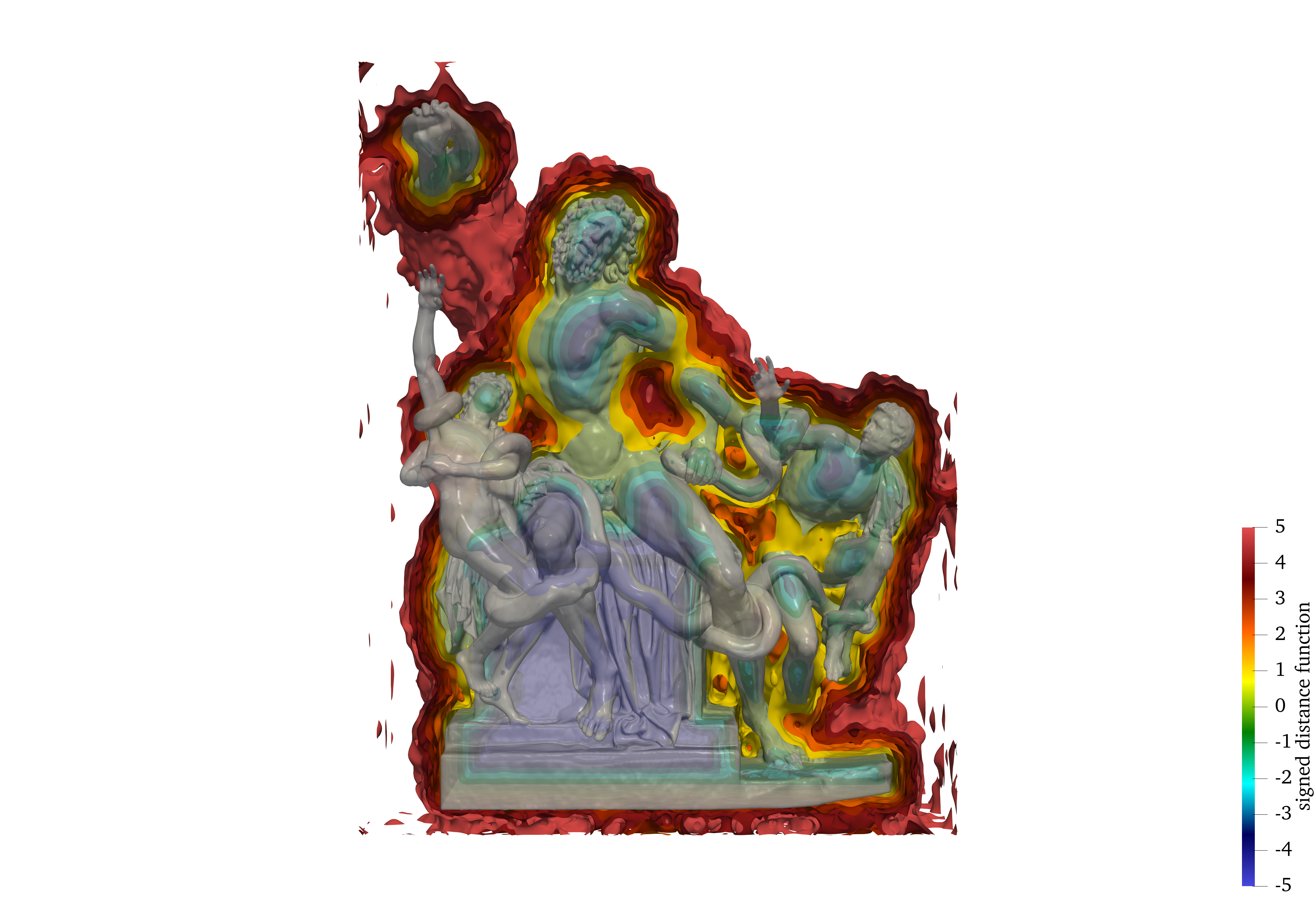}};
\draw(3.2,0)node{\includegraphics[scale=0.042,clip,trim= 1000 100 1000 100]{
  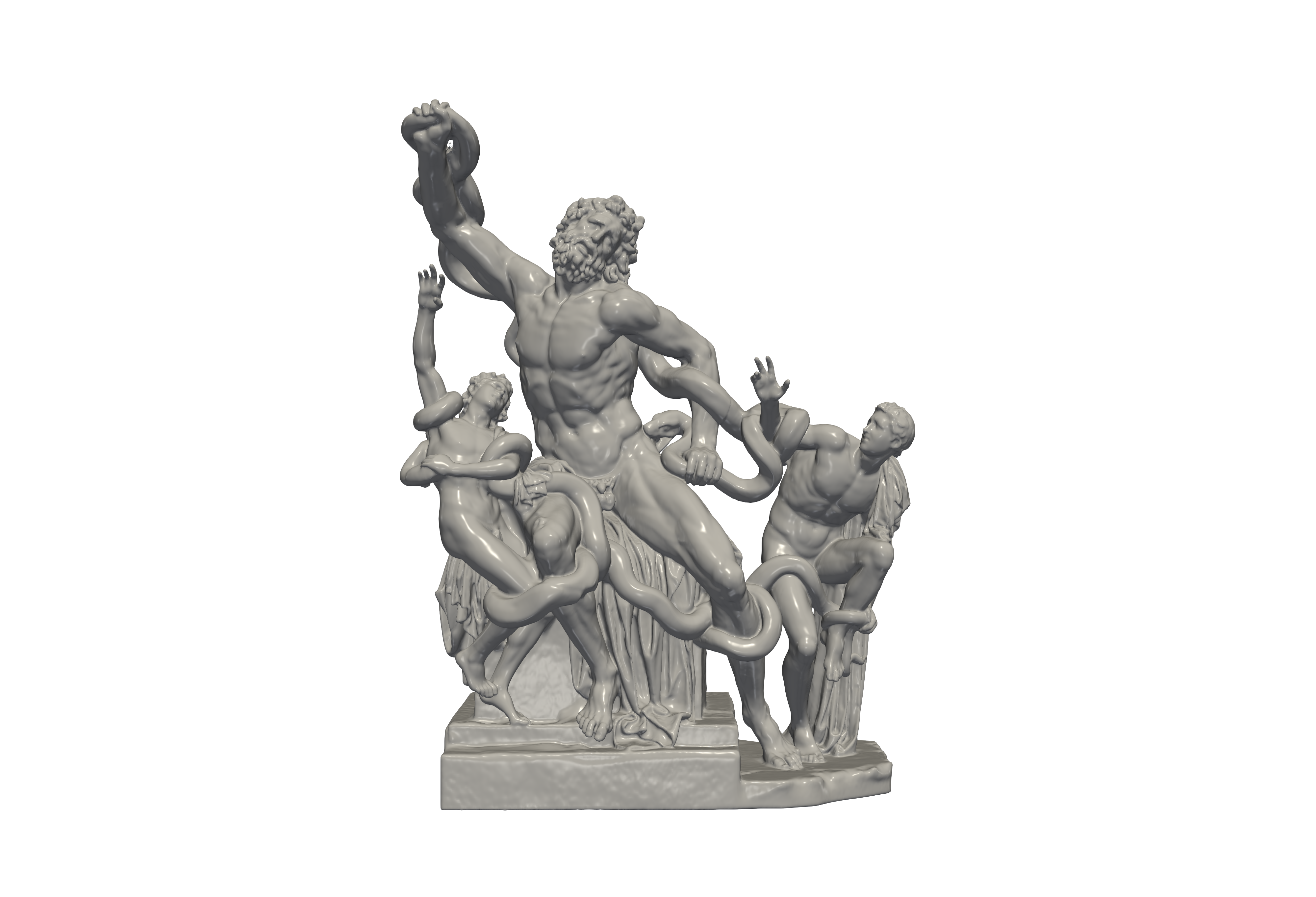}};
\draw(5,1.1)node{\includegraphics[scale=0.02,clip,trim= 1400 850 1200 400,
    frame]{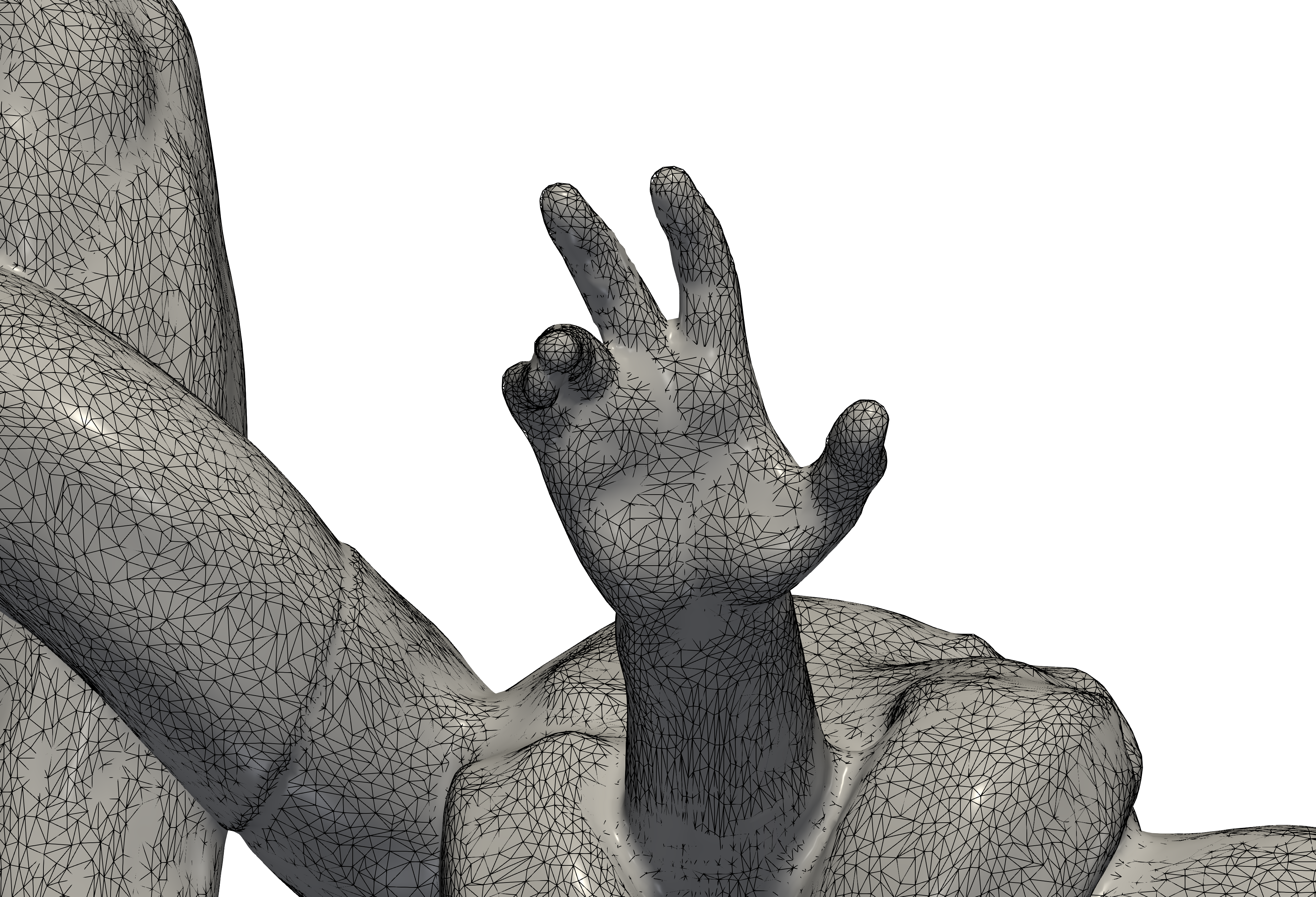}};
\draw[line width=0.5pt](3.58,0.2)rectangle(3.82,0.45);
\draw[line width=0.5pt](3.82,0.2)--(4.5,0.57);
\draw[line width=0.5pt](3.82,0.45)--(4.5,1.635);
\end{tikzpicture}
\caption{\label{HM_fig:SDFInterpolation}Signed distance function interpolation
from a surface mesh of the Laokoon group. The left image shows a subsample of
the signed distance function used for approximation, the image in the middle
shows the levelsets of the signed distance function for the values
\(\{-5,-4,\ldots,5\}\). The image on the right is the zero levelset with a zoom 
of the original surface mesh and the reconstruction of the right son's left
hand.}
\end{center}
\end{figure}

To demonstrate the capabilities of the samplet matrix compression, we consider a
surface reconstruction problem, similarly to \cite{HM_CarrEtal01}. Given a planar 
triangulation resulting from a 3D scan of the Laokoon group (the scan is 
provided by the Statens Museum for Kunst), we generate uniform samples of the
signed distance by using the about 500\,000 vertices the surface mesh and
another 250\,000 random points within the bounding box of the Laokoon group.
This results in \(N=750\,000\) data sites in total. The left image in
Figure~\ref{HM_fig:SDFInterpolation} shows a uniform subsample of size 100\,000
of the data points. For the interpolation, we consider the exponential kernel 
\(k_{1/2}\) from \eqref{HM_eg:1/2} with lengthscale parameter \(\ell=0.01\),
where the data sites are rescaled to the hypercube $[0,1]^3$. The kernel matrix
is compressed by using samplets with \(q+1=4\) vanishing moments and the linear 
system \eqref{HM_eg:LSEreg} of equations is solved with the regularization
parameter \(\mu=10^{-8}\) by the conjugate gradient method. 

The interpolated signed distance function is then evaluated at a uniform grid of
125\,000\,000 points, that is 500 points per axis direction. The evaluation is
performed by the fast multipole method developed in \cite{HM_HMQ24} using total degree
polynomials of degree 6 for the interpolation of the kernel function in the
farfield. The image in the middle of Figure~\ref{HM_fig:SDFInterpolation} shows
the levelsets of the signed distance function for the values 
\(\{-5,-4,\ldots,5\}\). The right image shows the zero levelset, which
corresponds to the desired surface. From the zoom of the  right son's hand, it
can be seen how the surface is smoothened in comparison to the original surface
mesh.

%===============================================================================
\subsection{Multiresolution kernel matrix algebra}
%===============================================================================
Let $\Omega\subset\mathbb{R}^d$ be a smooth domain and assume that the set of
sites $X\subset\Omega$ is \emph{asymptotically distributed uniformly modulo
one}, see \cite{HM_LP10}. This means, we have for every Riemann integrable function
$f:\Omega\to\mathbb{R}$ that
\begin{equation*}%\label{HM_eg:modulo_one}
%===============================================================================
  \bigg|\int_\Omega f({\bs x})\d{\bs x}
  -\frac{|\Omega|}{N}\sum_{i=1}^N f({\bs x}_i)\bigg|
  \to 0\quad\text{as $N\to\infty$}.
\end{equation*}
Under this condition, the kernel matrix ${\bs K}$ corresponds to a Nystr\"om
discretization of an associated integral operator $K$. Namely, the reproducing
kernel $\kernel(\cdot,\cdot)$ gives rise to the compact integral operator
\begin{equation}\label{HM_eg:pseudo}
%===============================================================================
  K:L^2(\Omega)\to L^2(\Omega), 
  \quad u \mapsto \int_\Omega \kernel(\cdot,{\bs y}) u({\bs y})\d{\bs y}.
\end{equation}

For many practical applications, \eqref{HM_eg:pseudo} constitutes a classical 
pseudo-differential operator of negative order $s<0$, especially in case of
Mat\'ern kernels. We refer to, for example, \cite{HM_HorIII,HM_Rodino,HM_Seeley,HM_Taylor81}
for the details of this theory, including the subsequent developments. We are
interested here in pseudo-differential operators $K$ which belong to the
subclass $OPS_{cl,1}^s$ of \emph{analytic pseudo-differential operators}, see 
\cite{HM_Rodino}. Their kernels are known to be asymptotically smooth, satisfying
for all $({\bs x},{\bs y})\in(\Omega\times\Omega)\setminus\Delta$ the refined 
decay property
\begin{equation}\label{HM_eg:kerneldecay}
%===============================================================================
\bigg|\frac{\partial^{|\bs\alpha|+|\bs\beta|}}
  	{\partial{\bs x}^{\bs\alpha}
  	\partial{\bs y}^{\bs\beta}} \kernel({\bs x},{\bs y})\bigg|
\lesssim \frac{(|\bs\alpha|+|\bs\beta|)!}
{\rho^{|\bs\alpha|+|\bs\beta|}
\|{\bs x}-{\bs y}\|_2^{s+d+|\bs\alpha|+|\bs\beta|}}
\end{equation}
uniformly in $\bs\alpha,\bs\beta\in\mathbb{N}^d$. Note that, in case of 
reproducing kernels, we always have $s+d<0$ due to the continuity of the kernel
such that \eqref{HM_eg:kerneldecay} implies also \eqref{HM_eg:kernel_estimate},
required for the samplet matrix compression in 
Section~\ref{HMsec:matrix_compression}.

The addition ${\bf K}+{\bf K}'$ of two compressed kernel matrices in compressed
format is obvious. For the computation of the matrix product 
${\bf K}\cdot{\bf K}'$, we make use of the fact that it corresponds to the 
concatenation $K\circ K'$ of the underlying pseudo-differential operators since
\begin{equation}\label{HM_eg:multiplication}
%===============================================================================
  \int_\Omega \kernel({\bs x}_i,{\bs z})\kernel'({\bs z},{\bs x}_j)\d{\bs z}
  \approx \frac{|\Omega|}{N}\sum_{k=1}^N \kernel({\bs x}_i,{\bs x}_k)
  \kernel'({\bs x}_k,{\bs x}_j).
\end{equation}
If $K\in OPS_{cl,1}^s$ and $K'\in OPS_{cl,1}^{s'}$, then there holds 
$K\circ K'\in OPS_{cl,1}^{s+s'}$ and, thus, $K\circ K'$ is compressible. In 
accordance with \cite{HM_HMSS}, we can therefore compute the matrix product 
${\bf K}'\cdot{\bf K}$ in loglinear complexity on the given matrix pattern.

Consider next a symmetric and positive definite kernel function 
$\kernel(\cdot,\cdot)$ such that the associated pseudo-differential operator 
satisfies $K\in OPS_{cl,1}^s$ with $s+d<0$. As shown in \cite{HM_HMSS}, the inverse
of $K+\mu\operatorname{Id}$ is of the form $\mu^{-1}\operatorname{Id}-K'$ with 
$K'$ being likewise a pseudo-differential operator of class $OPS_{cl,1}^s$. 
Especially, the kernel function $\kernel'$ which underlies the operator $K'$ by
the Schwartz kernel theorem, see \cite{HM_HorI} for instance, is symmetric,
positive definite, and likewise asymptotically smooth. Thus, the inverse kernel
matrix $({\bs K}+\mu {\bs I})^{-1}$ is also compressible and can, in view of
\eqref{HM_eg:multiplication}, be efficiently approximated by selected inversion,
see \cite{HM_selinv}, of the associated pattern, see \cite{HM_HMSS} for the details.

More complicated matrix functions like powers of the kernel matrix or the matrix
exponential become accessible, too, by using contour integrals, see 
\cite{HM_Hale2008}. Indeed, the kernel matrix algebra proposed in \cite{HM_HMSS} has
the property that the arithmetics is exact on the prescribed (fixed) pattern.
We refer the reader to \cite{HM_HMSS} for specific examples and results.

%===============================================================================
\subsection{Samplet basis pursuit}\label{HMsec:sparsity}
%===============================================================================
Sparsity constraints are widely used in computational learning, statistics, and
signal processing such as deblurring, feature selection and compressive sensing,
see \cite{HM_Candes,HM_BasisPursuit,HM_Donoho,HM_learning,HM_Tao} for example. Sparsity
constraints are imposed by adding an $\ell^1$-penalty term to the functional to
be minimized. However, such sparsity constraints make only sense if a basis is
used for the discretization where the data become sparse. In the past, mostly
wavelets bases, Fourier bases, or frames like curvelets, contourlets, and
shearlets have been used as they are known to give raise to sparse
representations, see \cite{HM_curvelets,HM_ISTA,HM_contourlets,HM_rauhut,HM_shearlets} for
example. However, such discretization concepts are based on regular grids and it
is not obvious how to extend them to scattered data as they appear often in
machine leaning.

In this section, we shall therefore discuss \(\ell^1\)-regularization for
scattered data approximation with respect to samplet coordinates. To this end,
we consider the functional
\begin{equation}\label{HM_eg:MSLASSO}
%===============================================================================
\min_{{\bs\alpha}\in\Rbb^N}\frac 1 2\|{\bs h}-{\bs K}{\bs\alpha}\|_2^2
+\sum_{i=1}^N w_i|\beta_i|,\ \text{where}\ 
{\bs\beta}={\bs T}{\bs\alpha}.
\end{equation}
The weight vector ${\bs w}=[w_i]_i\in\Rbb^N$ plays the role of the
regularization parameter, where each coefficient is regularized individually. We
refer to, e.g., \cite{HM_ISTA,HM_Lorenz,HM_RT} for the analysis of the regularizing
properties and for appropriate parameter choice rules. 

Numerical algorithms to solve the optimization problem \eqref{HM_eg:MSLASSO} are
based on soft-thresholding. Indeed, the sparsity constrained minimization
problem \eqref{HM_eg:MSLASSO} can be recast into the root finding problem
\begin{equation}\label{HM_eg:SSN}
%================================================
{\bs 0}=
{\bs\beta}^\star-\operatorname{SS}_{\gamma{\bs w}}\big(
{\bs\beta}^\star+\gamma({\bs K}^{\Sigma})^\intercal({\bs h}^{\Sigma}
-{\bs K}^{\Sigma}{\bs\beta}^\star)\big),
\end{equation}
where $\gamma>0$ and 
\[
\operatorname{SS}_{{\bs w}}({\bs v})
\isdef\operatorname{sign}({\bs v})\max\{{\bs 0},|\bs v|-{\bs w}\}
\]
is the soft-shrinkage operator. Problem \eqref{HM_eg:SSN} can efficiently be
solved by the semi-smooth Newton method, 
see \cite{HM_BHM,HM_LorenzGriesse}.

In the spirit of \cite{HM_CDS98,HM_GCG05}, also a dictionary of multiple kernels can
be employed in \eqref{HM_eg:MSLASSO}. Given the kernels 
$\kernel_1,\ldots,\kernel_L$, we are then looking for a sparse representation of
the form
\[
s_h=\sum_{\ell=1}^L\sum_{i=1}^N\alpha_{i}^{({\ell})}
\kernel_\ell({\bs x}_i,\cdot)
=\sum_{\ell=1}^L\sum_{i=1}^N\beta_{i}^{({\ell})}\psi_{i}^{(\ell)}.
\]
Setting 
\[
\bs K\isdef[{\bs K}_1,\ldots,{\bs K}_L],\quad 
{\bs K}_\ell\isdef[\kernel_\ell({\bs x}_i,{\bs x}_j)]_{i,j=1}^N,
\] 
and 
\[
\bs\alpha^\intercal\isdef
[{\bs\alpha}_{1}^\intercal,\ldots,{\bs\alpha}_{L}^\intercal],
\quad {\bs\alpha}_j\in\Rbb^N,
\]
this approach also amounts to \eqref{HM_eg:MSLASSO} with obvious modifications.
The most important difference to the original problem is, of course, that the
matrix $\bs K$ is not quadratic any more.This means that the underlying linear
system $\bs K\bs\alpha = \bs h$ of equations is underdetermined.

\begin{figure}[htb]
\begin{center}
\begin{tikzpicture}
\draw(0,1.1)node{\tiny February};
\draw(3.4,1.1)node{\tiny June};
\draw(6.8,1.1)node{\tiny October};
\draw(-1.9,0)node[rotate=90]{\tiny reconstruction};
\draw(-1.9,-1.9)node[rotate=90]{\tiny contribution \(k_1\)};
\draw(-1.9,-3.8)node[rotate=90]{\tiny contribution \(k_2\)};
\draw(-1.9,-5.7)node[rotate=90]{\tiny pointwise relative error};
\draw(0,0)node{\includegraphics[scale=0.032,clip,trim=500 350 500 350]{
  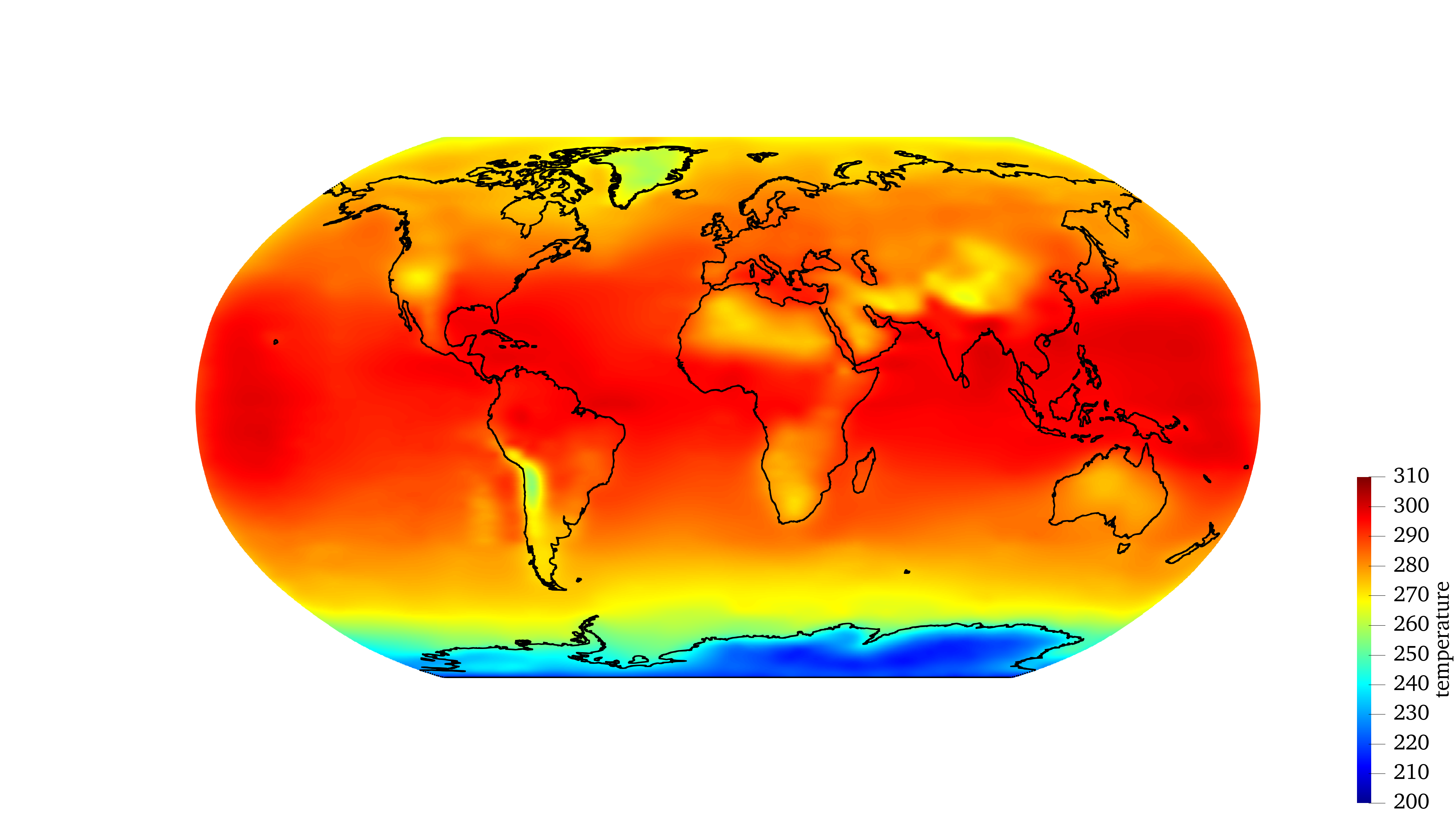}};
\draw(3.4,0)node{\includegraphics[scale=0.032,clip,trim=500 350 500 350]{
  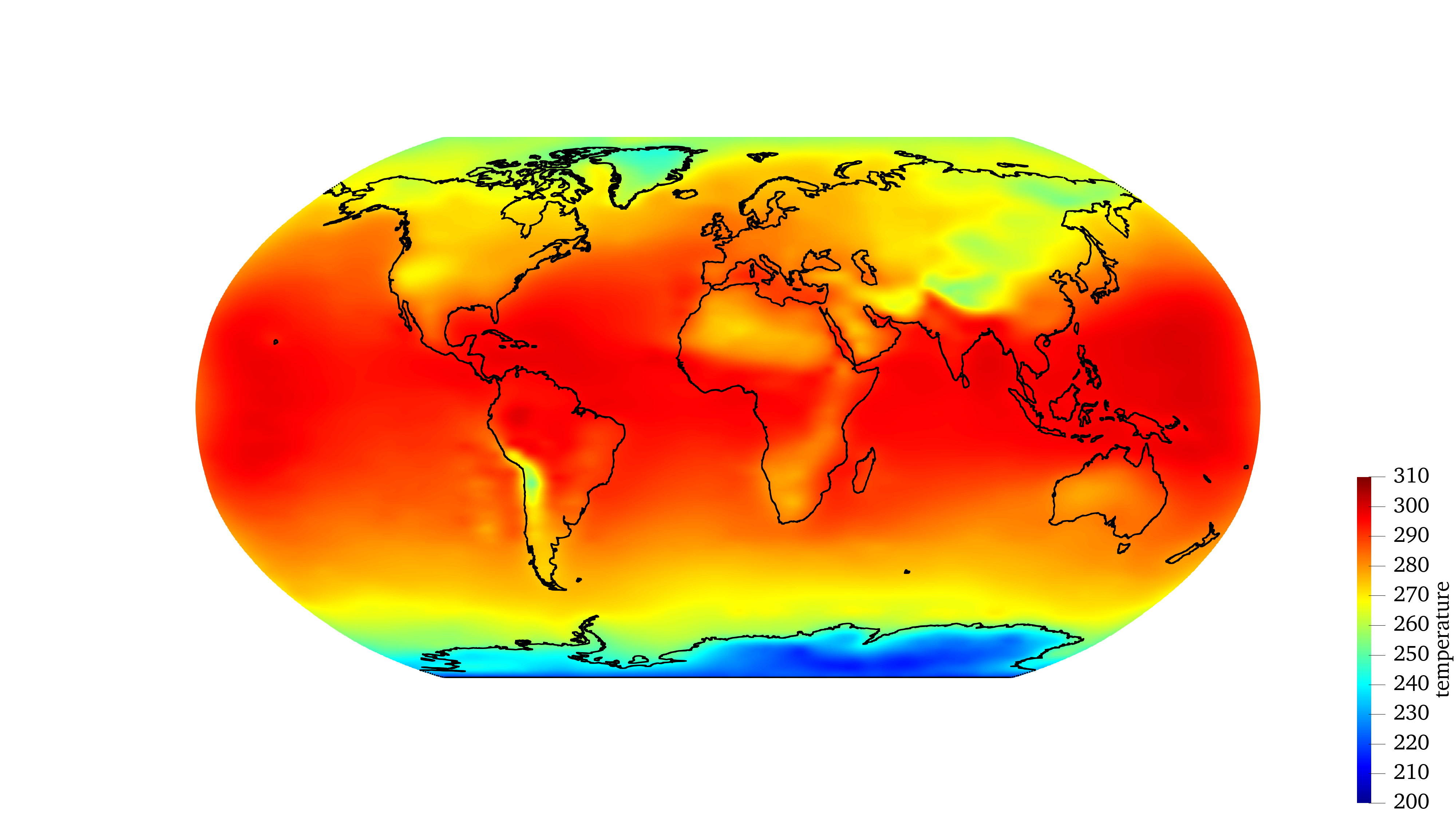}};
\draw(6.8,0)node{\includegraphics[scale=0.032,clip,trim=500 350 500 350]{
  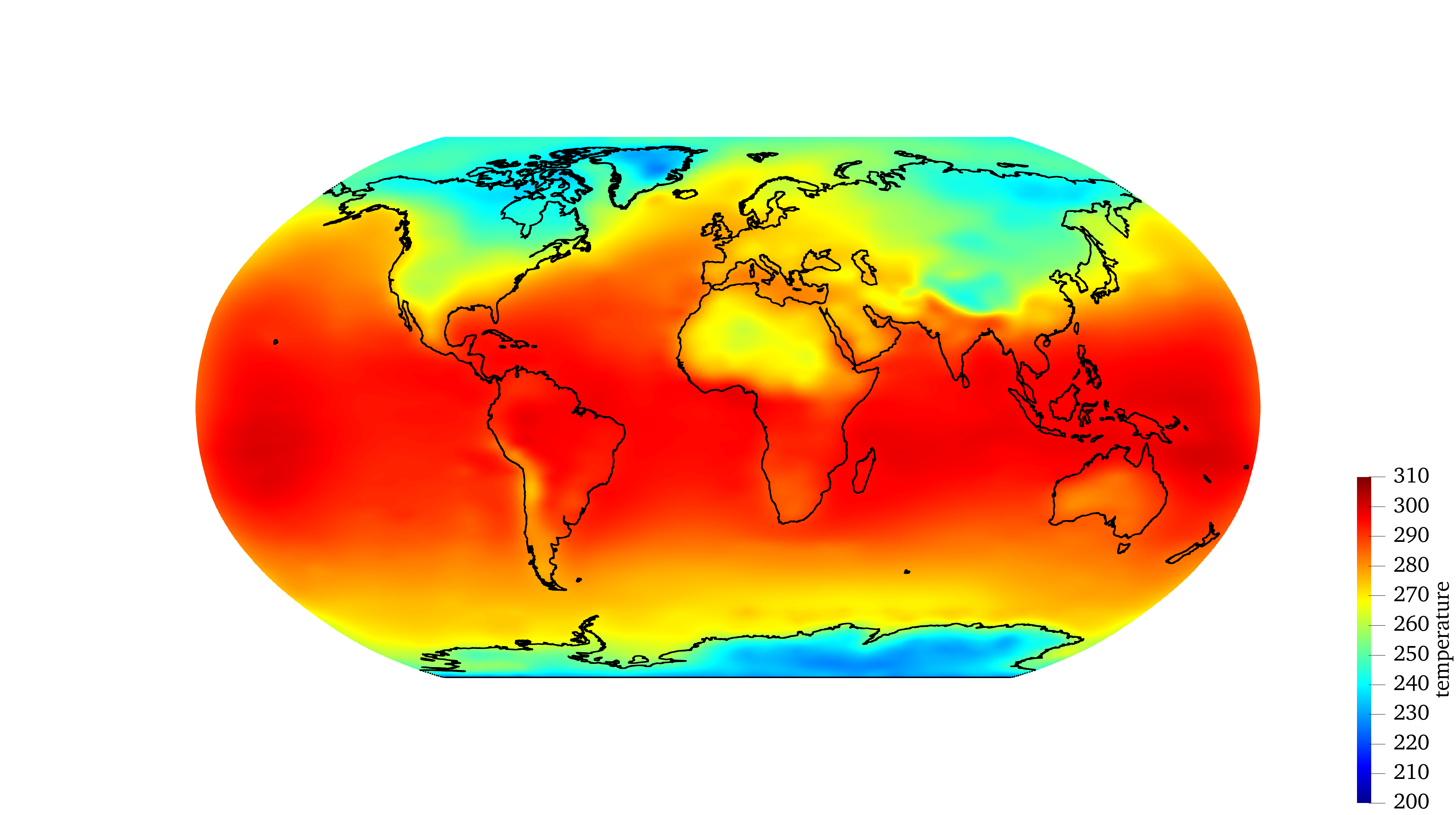}};
\draw(0,-1.9)node{\includegraphics[scale=0.032,clip,trim=500 350 500 350]{
  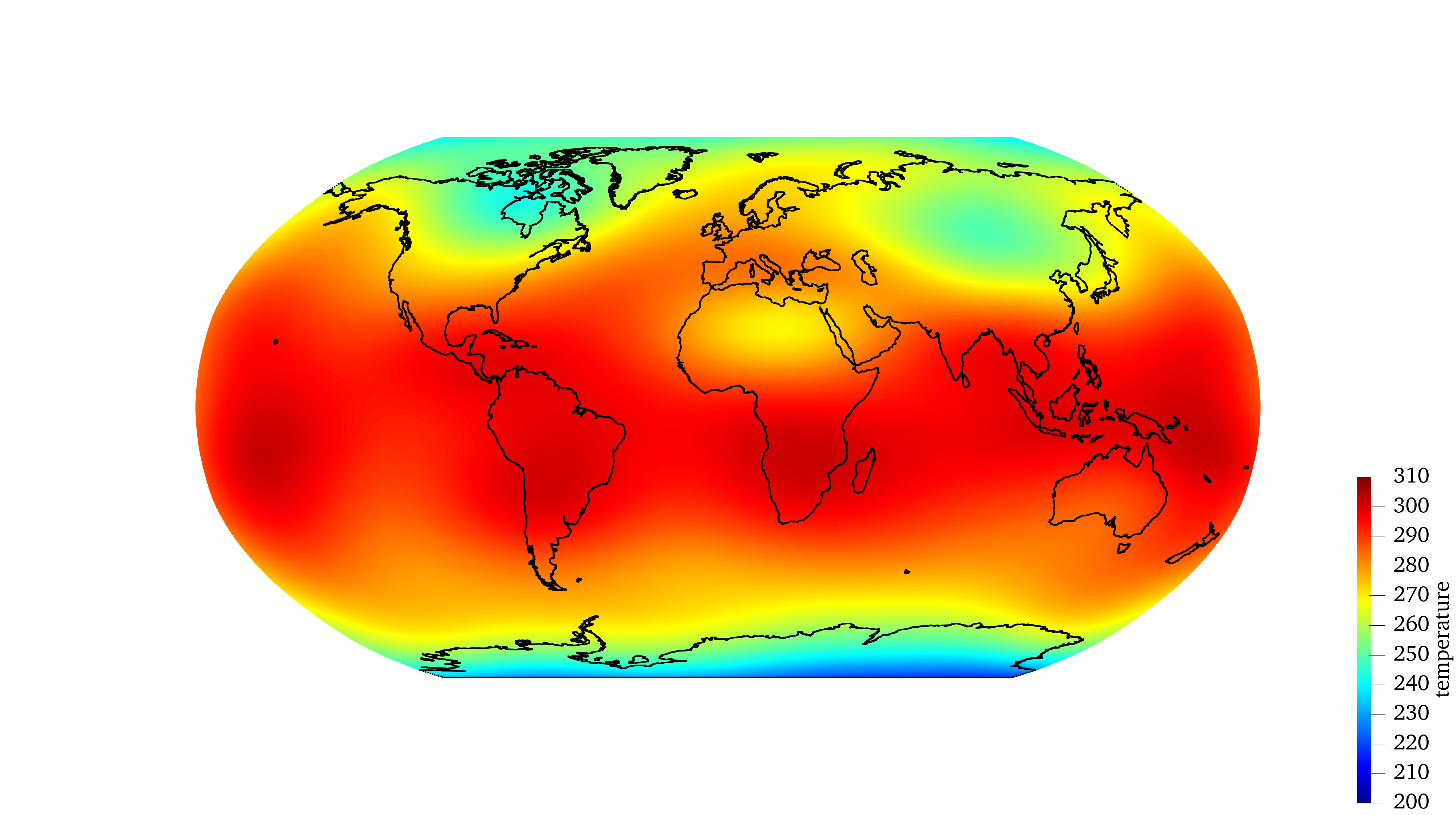}};
\draw(3.4,-1.9)node{\includegraphics[scale=0.032,clip,trim=500 350 500 350]{
  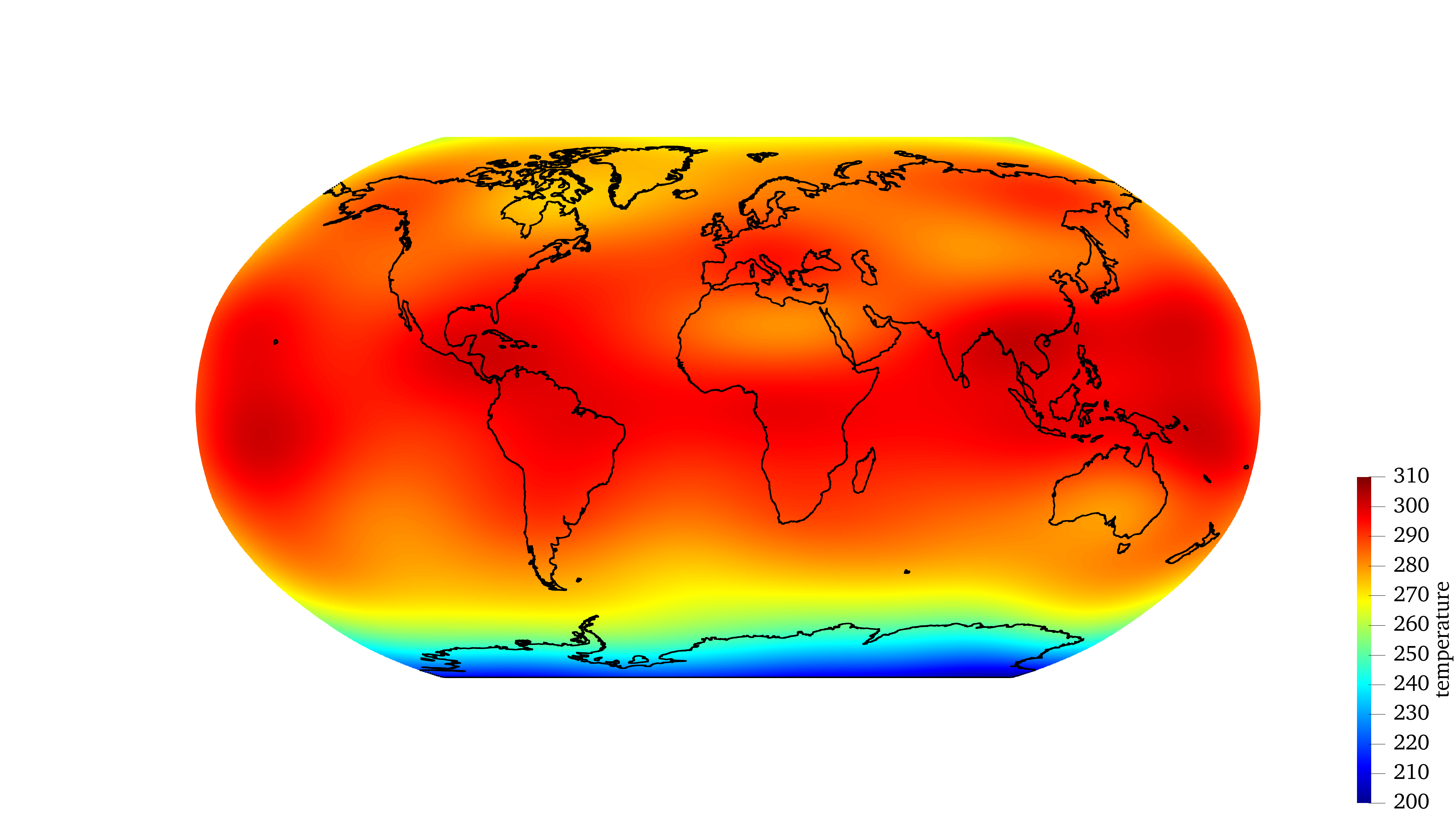}};
\draw(6.8,-1.9)node{\includegraphics[scale=0.032,clip,trim=500 350 500 350]{
  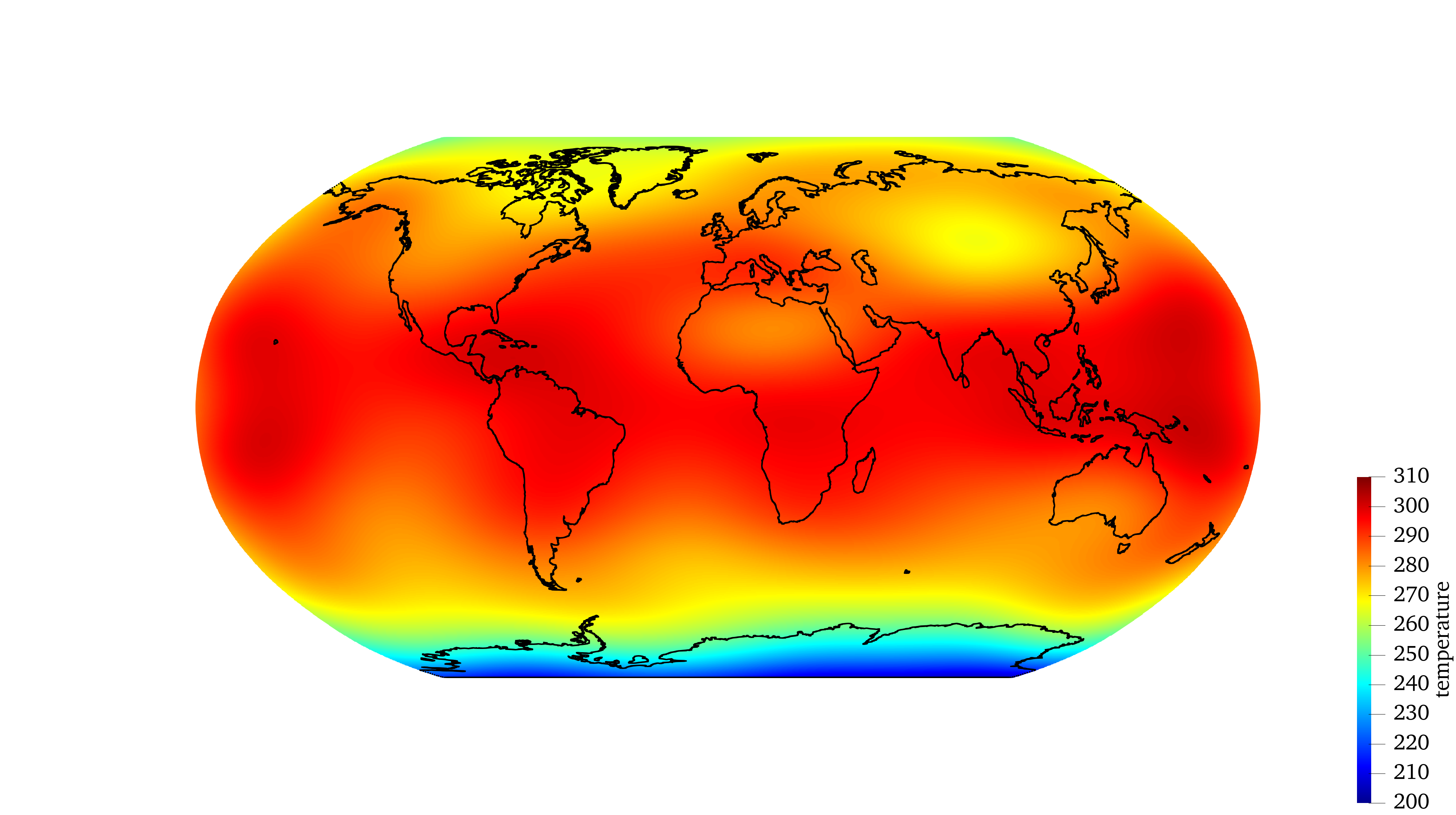}};
\draw(0,-3.8)node{\includegraphics[scale=0.032,clip,trim=500 350 500 350]{
  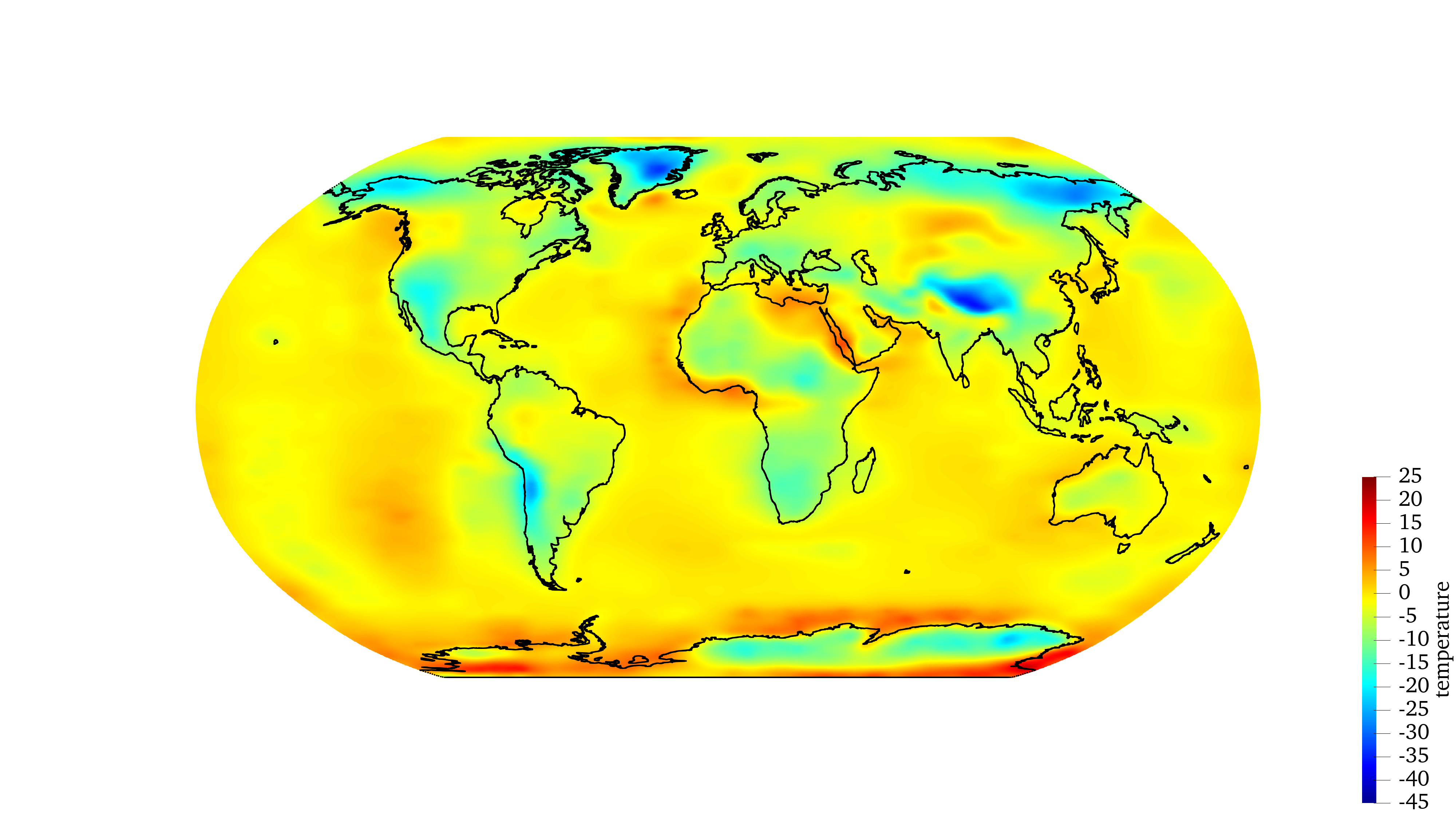}};
\draw(3.4,-3.8)node{\includegraphics[scale=0.032,clip,trim=500 350 500 350]{
  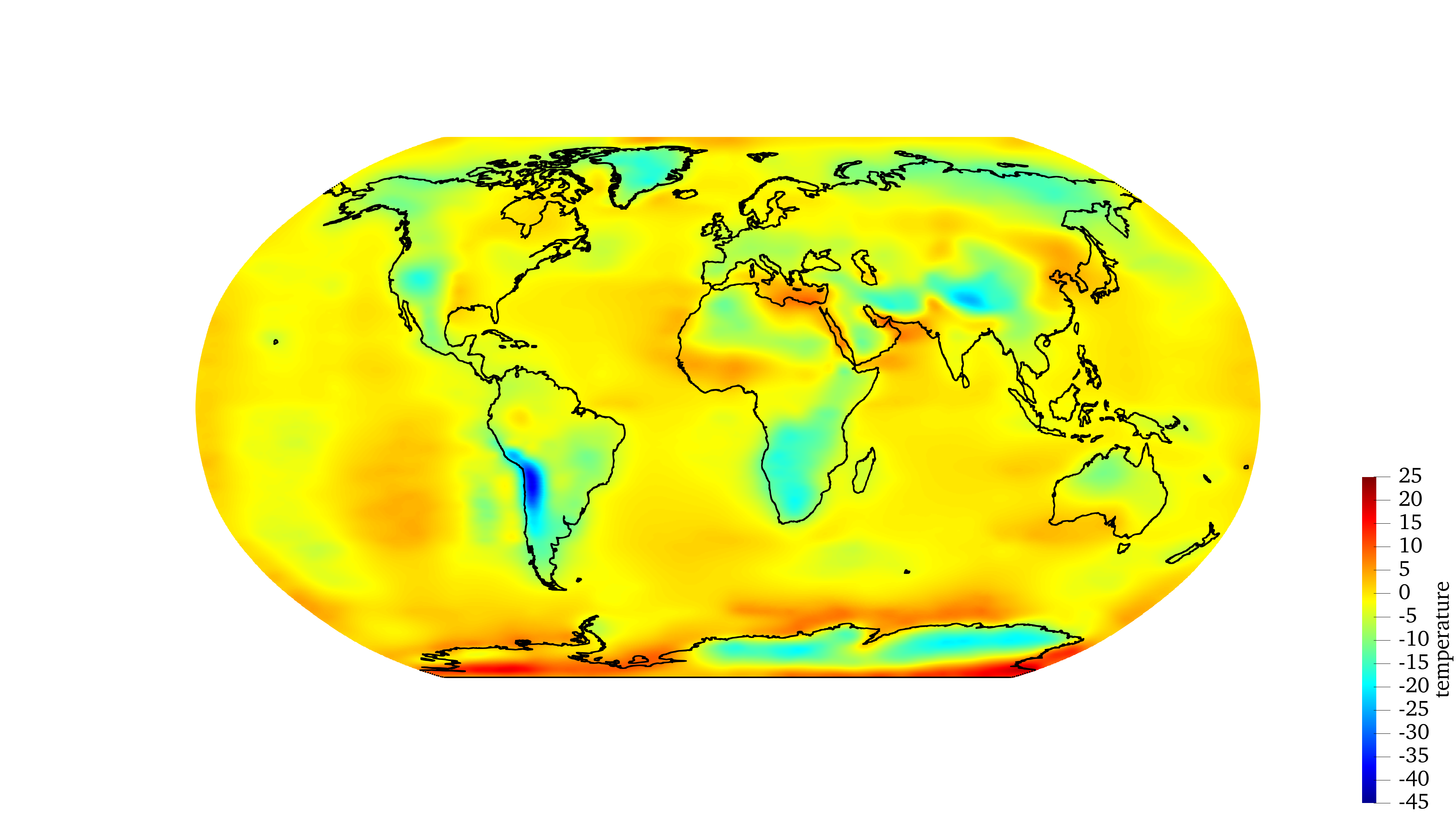}};
\draw(6.8,-3.8)node{\includegraphics[scale=0.032,clip,trim=500 350 500 350]{
  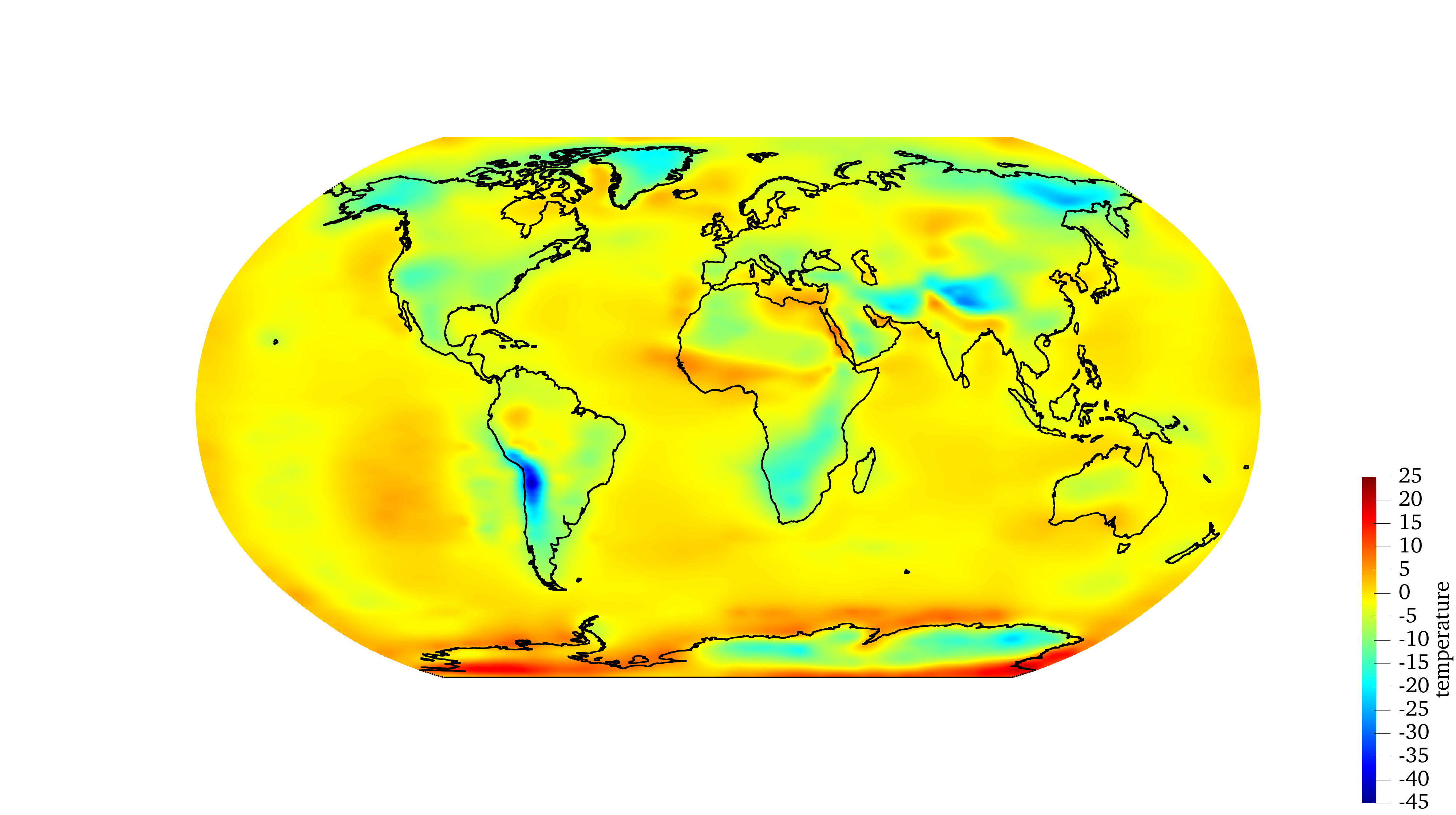}};
\draw(0,-5.7)node{\includegraphics[scale=0.032,clip,trim=500 350 500 350]{
  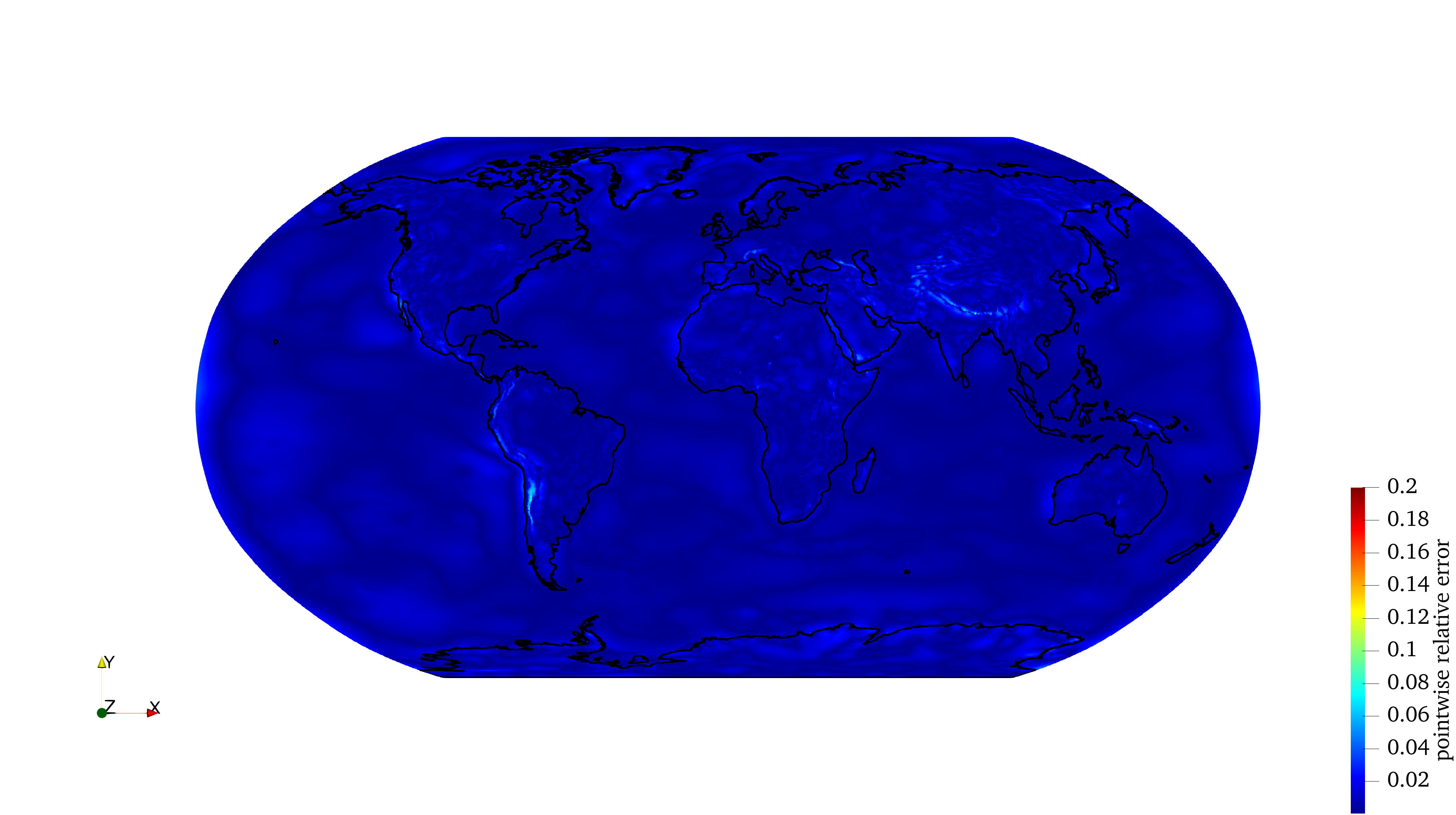}};
\draw(3.4,-5.7)node{\includegraphics[scale=0.032,clip,trim=500 350 500 350]{
  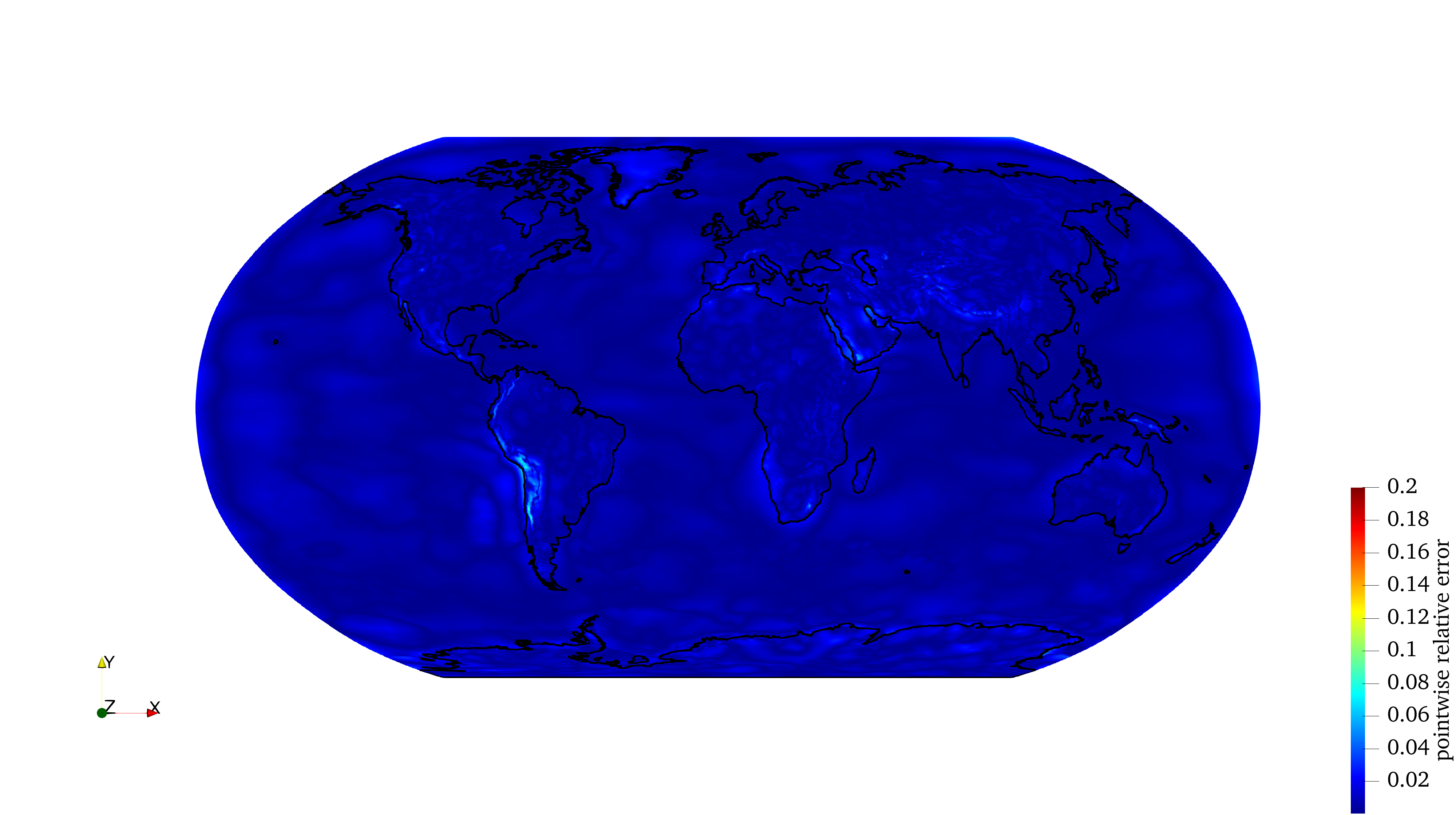}};
\draw(6.8,-5.7)node{\includegraphics[scale=0.032,clip,trim=500 350 500 350]{
  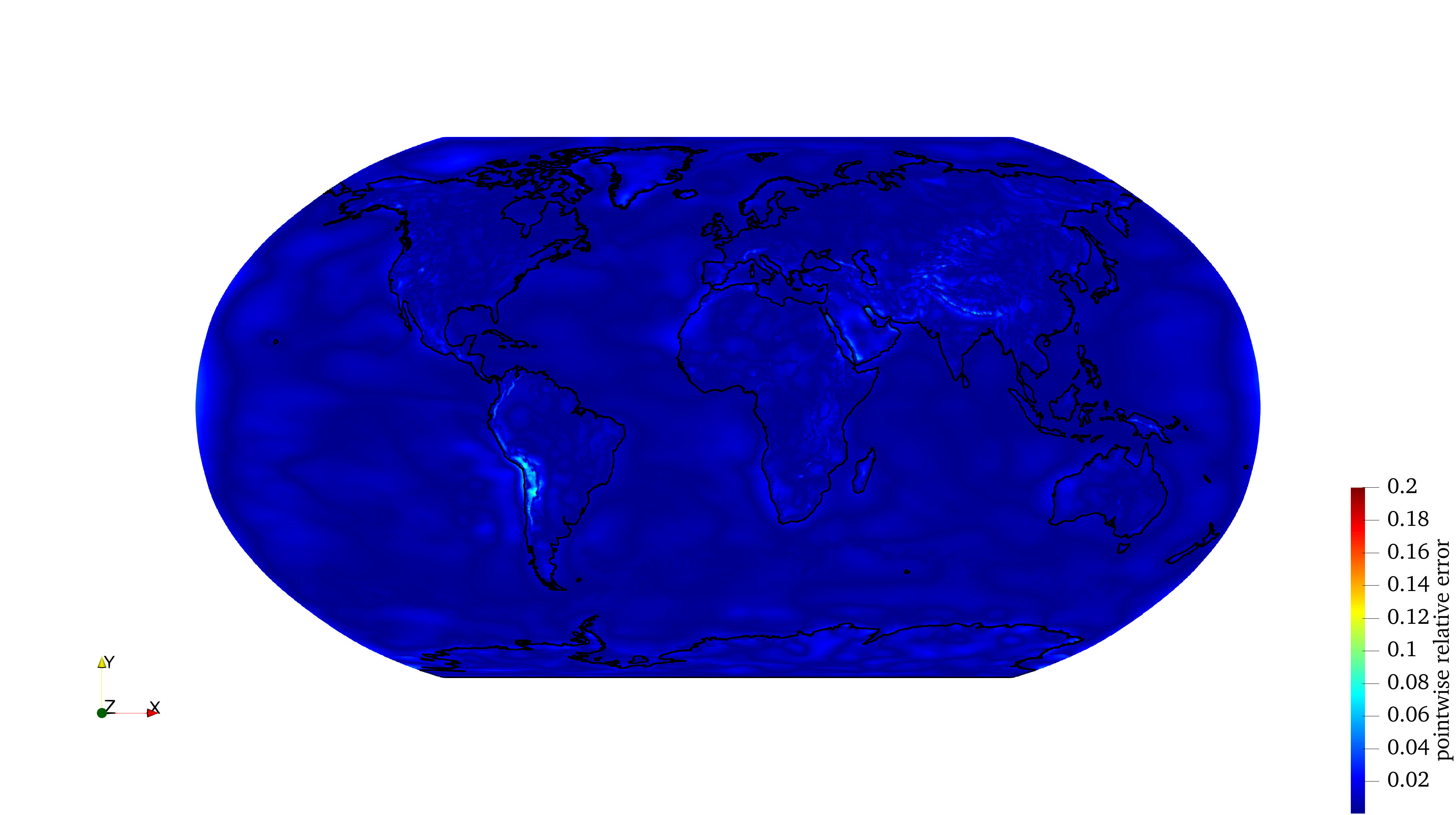}};
\draw(8.9,-5.7)node{\includegraphics[scale=0.052,clip,trim=3620 0 0 1250]{
  ./Images/SSN_err_9_lins.png}};
\draw(8.9,-3.8)node{\includegraphics[scale=0.052,clip,trim=3620 0 0 1250]{
  ./Images/SSN_k2_9.png}};
\draw(8.9,-1.9)node{\includegraphics[scale=0.052,clip,trim=3620 0 0 1250]{
  ./Images/SSN_k1_9.png}};
\draw(8.9,0)node{\includegraphics[scale=0.052,clip,trim=3620 0 0 1250]{
  ./Images/SSN_rec_9.png}};
\end{tikzpicture}
\caption{\label{HM_fig:SSNrec} Reconstruction of the average temperature for
February, June, October 2022 (top row), contributions of \(k_1,k_2\) (second
and third row) and pointwise relative error (bottom row).}
\end{center}
\end{figure}

To illustrate the approach, we consider the sparse reconstruction of the
temperature data from the ERA5 data set in space-time. We employ two space-time
kernels, namely
\[
\kernel_1({\bs z},{\bs z}')\isdef
k_{3/2}(\|{\bs x}-{\bs y}\|_2)k_{\text{per}}(|t-t'|)
\]
and
\[
\kernel_2({\bs z},{\bs z}')\isdef k_{1/2}(\|{\bs z}-{\bs z}'\|_2),
\]
where ${\bs z}\isdef({\bs x},t)$ and ${\bs z}'\isdef({\bs y},t')$.
Herein, the lengthscale parameter is set to \(\ell=0.2\) for \(k_{3/2}\) and to
\(\ell=0.01\) for \(k_{1/2}\), compare \eqref{HM_eg:MaternKernels}. Moreover, we
define the periodic Gaussian kernel 
\(k_{\text{per}}(r)\isdef e^{-50\sin^2(\pi r)}\). The kernel \(\kernel_1\) is a
tensor product kernel with relatively large lengthscale parameter. This kernel
is intended to capture the smooth parts of the temperature distribution over
time. The second kernel \(\kernel_2\) is a quite rough exponential kernel in
space time and chosen to capture sharp features. 

For the numerical computation, the set of data sites is obtained by adaptively
sampling \(100\,000\) data sites per time step, compare 
Section~\ref{HMsubsec:AdaptiveSubs}. This results in \(N=1\,200\,000\) data
sites in total. These data sites are rescaled to the unit hyper-cube 
\([0,1]^3\). For the compression of the respective kernel matrices, we apply
samplets with \(q+1=4\) vanishing moments. We compute the coefficients in 
\eqref{HM_eg:MSLASSO} by using an iteratively regularized version of the
semi-smooth Newton method, see \cite{HM_BHM} for details. The weight vector
${\bs w}$ is set to \(w_i=10^{-6}\) for all $i=1,\ldots,N$. We obtain 
\(\|{\bs\beta}_1\|_0=671\) and \(\|{\bs\beta}_2\|_0=5\,883\) non-vanishing
coefficients, resulting in a relative in-sample error of \(8.18\cdot 
10^{-3}\) in the Euclidean norm, while the value of the functional is 
\(4.76\cdot 10^{-5}\). The pointwise relative error which is obtained on the
full data set  is smaller than \(9.73\%\) for all time steps.

The top row of Figure~\ref{HM_fig:SSNrec} shows the reconstruction evaluated for
the full data set for the months February, June, and October. The second row 
shows the contribution of \(\kernel_1\). As can be seen, \(\kernel_1\) captures
the coarse-scale structure of the temperature distribution with relatively few
coefficients. The third row shows the contribution of \(\kernel_2\). The kernel
localizes at the coast lines and mountains and represents sharp features, which
requires the major part of the non-zero coefficients. The pointwise relative
approximation error is finally found in the last row.
%===============================================================================
\bibliographystyle{plain}
\bibliography{literature}

\begin{thebibliography}{10}

\bibitem{HM_AHK14}
D.~Alm, H.~Harbrecht, and U.~Kr{\"a}mer.
\newblock The \(\mathcal{H}^2\)-wavelet method.
\newblock {\em J. Comput. Appl. Math.}, 267:131--159, 2014.

\bibitem{HM_Alp93}
B.K. Alpert.
\newblock A class of bases in \({L}^2\) for the sparse representation of
  integral operators.
\newblock {\em SIAM J. Math. Anal.}, 24(1):246--262, 1993.

\bibitem{HM_Aronszajn50}
N.~Aronszajn.
\newblock Theory of reproducing kernels.
\newblock {\em Trans. Amer. Math. Soc.}, 68(3):337--404, 1950.

\bibitem{HM_BHM}
D.~Baroli, H.~Harbrecht, and Multerer.
\newblock Samplet basis pursuit: Multiresolution scattered data approximation
  with sparsity constraints.
\newblock {\em IEEE Trans. Sign. Proc.}, 72:1813--1823, 2024.

\bibitem{HM_BA04}
A.~Berlinet and C.~Thomas-Agnan.
\newblock {\em Reproducing Kernel {H}ilbert Spaces in Probability and
  Statistics}.
\newblock Springer Science \& Business Media, New York, 2004.

\bibitem{HM_BCR}
G.~Beylkin, R.~Coifman, and V.~Rokhlin.
\newblock {T}he fast wavelet transform and numerical algorithm.
\newblock {\em Comm. Pure Appl. Math}, 44:141--183, 1991.

\bibitem{HM_BD97}
P.~Binev and R.A. DeVore.
\newblock Fast computation in adaptive tree approximation.
\newblock {\em Numer. Math.}, 97:193--217, 2004.

\bibitem{HM_BFPR+73}
M.~Blum, R.W. Floyd, V.R. Pratt, R.L. Rivest, and R.E. Tarjan.
\newblock Time bounds for selection.
\newblock {\em J. Comput. Syst. Sci.}, 7(4):448--461, 1973.

\bibitem{HM_Boe10}
S.~B{\"o}rm.
\newblock {\em Efficient Numerical Methods for Mon-Local Operators:
  \(\mathcal{H}^2\)-Matrix Compression, algorithms and analysis}.
\newblock European Mathematical Society, Z{\"u}rich, 2010.

\bibitem{HM_Bra13}
D.~Braess.
\newblock {\em Finite Elemente. Theorie, schnelle L{\"o}ser und Anwendungen in
  der Elastizit{\"a}ts\-theorie}.
\newblock Springer, Berlin-Heidelberg, 2013.

\bibitem{HM_Candes}
E.~Cand{\`e}s, J.~Romberg, and T.~Tao.
\newblock Stable signal recovery from incomplete and inaccurate measurements.
\newblock {\em Comm. Pure Appl. Math.}, 59(8):1207--1223, 2006.

\bibitem{HM_curvelets}
E.J. Candes and D.L. Donoho.
\newblock Curvelets: {A} surprisingly effective nonadaptive representation for
  objects with edges.
\newblock In A.~Cohen, C.~Rabut, and L.~Schumaker, editors, {\em Curves and
  Surface Fitting: Saint-Malo 1999}, page 105–120, Nashville, 2000.
  Vanderbilt University Press.

\bibitem{HM_CarrEtal01}
J.C. Carr, R.K. Beatson, J.B. Cherrie, T.J. Mitchell, W.R. Fright, B.C.
  McCallum, and T.R. Evans.
\newblock Reconstruction and representation of {3D} objects with radial basis
  functions.
\newblock In {\em Proceedings of the 28th annual conference on Computer
  graphics and interactive techniques}, SIGGRAPH '01, pages 67--76, New York,
  2001. Association for Computing Machinery.

\bibitem{HM_BasisPursuit}
S.~Chen and D.L. Donoho.
\newblock Basis pursuit.
\newblock In {\em Proceedings of 1994 28th Asilomar Conference on Signals,
  Systems and Computers}, volume~1, pages 41--44, Pacific Grove, CA, USA, 1994.
  IEEE.

\bibitem{HM_CDS98}
S.S. Chen, D.L. Donoho, and M.A. Saunders.
\newblock Atomic decomposition by basis pursuit.
\newblock {\em SIAM J. Sci. Comput.}, 20(1):33--61, 1998.

\bibitem{HM_Chui}
C.K. Chui.
\newblock {\em An Introduction to Wavelets}.
\newblock Academic Press, San Diego, 1992.

\bibitem{HM_Quak}
C.K. Chui and E.~Quak.
\newblock Wavelets on a bounded interval.
\newblock {\em Numer. Meth. Approx. Theory}, 9:53--75, 1992.

\bibitem{HM_Coh03}
A.~Cohen.
\newblock {\em Numerical Analysis of Wavelet Methods}.
\newblock Elsevier, Amsterdam, 2003.

\bibitem{HM_CM06}
R.R. Coifman and M.~Maggioni.
\newblock Diffusion wavelets.
\newblock {\em Appl. Comput. Harmon. Anal.}, 21(1):53--94, 2006.

\bibitem{HM_Dahmen}
W.~Dahmen.
\newblock Wavelet and multiscale methods for operator equations.
\newblock {\em Acta Numer.}, 6:55--228, 1997.

\bibitem{HM_DHS}
W.~Dahmen, H.~Harbrecht, and R.~Schneider.
\newblock {C}ompression techniques for boundary integral equations. {O}ptimal
  complexity estimates.
\newblock {\em SIAM J. Numer. Anal.}, 43(6):2251--2271, 2006.

\bibitem{HM_DKU}
W.~Dahmen, A.~Kunoth, and K.~Urban.
\newblock Biorthogonal spline wavelets on the interval -- stability and moment
  conditions.
\newblock {\em Appl. Comp. Harm. Anal.}, 6(2):132--196, 1999.

\bibitem{HM_DPS}
W.~Dahmen, S.~Pr{\"o}{\ss}dorf, and R.~Schneider.
\newblock Wavelet approximation methods for pseudodifferential equations {II}:
  Matrix compression and fast solution.
\newblock {\em Adv. Comput. Math.}, 1(3):259--335, 1993.

\bibitem{HM_STE}
W.~Dahmen and R.~Stevenson.
\newblock Element-by-element construction of wavelets satisfying stability and
  moment conditions.
\newblock {\em SIAM J. Numer. Anal.}, 37(1):319--352, 1999.

\bibitem{HM_Daubechies}
I.~Daubechies.
\newblock {\em Ten Lectures on Wavelets}.
\newblock Society of Industrial and Applied Mathematics, Philadelphia, 1992.

\bibitem{HM_ISTA}
I.~Daubechies, M.~Defrise, and C.~De~Mol.
\newblock An iterative thresholding algorithm for linear inverse problems with
  a sparsity constraint.
\newblock {\em {C}omm. {P}ure {A}ppl. {M}ath.}, 57:1413--1457, 2004.

\bibitem{HM_Dev98}
R.A. DeVore.
\newblock Nonlinear approximation.
\newblock {\em Acta Numer.}, 7:51--150, 1998.

\bibitem{HM_contourlets}
M.N. Do and M.~Vetterli.
\newblock The contourlet transform: An efficient directional multiresolution
  image representation.
\newblock {\em IEEE Trans. Image Proc.}, 14:2091--2106, 2005.

\bibitem{HM_Donoho}
D.L. Donoho.
\newblock Compressed sensing.
\newblock {\em IEEE Trans. Inf. Theory}, 52(4):1289--1306, 2006.

\bibitem{HM_Fasshauer2007}
G.E. Fasshauer.
\newblock {\em Meshfree Approximation Methods with {MATLAB}}.
\newblock World Scientific, River Edge, 2007.

\bibitem{HM_rauhut}
S.~Foucart and H.~Rauhut.
\newblock {\em A Mathematical Introduction to Compressive Sensing}.
\newblock Applied and Numerical Harmonic Analysis. Birkh\"auser, New York,
  2013.

\bibitem{HM_Geo73}
A.~George.
\newblock Nested dissection of a regular finite element mesh.
\newblock {\em SIAM J. Numer. Anal.}, 10(2):345--363, 1973.

\bibitem{HM_Gie01}
K.~Giebermann.
\newblock Multilevel approximation of boundary integral operators.
\newblock {\em Computing}, 67:183--207, 2001.

\bibitem{HM_LorenzGriesse}
R.~Griesse and D.A. Lorenz.
\newblock A semismooth {N}ewton method for {T}ikhonov functionals with sparsity
  constraints.
\newblock {\em Inverse Problems}, 24(3):035007, 2008.

\bibitem{HM_GCG05}
V.~Guigue, A.~Rakotomamonjy, and S.~Canu.
\newblock Kernel basis pursuit.
\newblock In J.~Gama, R.~Camacho, P.B. Brazdil, A.M. Jorge, and L.~Torgo,
  editors, {\em Machine Learning: ECML 2005}, pages 146--157,
  Berlin-Heidelberg, 2005. Springer.

\bibitem{HM_shearlets}
K.~Guo and D.~Labate.
\newblock Optimally sparse multidimensional representation using shearlets.
\newblock {\em SIAM J. Math. Anal.}, 39(1):298--318, 2007.

\bibitem{HM_HB02}
W.~Hackbusch and S.~B{\"o}rm.
\newblock $\mathcal{H}^2$-matrix approximation of integral operators by
  interpolation.
\newblock {\em Appl. Numer. Math.}, 43(1-2):129--143, 2002.

\bibitem{HM_Hale2008}
N.~Hale, N.J. Higham, and L.N. Trefethen.
\newblock Computing {${\bf A}^\alpha,\ \log({\bf A})$}, and related matrix
  functions by contour integrals.
\newblock {\em SIAM J. Numer. Anal.}, 46(5):2505--2523, 2008.

\bibitem{HM_HKS05}
H.~Harbrecht, U.~K{\"a}hler, and R.~Schneider.
\newblock Wavelet {G}alerkin {BEM} on unstructured meshes.
\newblock {\em Comput. Vis. Sci.}, 8(3-4):189--199, 2005.

\bibitem{HM_HM1}
H.~Harbrecht and M.~Multerer.
\newblock A fast direct solver for nonlocal operators in wavelet coordinates.
\newblock {\em J.\ Comput.\ Phys.}, 428:110056, 2021.

\bibitem{HM_HM2}
H.~Harbrecht and M.~Multerer.
\newblock Samplets: Construction and scattered data compression.
\newblock {\em J.\ Comput.\ Phys.}, 471:111616, 2022.

\bibitem{HM_HMQ24}
H.~Harbrecht, M.~Multerer, and J.~Quizi.
\newblock The dimension weighted fast multipole method for scattered data
  approximation.
\newblock {\em arXiv:2402.09531}, 2024.

\bibitem{HM_HMSS}
H.~Harbrecht, M.~Multerer, O.~Schenk, and C.~Schwab.
\newblock Multiresolution kernel matrix algebra.
\newblock {\em Numer. Math.}, 156(3):1085–1114, 2024.

\bibitem{HM_HS}
H.~Harbrecht and R.~Schneider.
\newblock Biorthogonal wavelet bases for the boundary element method.
\newblock {\em Math. Nachr.}, 269(1):167--188, 2004.

\bibitem{HM_ERA5}
H.~Hersbach, B.~Bell, P.~Berrisford, S.~Hirahara, A.~Hor{\'a}nyi,
  J.~Muñoz-Sabater, J.~Nicolas, C.~Peubey, R.~Radu, D.~Schepers, A.~Simmons,
  C.~Soci, S.~Abdalla, X.~Abellan, G.~Balsamo, P.~Bechtold, G.~Biavati,
  J.~Bidlot, M.~Bonavita, G.~De~Chiara, P.~Dahlgren, D.~Dee, M.~Diamantakis,
  R.~Dragani, J.~Flemming, R.~Forbes, M.~Fuentes, A.~Geer, L.~Haimberger,
  S.~Healy, R.J. Hogan, E.~Hólm, M.~Janisková, S.~Keeley, P.~Laloyaux,
  P.~Lopez, C.~Lupu, G.~Radnoti, P.~de~Rosnay, I.~Rozum, F.~Vamborg,
  S.~Villaume, and J.-N. Thépaut.
\newblock The {ERA5} global reanalysis.
\newblock {\em Q. J. R. Meteorol.}, 146(730):1999--2049, 2020.

\bibitem{HM_HSS08}
T.~Hofmann, B.~Sch{\"o}lkopf, and A.J. Smola.
\newblock Kernel methods in machine learning.
\newblock {\em Ann. Stat.}, 36(3):1171--1220, 2008.

\bibitem{HM_HorI}
L.~H{\"o}rmander.
\newblock {\em The Analysis of Linear Partial Differential Operators. {I}}.
\newblock Classics in Mathematics. Springer, Berlin, 2003.
\newblock Distribution theory and Fourier analysis, Reprint of the 1990
  edition.

\bibitem{HM_HorIII}
L.~H{\"o}rmander.
\newblock {\em The Analysis of Linear Partial Differential Operators. {III}}.
\newblock Classics in Mathematics. Springer, Berlin, 2007.
\newblock Pseudo-differential operators, Reprint of the 1994 edition.

\bibitem{HM_learning}
G.~James, D.~Witten, T.~Hastie, and R.~Tibshirani.
\newblock {\em An {I}ntroduction to {S}tatistical {L}earning}.
\newblock Springer Texts in Statistics. Springer, New York, 2013.

\bibitem{HM_Kae07}
U.~K\"ahler.
\newblock {\em \(\mathcal{H}^2\)-Wavelet Galerkin BEM and Its Application to
  the Radiosity Equation}.
\newblock Dissertation TU Chemnitz, Chemnitz, 2007.

\bibitem{HM_LP10}
G.~Leobacher and F.~Pillichshammer.
\newblock {\em Introduction to Quasi-{M}onte {C}arlo Integration and
  Applications}.
\newblock Springer International Publishing, Cham, 2010.

\bibitem{HM_selinv}
L.~Lin, C.~Yang, J.~C. Meza, J.~Lu, L.~Ying, and W.~E.
\newblock Selinv---an algorithm for selected inversion of a sparse symmetric
  matrix.
\newblock {\em ACM Trans. Math. Softw.}, 37(4):40, 2011.

\bibitem{HM_Lorenz}
D.A. Lorenz.
\newblock Convergence rates and source conditions for tikhonov regularization
  with sparsity constraints.
\newblock {\em J. Inverse Ill-Posed Probl.}, 16(5):463--478, 2008.

\bibitem{HM_FWT}
S.G. Mallat.
\newblock A theory for multiresolution signal decomposition: The wavelet
  representation.
\newblock {\em IEEE Trans. Pattern Anal. Mach. Intell.}, 2(7), 1989.

\bibitem{HM_Mallat}
S.G. Mallat.
\newblock {\em A Wavelet Tour of Signal Processing}.
\newblock Academic Press, San Diego, 1999.

\bibitem{HM_Mallat2016}
S.G. Mallat.
\newblock Understanding deep convolutional networks.
\newblock {\em Philos. Trans. R. Soc. A}, 374(2065):2015.0203, 2016.

\bibitem{HM_MZ}
S.G. Mallat and Z.~Zhang.
\newblock Matching pursuits with time-frequency dictionaries.
\newblock {\em IEEE Trans. Sign. Proc.}, 41(12):3397--3415, 1993.

\bibitem{HM_Mic84}
C.A. Micchelli.
\newblock Orthogonal projections are optimal algorithms.
\newblock {\em J. Approx. Theory}, 40(2):101--110, 1984.

\bibitem{HM_MR77}
C.A. Micchelli and T.J. Rivlin.
\newblock A survey of optimal recovery.
\newblock In C.A. Micchelli and T.J. Rivlin, editors, {\em Optimal Estimation
  in Approximation Theory}, pages 1--54, New York, 1977. Springer.

\bibitem{HM_MFSS17}
K.~Muandet, K.~Fukumizu, B.~Sriperumbudur, and B.~Sch{\"o}lkopf.
\newblock Kernel mean embedding of distributions: {A} review and beyond.
\newblock {\em Found. Trends Mach. Learn.}, 10(1-2):1--141, 2017.

\bibitem{HM_DS58}
J.T.~Schwartz N.J.~Dunford.
\newblock {\em Linear Operators. Part I: General Theory}.
\newblock Interscience Publishers, New York, 1958.

\bibitem{HM_RE11}
I.~Ram, M.~Elad, and I.~Cohen.
\newblock Generalized tree-based wavelet transform.
\newblock {\em IEEE Trans. Signal Process.}, 59(9):4199--4209, 2011.

\bibitem{HM_RT}
R.~Ramlau and G.~Teschke.
\newblock A {T}ikhonov-based projection iteration for nonlinear ill-posed
  problems with sparsity constraints.
\newblock {\em Numer. Math.}, 104(2):177–203, 2006.

\bibitem{HM_Rob74}
A.H. Robinson.
\newblock A new map projection: Its development and characteristics.
\newblock {\em International yearbook of cartography}, 14(1974):145--155, 1974.

\bibitem{HM_Rodino}
L.~Rodino.
\newblock {\em Linear Partial Differential Operators in {G}evrey Spaces}.
\newblock World Scientific Publishing Co. Inc., River Edge, NJ, 1993.

\bibitem{HM_Schaback2006}
R.~Schaback and H.~Wendland.
\newblock Kernel techniques: from machine learning to meshless methods.
\newblock {\em Acta Numer.}, 15:543--639, 2006.

\bibitem{HM_SCHN}
R.~Schneider.
\newblock {\em {M}ultiskalen- und {W}avelet-{M}atrixkompression:
  {A}nalysisbasierte {M}ethoden zur {L}{\"o}sung gro{\ss{}}er vollbesetzter
  {G}leichungssysteme}.
\newblock B.G.~Teubner, Stuttgart, 1998.

\bibitem{HM_Seeley}
R.T. Seeley.
\newblock Topics in pseudo-differential operators.
\newblock In {\em Pseudo-Differential Operators (C.I.M.E., Stresa, 1968)}, page
  167–305. Edizioni Cremonese, Rome, 1969.

\bibitem{HM_Sha48}
C.E. Shannon.
\newblock A mathematical theory of communication.
\newblock {\em Bell Syst. Tech. J.}, 27(3):379--423, 1948.

\bibitem{HM_Tao}
T.~Tao and E.J. Cand{\`e}s.
\newblock Near-optimal signal recovery from random projections: universal
  encoding strategies?
\newblock {\em IEEE Trans. Inf. Theory}, 52(12):5406–5425, 2006.

\bibitem{HM_TW03}
J.~Tausch and J.~White.
\newblock Multiscale bases for the sparse representation of boundary integral
  operators on complex geometry.
\newblock {\em SIAM J. Sci. Comput.}, 24(5):1610--1629, 2003.

\bibitem{HM_Taylor81}
M.E. Taylor.
\newblock {\em Pseudodifferential Operators}, volume~34 of {\em Princeton
  Mathematical Series}.
\newblock Princeton University Press, Princeton, NJ, 1981.

\bibitem{HM_tropp}
J.A. Tropp.
\newblock Greed is good: algorithmic results for sparse approximation.
\newblock {\em IEEE Trans. Inf. Theory}, 50(10):2231--2242, 2004.

\bibitem{HM_PS}
T.~von Petersdorff and C.~Schwab.
\newblock Fully discrete multiscale {G}alerkin {BEM}.
\newblock In W.~Dahmen, A.~Kurdila, and P.~Oswald, editors, {\em Multiscale
  wavelet methods for PDEs}, pages 287--346. Academic Press, San Diego, 1997.

\bibitem{HM_PSS97}
T.~von Petersdorff, C.~Schwab, and R.~Schneider.
\newblock Multiwavelets for second-kind integral equations.
\newblock {\em SIAM J. Numer. Anal.}, 34(6):2212--2227, 1997.

\bibitem{HM_Wendland2004}
H.~Wendland.
\newblock {\em Scattered Data Approximation}.
\newblock Cambridge University Press, Cambridge, 2004.

\bibitem{HM_Rasmussen2006}
C.K. Williams and C.E. Rasmussen.
\newblock {\em Gaussian Processes for Machine Learning}.
\newblock MIT Press, Cambridge, 2006.

\bibitem{HM_Williams1998}
C.K.I. Williams.
\newblock Prediction with {G}aussian processes. {F}rom linear regression to
  linear prediction and beyond.
\newblock In M.I. Jordan, editor, {\em Learning in Graphical Models}, volume~89
  of {\em NATO ASI Series (Series D: Behavioural and Social Sciences)}.
  Springer, Dordrecht, 1998.

\end{thebibliography}
\end{document}